\documentclass{ametsoc}

\journal{mwr}

\bibpunct{(}{)}{;}{a}{}{,}

\usepackage{amsmath,amssymb}
\usepackage{textcomp}
\usepackage{fancyhdr}
\usepackage{multicol}
\usepackage{color}
\usepackage{hyphenat}
\usepackage[normalem]{ulem}
\usepackage[utf8]{inputenc}
\usepackage{caption}

\usepackage{lineno}
\nolinenumbers
\usepackage{epsfig,subfig}
\usepackage{subfig}
\usepackage{IEEEtrantools}
\usepackage{epstopdf}
\usepackage{multirow}
\usepackage{url}
\usepackage{setspace}
\DeclareMathOperator*{\argmin}{arg\,min}
\newcommand{\nScale}{0.4}

\title{Ensemble-based kernel learning for a class of data assimilation problems with imperfect forward simulators}

\authors{Xiaodong Luo\correspondingauthor{Xiaodong Luo, Norwegian Research Centre (NORCE), Nyg\r{a}dsgaten 112, 5008 Bergen, Norway}}

\affiliation{Norwegian Research Centre (NORCE), Norway} 

\email{xluo@norceresearch.no}

\abstract{Simulator imperfection, often known as model error, is ubiquitous in practical data assimilation problems. Despite the enormous efforts dedicated to addressing this problem, properly handling simulator imperfection in data assimilation remains to be a challenging task. In this work, we propose an approach to dealing with simulator imperfection from a point of view of functional approximation that can be implemented through a certain machine learning method, such as kernel-based learning adopted in the current work. To this end, we start from considering a class of supervised learning problems, and then identify similarities between supervised learning and variational data assimilation. These similarities found the basis for us to develop an ensemble-based learning framework to tackle supervised learning problems, while achieving various advantages of ensemble-based methods over the variational ones. After establishing the ensemble-based learning framework, we proceed to investigate the integration of ensemble-based learning into an ensemble-based data assimilation framework to handle simulator imperfection. In the course of our investigations, we also develop a strategy to tackle the issue of multi-modality in supervised-learning problems, and transfer this strategy to data assimilation problems to help improve assimilation performance. For demonstration, we apply the ensemble-based learning framework and the integrated, ensemble-based data assimilation framework to a supervised learning problem and a data assimilation problem with an imperfect forward simulator, respectively. The experiment results indicate that both frameworks achieve good performance in relevant case studies, and that functional approximation through machine learning may serve as a viable way to account for simulator imperfection in data assimilation problems.}

\begin{document}

\maketitle

\section*{Introduction}\label{sec:introduction}
In recent years, the advent of big data era has led to surging interest in handling big data assimilation problems in data assimilation community \citep{Miyoshi2016big,luo2018eficient}. In reservoir engineering, using 4D seismic data for reservoir characterization through a certain history matching method constitutes a big data assimilation problem crucial to the industry. 4D seismic data contain spatially rich information of the hydrocarbon reservoir, while such information is often unavailable -- or at least, extremely difficult to extract -- from the conventional production data. Qualitative use of 4D seismic for reservoir monitoring \citep{lumley2001time} has now become a standard tool in the industry, yet quantitative utilization of 4D seismic for reservoir characterization, often under the name of \textit{4D seismic history matching} (4D SHM), still appears to be unmatured.

In the past few years, there have been a series of investigations \citep{lorentzen2017history,lorentzen2018ECMOR,luo2016sparse2d_spej,luo2018correlation_spej,luo2018eficient,Luo2018CorrLoc_Norne,luo2018ECMOR_towards,luo2016estimating_cg} inside the author's group at the International Research Institute of Stavanger (IRIS, now a part of Norwegian Research Centre, NORCE), which were dedicated to the research and development of an efficient workflow for 4D SHM through ensemble-based data assimilation \citep{Evensen2009}. In our investigations, we encountered a few major challenges, namely, \textit{big data}, \textit{uncertainty quantification} and \textit{imperfection} in forward seismic simulators \citep{luoIORNorway2018_bigData}. Driven by the needs to address these identified challenges in our studies, certain (multidisciplinary) methods, such as image processing \citep{lorentzen2017history,luo2016estimating_cg}, sparse data representation \citep{lorentzen2017history,luo2016sparse2d_spej,luo2018eficient}, adaptive localization \citep{luo2018correlation_spej,Luo2018CorrLoc_Norne,luo2018ECMOR_towards}, have been exploited or developed, and integrated into an ensemble-based SHM workflow, whose efficacy is now demonstrated in a full Norne field case study using real production and seismic data \citep{lorentzen2018ECMOR}. 

So far, our investigations have been mainly dedicated to tackling the issues of \textit{big data} and \textit{uncertainty quantification}, while leaving the issue of \textit{imperfection} largely untouched. As an attempt towards addressing this remaining challenge, in the current work, we propose a method that treats simulator imperfection from the perspective of functional approximation through machine learning, and investigate the integration of this method into an ensemble-based data assimilation framework. 

Imperfection in forward simulators (also known as model error) is a ubiquitous problem in geophysical data assimilation practices. Imperfection will arise when there are, e.g., unresolved fine-scale resolutions, missing or mis-specified physical processes, incorrect boundary conditions and so on, in the course of developing a physics-based forward simulator. In the context of practical SHM, for instance, one may expect that the rock physics model (RPM), as an essential part of the forward seismic simulator, is prone to imperfection, since the RPM is often built upon simplified assumptions on rock physics, and calibrated using core or well log data at a few locations. 

In the course of identifying and handling simulator imperfection during data assimilation, a challenge involving the combined effects of imperfection and uncertain model state and/or parameters will arise. For instance, when there are substantial gaps (residuals) between real and simulated observations, they may be attributed to simulator imperfection, or the inability of the assimilation algorithm to obtain globally optimal estimations of model state and/or parameters, or both. As a result, a prerequisite for addressing the issue of simulator imperfection would be to choose a method that helps untangle the gross effects of simulator imperfection and uncertain model state and/or parameters. 

Currently, a common practice in this regard is to add some (typically) additive stochastic term into the forward simulator, as a simple way to represent simulator imperfection (see, for example, \citealp{berry2017correcting,carrassi2010accounting,dee1995line,evensen2018accounting,griffith2000adjoint,howes2017accounting,oliver2018calibration,sakov2018iterative,sommer2018flexible}). For practical convenience, one may presume that the stochastic term follows a Gaussian distribution, so that the effect of simulator imperfection is taken into account by including the mean and covariance matrix of the stochastic term into the assimilation algorithm. There are a few simplifying assumptions, e.g., whiteness, stationarity, absence of bias and normality \citep{dee1995line}, involved in this way of treating simulator imperfection, which may not necessarily be valid in practice. There is also some recent work that aims to account for simulator imperfection from other perspectives. For instance, in \cite{kopke2018accounting}, the authors assume that there is an orthogonality between residuals due to simulator imperfection and those due to uncertain model state and/or parameters. Based on this assumption, local basis functions can be constructed and used to describe simulator imperfection. In practice, however, it is not clear yet to what extent the orthogonality assumption may be valid.

In the current work, we consider an approach that treats the modelling of simulator imperfection as a functional approximation problem, which can be solved using a certain machine learning method. To this end, we start from a supervised learning problem, in which one aims to optimize a certain function that maps a set of training inputs to a corresponding set of training outputs. We first show similarities between supervised learning and variational data assimilation. Motivated by this observation, we then proceed to develop a derivative-free, ensemble-based learning framework to tackle a class of supervised learning problems. In doing so, we are able to not only achieve all the benefits in using ensemble-based methods (which will be discussed later), but also facilitate the integration of the proposed imperfection-handling method into an ensemble-based data assimilation framework, which is presented after introducing the ensemble-based learning framework. 

For demonstration, we investigate the performance of the ensemble-based learning framework in a supervised learning problem. We identify a challenge which may arise when multi-modal training inputs are present in the learning process, and propose a strategy that helps overcome this problem. After that, we study a data assimilation problem with an imperfect forward simulator. Ensemble-based learning is then incorporated into an ensemble-based assimilation algorithm to tackle the data assimilation problem, while the insights and experience gained in the supervised learning problem are transferred to the data assimilation problem, helping improve the performance of data assimilation. Based on the results obtained in these two experiments, we conclude the current work with discussions and some thoughts of future work.                               

\section*{An ensemble-based kernel learning algorithm for a class of supervised learning problems}\label{sec:ensemble_supervised_learning}

\subsection*{Supervised learning as a variational data assimilation problem}\label{subsec:var_supervised_learning}
We consider a class of supervised learning problems (SLP), in which we are given a set of $N_s$ inputs, denoted by $X \equiv \{ x_i: x_i \in \mathbb{D}_x \subseteq \mathbb{R} \}_{i=1}^{N_s}$; and the corresponding set of $N_s$ outputs, denoted by $Y \equiv \{ y_i: y_i \in \mathbb{D}_y \subseteq  \mathbb{R} \}_{i=1}^{N_s}$, with $\mathbb{D}_x$ and $\mathbb{D}_y$ being the domains with respect to the inputs and outputs, respectively. Here, our objective is to learn a certain function $h: \mathbb{D}_x \rightarrow \mathbb{D}_y$, such that $h(x_i)$ match $y_i$ ($i = 1, 2, \dotsb, N_s$) to a good extent. Note that, in general, the outputs $y_i$ may be contaminated by certain noise.

To achieve the above objective, one can solve the SLP as a regularized empirical risk minimization (ERM) problem \citep{scholkopf2002learning}), as defined below    
\begin{linenomath*}  
	\begin{IEEEeqnarray}{r} \label{eq:functional_ERM}
		h^* = \argmin_{h} \, \dfrac{1}{N_s} \, \sum_{i=1}^{N_s} L(y_i - h(x_i)) + \gamma R(\Vert h \Vert) \, ,
	\end{IEEEeqnarray}
\end{linenomath*}   
where $L$ is a suitable loss function that measures the distance between $y_i$ and $h(x_i)$, $\gamma R(\Vert h \Vert)$ is a regularization term, with $\gamma$ being the regularization parameter, $\Vert h \Vert$ the norm of $h$ with respect to a certain metric space, and $R$ the regularization operator. From the perspective of inverse problem theory \citep{Engl2000-regularization}, the regularization term is typically introduced to prevent the estimated function $h^*$ from over-fitting the training data, as well as avoid potential numerical issues in the course of solving the minimization problem. 

Clearly, without imposing any constraint on the functional $h$, the regularized ERM problem in Eq. (\ref{eq:functional_ERM}) is intractable. In practice, it is customary to assume that $h$ belongs to a certain function space, and can be approximated through some parametric model, e.g., in the form of
\begin{linenomath*}  
	\begin{IEEEeqnarray}{r} \label{eq:functional_approximation}
		h(\bullet) \approx \hat{h}(\bullet; \boldsymbol{\theta}) \, ,
	\end{IEEEeqnarray}
\end{linenomath*}      
where $\boldsymbol{\theta}$ is a set of parameters in the (parametric) functional $\hat{h}$. Since $\hat{h}$ is parametrized by $\boldsymbol{\theta}$, replacing $h$ by $\hat{h}$, then the regularized ERM problem in Eq. (\ref{eq:functional_ERM}) becomes a parameter estimation problem, in the form of
\begin{linenomath*}  
	\begin{IEEEeqnarray}{r} \label{eq:parameter_estimation_ERM}
		\argmin_{\boldsymbol{\theta}} \, \dfrac{1}{N_s} \, \sum_{i=1}^{N_s} L(y_i - \hat{h}(x_i;\boldsymbol{\theta})) + \gamma R(\boldsymbol{\theta}) \, .
	\end{IEEEeqnarray}
\end{linenomath*}

In addition, let us define 
\begin{linenomath*}  
	\begin{IEEEeqnarray}{ll} \label{eq:vecterization}
		& \mathbf{Y}^s = \left[y_1, y_2,\dotsb, y_{N_s}\right]^T \, ; \\
		& \mathbf{X}^s = \left[x_1, x_2,\dotsb, x_{N_s}\right]^T \, ; \\
		& \hat{\mathbf{H}} \left( \boldsymbol{\theta}; \mathbf{X}^s\right) = \left[\hat{h}(x_1;\boldsymbol{\theta}),\hat{h}(x_2;\boldsymbol{\theta}),\dotsb, \hat{h}(x_{N_s};\boldsymbol{\theta})\right]^T \, ,
	\end{IEEEeqnarray}
\end{linenomath*}    
where in $\hat{\mathbf{H}}$ we place the argument $\boldsymbol{\theta}$ in front of $\mathbf{X}^s$ to emphasize that now $\boldsymbol{\theta}$ is the quantity in estimation, and we use a semicolon to separate quantities in estimation (i.e., $\boldsymbol{\theta}$) and those that are given (i.e., $\mathbf{X}^s$). Similar custom  will be adopted later for notational convenience. 

To facilitate the introduction to our idea, we first consider the situation in which the training inputs $x_i$ follow a certain unimodal distribution. In this case, we choose the functionals $L$ and $R$ in such a way that 
\begin{linenomath*}  
	\begin{IEEEeqnarray}{ll} \label{eq:L_and_R_form}
		& \sum_{i=1}^{N_s} L(y_i - \hat{h}(x_i;\boldsymbol{\theta})) = \left(\mathbf{Y}^s - \hat{\mathbf{H}} \left( \boldsymbol{\theta}; \mathbf{X}^s\right) \right)^T \mathbf{C}_y^{-1} \left(\mathbf{Y}^s - \hat{\mathbf{H}} \left( \boldsymbol{\theta}; \mathbf{X}^s\right) \right) \, ; \\
		& R(\boldsymbol{\theta}) = (\boldsymbol{\theta} - \boldsymbol{\theta}^b)^T \mathbf{C}_{\theta}^{-1} (\boldsymbol{\theta} - \boldsymbol{\theta}^b) \, ,
	\end{IEEEeqnarray}
\end{linenomath*}  
where $\mathbf{C}_y^{-1}$ and $\mathbf{C}_{\theta}^{-1}$ are some pre-chosen weight matrices associated with $L$ and $R$, respectively, and $\boldsymbol{\theta}^b$ stands for a (pre-chosen) initial guess (called ``background'' hereafter) of $\boldsymbol{\theta}$. Under these settings, the regularized ERM problem in Eq. (\ref{eq:parameter_estimation_ERM}) is equivalent to
\begin{linenomath*}  
	\begin{IEEEeqnarray}{r} \label{eq:vecterized_ERM}
		\argmin_{\boldsymbol{\theta}} \,  \left(\mathbf{Y}^s - \hat{\mathbf{H}} \left( \boldsymbol{\theta}; \mathbf{X}^s\right) \right)^T \mathbf{C}_y^{-1} \left(\mathbf{Y}^s - \hat{\mathbf{H}} \left( \boldsymbol{\theta}; \mathbf{X}^s\right) \right) + \gamma (\boldsymbol{\theta} - \boldsymbol{\theta}^b)^T \mathbf{C}_{\theta}^{-1} (\boldsymbol{\theta} - \boldsymbol{\theta}^b) \, .
	\end{IEEEeqnarray}
\end{linenomath*}
Comparing Eqs. (\ref{eq:parameter_estimation_ERM}) and (\ref{eq:vecterized_ERM}), the scalar factor $1/N_s$ is dropped in Eq. (\ref{eq:vecterized_ERM}), with its impact being absorbed into the regularization parameter $\gamma$. From a perspective of data assimilation, Eq. (\ref{eq:vecterized_ERM}) constitutes a conventional variational data assimilation (VAR-DA) problem, which can be solved through, e.g., optimal interpolation (OI) or three-dimensional variational (3D-VAR) method \citep{Kalnay-atmospheric}.

When the training inputs follow a multi-modal distribution, it may be necessary to cluster the training inputs into different groups (so that each group contains unimodal training inputs), and then estimate a set of parameters $\boldsymbol{\theta}$ for each group, using an estimation method developed for unimodal cases. In this sense, a parameter-estimation method developed for unimodal cases can serve as the building block of a method for multi-modal cases, similar to the work of \cite{Hoteit2012,Luo2008-spgsf1}. For this reason, in what follows, we focus on presenting an ensemble-based estimation method for unimodal cases. We will discuss how one can adapt the developed method to multi-modal cases, when we come to a concrete SLP problem with multi-modal training inputs.     

\subsection*{An ensemble-based approach to solving the supervised learning problem}\label{subsec:ensemble_supervised_learning}
In analogy to the advance of assimilation approaches from the conventional variational methods \citep{Kalnay-atmospheric} to the more recent, ensemble-based methods \citep{Evensen2009}, it is natural for us to develop a certain ensemble-based method to tackle the SLP. To this end, instead of solving Eq. (\ref{eq:vecterized_ERM}) to obtain a single set of estimated parameters, we aim to estimate an ensemble of such parameters. By doing so, we will obtain all the intrinsic benefits in using ensemble-based methods, which includes, for instance \citep{luo2016sparse2d_spej}, 
\begin{itemize}
	\item no need to develop a complicated and time-consuming adjoint system (``adjoint free'');
	\item the capacity to provide a means of uncertainty quantification for the estimated results (``uncertainty quantification'');
	\item the ability to handle large numbers of state and/or parameter variables (``algorithm scalability'');
	\item straightforward and fast implementation (``implementation convenience'').   
\end{itemize}

Employing this ``ensemblizing'' strategy, we reformulate the regularized ERM problem in Eq. (\ref{eq:vecterized_ERM}) as an minimum-average-cost (MAC) problem \citep{luo2015Iterative}, in terms of
\begin{linenomath*}  
	\begin{IEEEeqnarray}{r} \label{eq:ensemblize_ERM}
		\argmin_{\{\boldsymbol{\theta}_j\}_{j=1}^{N_e}} \, \dfrac{1}{N_e} \sum_{j=1}^{N_e} \left\{ \left(\mathbf{Y}^s - \hat{\mathbf{H}} \left( \boldsymbol{\theta}_j; \mathbf{X}^s\right) \right)^T \mathbf{C}_y^{-1} \left(\mathbf{Y}^s - \hat{\mathbf{H}} \left( \boldsymbol{\theta}_j; \mathbf{X^s}\right) \right) + \gamma (\boldsymbol{\theta}_j - \boldsymbol{\theta}_j^b)^T \mathbf{C}_{\theta}^{-1} (\boldsymbol{\theta}_j - \boldsymbol{\theta}_j^b) \right\} \, ,
	\end{IEEEeqnarray}
\end{linenomath*}
where $N_e$ is the size of the ensemble $\boldsymbol{\Theta} \equiv \{\boldsymbol{\theta}_j\}_{j=1}^{N_e}$ of parameters in estimation. Note that each ensemble member $\boldsymbol{\theta}_j$ has its own associated background $\boldsymbol{\theta}_j^b$. Typically, the initial values of $\boldsymbol{\theta}_j^b$ are generated at random, thus $\boldsymbol{\theta}_j^b \neq \boldsymbol{\theta}_k^b$ almost surely if $j \neq k$. As a result, solving the MAC problem in Eq. (\ref{eq:ensemblize_ERM}) would result in an ensemble $\boldsymbol{\Theta}^a \equiv \{\boldsymbol{\theta}_j^a\}_{j=1}^{N_e}$ (called \textit{analysis ensemble} hereafter) of diversified estimates, and this naturally leads to a way of conducting uncertainty quantification for the estimated results.

Following the convention in ensemble-based methods, we choose $\mathbf{C}_y$ to be the covariance matrix of the observation noise in the outputs $y_i$, and $\mathbf{C}_{\theta}$ to be the sample covariance matrix with respect to the \textit{background ensemble} $\boldsymbol{\Theta}^b \equiv \{\boldsymbol{\theta}_j^b\}_{j=1}^{N_e}$, in the sense that
\begin{linenomath*}  
	\begin{IEEEeqnarray}{l} \label{eq:bg_statistics}
		\mathbf{C}_{\theta} = \mathbf{S}_{\theta}^b \left( \mathbf{S}_{\theta}^b \right)^T \, ; \\
		\mathbf{S}_{\theta}^b \equiv \dfrac{1}{\sqrt{N_e -1}} \left[\boldsymbol{\theta}_1^b - \bar{\boldsymbol{\theta}}^b, \boldsymbol{\theta}_2^b - \bar{\boldsymbol{\theta}}^b, \dotsb, \boldsymbol{\theta}_{N_e}^b - \bar{\boldsymbol{\theta}}^b \right] \, ; \\
		\bar{\boldsymbol{\theta}}^b = \dfrac{1}{N_e} \sum_{j=1}^{N_e} \boldsymbol{\theta}_j^b \, .  		
	\end{IEEEeqnarray}
\end{linenomath*}

As shown in \cite{luo2015Iterative}, with a linearization-based approximation strategy, a solution to the MAC problem in Eq. (\ref{eq:ensemblize_ERM}) is given by
\begin{linenomath*}  
	\begin{IEEEeqnarray}{l} 
		\label{eq:mac_solution}	\boldsymbol{\theta}_j^a = \boldsymbol{\theta}_j^b + \mathbf{K} \left(\mathbf{Y}^s - \hat{\mathbf{H}} \left( \boldsymbol{\theta}_j^b; \mathbf{X^s}\right) \right)	, \, j = 1, 2, \dotsb, N_e ; \\
		\label{eq:kalman_gain}	\mathbf{K} \equiv \mathbf{S}_{\theta}^b (\mathbf{S}_{h}^b)^T \left(\mathbf{S}_{h}^b(\mathbf{S}_{h}^b)^T + \gamma \, \mathbf{C}_y \right)^{-1} \, ; \\
		\mathbf{S}_{h}^b \equiv \dfrac{1}{\sqrt{N_e -1}} \left[\hat{\mathbf{H}} \left( \boldsymbol{\theta}_1^b; \mathbf{X^s}\right) - \bar{y}_h^b, \hat{\mathbf{H}} \left( \boldsymbol{\theta}_2^b; \mathbf{X^s}\right) - \bar{y}_h^b, \dotsb, \hat{\mathbf{H}} \left( \boldsymbol{\theta}_{N_e}^b; \mathbf{X^s}\right) - \bar{y}_h^b \right] \, ; \\	
		\label{eq:mean_prediction}	\bar{y}_h^b \equiv \hat{\mathbf{H}} \left( \bar{\boldsymbol{\theta}}^b; \mathbf{X^s}\right) \, .
	\end{IEEEeqnarray}
\end{linenomath*} 
Eqs. (\ref{eq:bg_statistics}) through (\ref{eq:mean_prediction}) essentially constitute the iterative ensemble smoother (iES) used in \cite{luo2015Iterative}. In a practical implementation of the iES update formula Eq. (\ref{eq:mac_solution}), one may choose to apply a truncated singular value decomposition (TSVD) to $\mathbf{S}_{h}^b$, so that the matrix inversion in Eq. (\ref{eq:kalman_gain}) can be carried out in a low-dimensional subspace (with the dimension less than $N_e$). For more information, see \cite{chen2013-levenberg,Evensen2009,Luo2018CorrLoc_Norne}. 

In addition, the update formula Eq. (\ref{eq:mac_solution}) often has to be iterated for a number of times to make sure that the estimated parameters would be able to achieve good data match. In such an iteration process, we adopt a ``warm restart'' strategy, in such a way that an analysis ensemble at one iteration step serves as the background ensemble at the next iteration step. The regularization parameter $\gamma$ also needs to adapt to the iteration process, and is chosen in such a way to avoid either too big or too small iteration steps. Details on the choice of $\gamma$ and the associated stopping criteria are elaborated in \cite{luo2015Iterative,Luo2018CorrLoc_Norne}, and are skipped in this work for succinctness.      

\subsection*{RBF kernel based functional approximation}\label{subsec:kernel_supervised_learning}
After establishing an ensemble-based framework to handle SLP, we go back to discuss the concrete approach to functional approximation in Eq. (\ref{eq:functional_approximation}). In this regard, there are many methods (see, e.g. \citealp{murphy2012machine,goodfellow2016deep}), such as generalized linear models (GLM), support vector machines (SVM), and various (shallow or deep) neural networks that one may exploit. In the current work, taking into account various factors like capacity, complexity and cost, we choose to adopt radial-basis-function (RBF) based kernels for functional approximation. The RBF kernel approach was previously proposed in a seminal work of \cite{broomhead1988radial} to solve the SLP in Eq. (\ref{eq:functional_ERM}) in a way similar to a VAR-DA method, and this led to the establishment of RBF networks \citep{haykin2008neural}. Recently, the RBF kernel approach is also adopted by \cite{guo2017physics} to build computationally cheap surrogate models for history matching. 

Specifically, following the RBF kernel approach to functional approximation in \cite{broomhead1988radial}, we have
\begin{linenomath*}  
	\begin{IEEEeqnarray}{l} 
		\label{eq:rbf_approximation} \hat{h}(x; \boldsymbol{\theta}) = \sum_{k=1}^{N_{cp}} c_k \, K \left( x - x_k^{cp} ;\beta_k\right) \, ; \\
		\label{eq:rbf_kernel} K \left( x - x_k^{cp} ;\beta_k\right) \equiv \exp\{-\beta_k^2 \,  (x - x_k^{cp}) ^2 / 2 \} \, ; \\
		\label{eq:kernel_para} \boldsymbol{\theta} = [c_1,c_2,\dotsb,c_{N_{cp}} \vert \beta_1,\beta_2,\dotsb,\beta_{N_{cp}}]^T \, . 
	\end{IEEEeqnarray}
\end{linenomath*}      
Note that in Eq. (\ref{eq:rbf_kernel}), we adopt the Gaussian RBF kernel, but other types of kernel functions can also be used, as long as they serve the purpose of functional approximation well. Hereafter, for the convenience of discussion, we may drop the word(s) ``Gaussian'' and/or ``RBF''.

Eq. (\ref{eq:rbf_approximation}) indicates that, the approximation functional $\hat{h}$ is composed of a set of $N_{cp}$ kernels $K$, which are associated different weights $c_k$, center points (CP) $x_k^{cp}$, and scale parameters $\beta_k$ that influence the spreads of the kernels. In the current work, for simplicity, we pre-choose $N_{cp}$ and $x_k^{cp}$, such that $\hat{h}$ is parametrized by a set $\boldsymbol{\theta}$ of parameters $c_k$ and $\beta_k$,       
as indicated in Eq. (\ref{eq:kernel_para}) (for ease of visualization, we use ``$\vert$'' to separate different groups of parameters in Eq. (\ref{eq:kernel_para})). 

Eqs. (\ref{eq:rbf_approximation}) through (\ref{eq:kernel_para}) are for univariate problems. To extend the kernel approach to multivariate problems (e.g., $\mathbf{x} \in \mathbb{D}_x \subseteq \mathbb{R}^m$), one may consider the following form:   
\begin{linenomath*}  
	\begin{IEEEeqnarray}{l} 
		\label{eq:mt_rbf_approximation} \hat{h}(\mathbf{x}; \boldsymbol{\theta}) = \sum_{k=1}^{N_{cp}} c_k \, K \left( \mathbf{x} - \mathbf{x}_k^{cp} ;\boldsymbol{\beta}_k\right) \, ; \\
		\label{eq:mt_rbf_kernel} K \left( \mathbf{x} - \mathbf{x}_k^{cp} ;\boldsymbol{\beta}_k\right) \equiv \exp\{-\left< \boldsymbol{\beta}_k^2, (\mathbf{x} - \mathbf{x}_k^{cp})^2 \right> \} \, ; \\
		\label{eq:mt_inner_product} \left< \boldsymbol{\beta}_k^2, (\mathbf{x} - \mathbf{x}_k^{cp})^2 \right> \equiv \dfrac{1}{2m} \sum_{\ell = 1}^{m} \beta_{k,\ell}^2 \,  (x_{\ell} - x_{k,\ell}^{cp})^2 \, ; \\
		\label{eq:mt_kernel_para} \boldsymbol{\theta} = [c_1,c_2,\dotsb,c_{N_{cp}} \vert \beta_{1,1},\beta_{2,1},\dotsb,\beta_{N_{cp},1} \vert \dotsb \vert \beta_{1,m},\beta_{2,m},\dotsb,\beta_{N_{cp},m}]^T \, . 
	\end{IEEEeqnarray}
\end{linenomath*} 
In reservoir history matching problems, $m$ can be interpreted as the number of different types of petrophysical parameters (e.g., permeability, porosity and so on) associated with each reservoir gridblock, hence typically it may not be very large.

Eq. (\ref{eq:mt_inner_product}) considers generic anisotropic scale parameters $\boldsymbol{\beta}_k$ that may have different values $\beta_{k,\ell}$ along different axes $x_{\ell}$ ($\ell = 1, 2, \dotsb, m$). In addition, the factor $1/m$ in Eq. (\ref{eq:mt_inner_product}) is adopted to mitigate the issue of arithmetic underflow, which may arise in case that $\sum_{\ell = 1}^{m} \beta_{k,\ell}^2 \,  (x_{\ell} - x_{k,\ell}^{cp})^2$ becomes sufficiently large.
Under the above settings, the total number (cardinality) of parameters in $\boldsymbol{\theta}$ is thus $(m+1) \times N_{cp}$, as indicated in Eq. (\ref{eq:mt_kernel_para}). Therefore, the cardinality of $\boldsymbol{\theta}$ is controlled by the number $N_{cp}$ of center points, while the dimension $m$ of $\mathbf{x}$ is typically fixed. 

In comparison to the previous work of \cite{broomhead1988radial,guo2017physics}, one feature of our proposed Kernel approach to functional approximation is that the scale parameters  $\boldsymbol{\beta}_k$ not only adapt to different center points $\mathbf{x}_k^{cp}$, but also vary along different coordinate axes. This kind of flexibility may be considered desirable in the context of machine learning, as it leads to additional parameters that may help improve the expressive power (or capacity) of a learning model to match the training data, and also reduce generalization errors \citep{goodfellow2016deep}.    

In addition, in many existing publications, the scale parameters are often manually chosen. In contrast, the ensemble-based approach (Eqs. (\ref{eq:mac_solution}) through (\ref{eq:mean_prediction})) renders an efficient and derivate-free framework to estimate multiple sets of such parameters, hence also provides a natural means of uncertainty quantification for the estimation results. 


\section*{Kernel-based learning workflow for a class of data assimilation problems with imperfect forward simulators}\label{sec:SLP_to_imperfect_DA}
After establishing ensemble-based kernel learning to deal with SLP, we investigate how this framework can be integrated into ensemble-based data assimilation to handle a class of data assimilation problems with imperfect forward simulators, which bear certain similarities to 4D SHM problems. We will first formulate a mathematical description of the assimilation problems, and then develop a solution that combines ensemble-based approaches to both supervised learning and data assimilation. Within this integrated, ensemble-based framework, the solution to the target data assimilation problems involves a certain joint estimation procedure, in which one aims to simultaneously estimate both model variables and parameters associated with a set of kernel functions. From this perspective, technically speaking, this type of data assimilation problems with imperfect forward simulators would not become substantially more complicated than the corresponding assimilation problems with perfect forward simulators. Indeed, as will be shown later, for data assimilation in the presence of an imperfect forward simulator, one can still use existing ensemble-based assimilation algorithms, although there is a need to modify the forward simulator by including a residual functional to account for possible imperfection.

\subsection*{Problem statement}\label{subsec:DA_imperfect_ps}  
We consider a data assimilation problem, in which the noisy observational data (observations) $\mathbf{d}^{o} \in \mathbf{D}_{\mathbf{d}^o}$ are obtained through the following observation system
\begin{linenomath*}  
	\begin{IEEEeqnarray}{l} 
		\label{eq:DA_true_obs_sys} \mathbf{d}^{o} = \mathbf{f} \left(\mathbf{z}^{tr}\right) + \boldsymbol{\epsilon} \, , 
	\end{IEEEeqnarray}
\end{linenomath*}       
where $\mathbf{z}^{tr} \in \mathbf{D}_{\mathbf{z}^{tr}} \subseteq \mathbb{R}^{m_z}$ represents a set of true model variables, $\mathbf{f}: \mathbf{D}_{\mathbf{z}^{tr}}  \rightarrow \mathbf{D}_{\mathbf{d}^o} $ the true forward simulator, and $\boldsymbol{\epsilon}$ additive observation noise, which is assumed to follow a Gaussian distribution with zero mean and covariance matrix $\mathbf{C}_{\mathbf{d}}$. For better comprehension, here we have deliberately avoided notational overlapping with those in the proceeding section as far as possible, since such distinctions would be useful for our discussions later.

In the current work, we assume that, for all $\mathbf{z} \in \mathbb{R}^{m_z}$, one has $\mathbf{f}(\mathbf{z}) = [f(z_1),f(z_2),\dotsb,f(z_{m_z})]^T$, where $f: R \rightarrow R$ is a scalar function. Thus, an immediate implication of this assumption is that the size of observations is also equal to $m_z$, i.e., $\mathbf{d}^{o} = [d_1^o,d_2^o,\dotsb,d_{m_z}^o]^T$. We note that, the assumption we made here aims to mimic the situation in seismic history matching problems (or other similar geophysical inversion problems which involve spatially distributed, image-like geophysical data), but with certain simplifications to facilitate computations and discussions later. Under this setting, one may treat the scalar function $f$ as an analogy to a rock physics model, which maps petrophysical and/or dynamical parameters (the inputs) to certain seismic attributes (the outputs), such as acoustic impedance, distributed over reservoir gridblocks.   

As a data assimilation problem, our objective is to estimate a set $\mathbf{z}$ of model variables, conditioned on the observations $\mathbf{d}^{o}$ and some initial guess (background) of $\mathbf{z}^b$, in such a way that $\mathbf{z}$ is as ``close'' to $\mathbf{z}^{tr}$ as possible. In a typical setting, we have access to a certain forward simulator $\mathbf{g}$ that maps $\mathbf{z}$ to some simulated (or predicted) observations $\mathbf{d}^{sim}$, i.e.,
\begin{linenomath*}  
	\begin{IEEEeqnarray}{l} 
		\label{eq:DA_sim_obs_sys} \mathbf{d}^{sim} = \mathbf{g} \left(\mathbf{z}\right) \, , 
	\end{IEEEeqnarray}
\end{linenomath*}  
with $\mathbf{g}(\mathbf{z}) = [g(z_1),g(z_2),\dotsb,g(z_{m_z})]^T$ for a scalar function $g: R \rightarrow R$. This simulator, $\mathbf{g}$, is often imperfect, and may not be exactly identical to the true forward simulator $\mathbf{f}$. 
In the next subsection, we address the issue of imperfection by integrating ensemble-based kernel approach to functional approximation 
into an ensemble-based data assimilation framework.

\subsection*{Integrating ensemble-based kernel learning into data assimilation}\label{subsec:DA_imperfect_plus_RBF_kernel}  
Based on Eqs. (\ref{eq:DA_true_obs_sys}) and (\ref{eq:DA_sim_obs_sys}), we have
\begin{linenomath*}  
	\begin{IEEEeqnarray}{l} 
		\label{eq:imperfect_obs_sys_plus_residual} \mathbf{d}^{o} = \mathbf{g} \left(\mathbf{z}\right) + \mathbf{r}\left(\mathbf{z}; \mathbf{d}^{o} \right); \\
		\mathbf{r}\left(\mathbf{z}; \mathbf{d}^{o}\right) \equiv  \mathbf{d}^{o} - \mathbf{g} \left(\mathbf{z}\right), 
	\end{IEEEeqnarray}
\end{linenomath*}   
where $\mathbf{r}$ represents a functional of residuals that measure the differences between real observations $\mathbf{d}^{o}$ and the simulations $\mathbf{g} \left(\mathbf{z}\right)$. As in the preceding subsection, we have $\mathbf{r}\left(\mathbf{z}; \mathbf{d}^{o}\right) = [r(z_1;d_1^o),r(z_2;d_2^o),\dotsb,r(z_{m_z};d_{m_z}^o)]^T$, with $r(z_{\ell};d_{\ell}^o) = d_{\ell}^o - g(z_{\ell}^o)$ for $\ell = 1, 2, \dotsb, m_z$. 

Following the idea of kernel approach to functional approximation (Eqs. (\ref{eq:rbf_approximation}) through (\ref{eq:mt_kernel_para})), we can approximate $r(z_{\ell};d_{\ell}^o)$ by
\begin{linenomath*}  
	\begin{IEEEeqnarray}{l} 
		\label{eq:rbf_residual_approx} r(z_{\ell};d_{\ell}^o) \approx \hat{r}\left(z_{\ell}, \boldsymbol{\eta} \right) \equiv \hat{r}\left(z_{\ell}, \boldsymbol{\eta} ; d_{\ell}^o, \mathbf{D}^{o,cp}, \mathbf{Z}^{cp}\right) \, ,
	\end{IEEEeqnarray}
\end{linenomath*}   
where $\hat{r}$ is composed of a set of kernels with their parameters contained in $\boldsymbol{\eta}$, in the form of
\begin{linenomath*}  
	\begin{IEEEeqnarray}{l} 
		\label{eq:residual_rbf_approximation} \hat{r}\left(z_{\ell}, \boldsymbol{\eta} \right) = \sum_{k=1}^{N_{cp}} c_k \, \exp\left\{-\left< \boldsymbol{\beta}_k^2, \left(\begin{bmatrix}
			z_{\ell} \\
			d_{\ell}^o - g(z_{\ell})
		\end{bmatrix} - 
		\begin{bmatrix}
			z_k^{cp} \\
			d_{\ell}^o - d_k^{o,cp}
		\end{bmatrix}
		\right)^2  \right> \right\} \, , 
	\end{IEEEeqnarray}
\end{linenomath*} 
with the operator $\left<\bullet,\bullet\right>$ being defined in Eq. (\ref{eq:mt_inner_product}); $\mathbf{Z}^{cp} \equiv \{ z_k^{cp} \}_{k=1}^{N_{cp}}$ represents a set of center points $z_k^{cp}$, and $\mathbf{D}^{o,cp} \equiv \{ d_k^{o,cp} \}_{k=1}^{N_{cp}}$ stands for the corresponding set of observations associated with $\mathbf{Z}^{cp}$. Likewise, we define $\hat{\mathbf{r}}\left(\mathbf{z}; \boldsymbol{\eta}\right) \equiv [\hat{r}\left(z_1, \boldsymbol{\eta} \right),\hat{r}\left(z_2, \boldsymbol{\eta} \right),\dotsb,\hat{r}\left(z_{m_z}, \boldsymbol{\eta} \right)]^T$.

It is worth noting an essential difference between SLP and data assimilation problems. In SLP (cf. Eqs. (\ref{eq:vecterized_ERM} and (\ref{eq:ensemblize_ERM})), one has multiple ``matched'' input-output pairs, $\mathbf{X}^s$ and $\mathbf{Y}^s$, respectively, as the training data; In data assimilation problems, however, typically we only have access to a single realization of the outputs (observations) $\mathbf{d}^{o}$ at a given time instance and a given spatial location, whereas our purpose is to infer possible inputs $\mathbf{z}$ given $\mathbf{d}^{o}$. Often, due to the limited capacity of the assimilation algorithm, $\mathbf{d}^{o}$ and $\mathbf{z}$ do not constitute a ``matched'' pair, or in other words, $\mathbf{z}$ would typically not be identical to the true model variables $\mathbf{z}^{tr}$ that generate the observations $\mathbf{d}^{o}$. Because of this inconsistency and the sample frequency of observations (at a given time instance and a given spatial location), data assimilation problems with imperfect forward simulators tend to be more challenging than SLP, as we will see later.

The aforementioned difference between SLP and data assimilation motivates us to take a slightly different form in Eq. (\ref{eq:residual_rbf_approximation}) for kernel-based functional approximation, in comparison to those in SLP (Eqs. (\ref{eq:rbf_approximation}) through (\ref{eq:mt_kernel_para})). Specifically, in Eq. (\ref{eq:residual_rbf_approximation}), we choose to augment both the model variables $z_{\ell}$ and the corresponding residuals $d_{\ell}^o - g(z_{\ell})$, and use the augmented vectors as the inputs to the kernel functions. In comparison to the settings in SLP, using $d_{\ell}^o - g(z_{\ell})$ in kernel functions allows us to tune additional scale parameters in data assimilation, which may be desirable in terms of flexibility. On the other hand, though, this also requires us to specify a set of observations $\mathbf{D}^{o,cp}$ associated with $\mathbf{Z}^{cp}$. In general, the choice of $\mathbf{Z}^{cp}$ and $\mathbf{D}^{o,cp}$ may be case-dependent. For instance, if one has a set of $(z_k^{cp},d_k^{o,cp})$ pairs from the hard data (e.g., those obtained from core analysis or well log data), then they can be included. In a case study later, we will give a specific implementation example on the choices of $\mathbf{Z}^{cp}$ and $\mathbf{D}^{o,cp}$. 

With kernel-based functional approximation to the residuals, similar to \cite{luo2015Iterative} (also see Eq. (\ref{eq:ensemblize_ERM})), the data assimilation problem with an imperfect forward simulator can then be addressed by solving the following optimization problem:
\begin{linenomath*}  
	\begin{IEEEeqnarray}{l} \label{eq:HM_with_SLP}
		\argmin_{\{\tilde{\boldsymbol{\theta}}_j\}_{j=1}^{N_e}} \, \dfrac{1}{N_e} \sum_{j=1}^{N_e} \left\{ \left(\mathbf{d}^{o} - \tilde{\mathbf{g}} \left(\tilde{\boldsymbol{\theta}}_j\right) \right)^T \mathbf{C}_{\mathbf{d}}^{-1} \left(\mathbf{d}^{o} - \tilde{\mathbf{g}} \left(\tilde{\boldsymbol{\theta}}_j\right) \right) + \gamma (\tilde{\boldsymbol{\theta}}_j - \tilde{\boldsymbol{\theta}}_j^b)^T \mathbf{C}_{\tilde{\theta}}^{-1} (\tilde{\boldsymbol{\theta}}_j - \tilde{\boldsymbol{\theta}}_j^b) \right\} \, , 
	\end{IEEEeqnarray}
\end{linenomath*}
with 
\begin{linenomath*}  
	\begin{IEEEeqnarray}{l} 
		\label{eq:HM_joint_vector} \tilde{\boldsymbol{\theta}} \equiv \left[ \mathbf{z}^T \, , \boldsymbol{\eta}^T \right]^T \, ; \\
		\label{eq:new_forward_simulator} \tilde{\mathbf{g}} \left(\tilde{\boldsymbol{\theta}}\right) \equiv \mathbf{g} \left(\mathbf{z}\right) + \hat{\mathbf{r}}\left(\mathbf{z}, \boldsymbol{\eta} \right) \, ,
	\end{IEEEeqnarray}
\end{linenomath*}         
where $\tilde{\boldsymbol{\theta}}$ is a joint vector that augments model variables $\mathbf{z}$ and parameters $\boldsymbol{\eta}$ associated with the set of kernels; $\mathbf{C}_{\tilde{\theta}}$ is the sample error covariance matrix with respect to an ensemble $\tilde{\mathbf{\Theta}}^b \equiv \left\{ \tilde{\boldsymbol{\theta}}_j^b \right\}_{j=1}^{N_e}$, similar to that in Eq. (\ref{eq:bg_statistics}); and $\tilde{\mathbf{g}} \left(\tilde{\boldsymbol{\theta}}\right)$ corresponds to the effective forward simulator. 

As in Eq. (\ref{eq:ensemblize_ERM}), Eq. (\ref{eq:HM_with_SLP}) also constitutes an MAC problem. As a result, Eqs. (\ref{eq:bg_statistics}) through (\ref{eq:mean_prediction}) provide an approximate solution to the data assimilation problem with an imperfect forward simulator, provided that one replaces $\hat{\mathbf{H}} \left( \boldsymbol{\theta}_j; \mathbf{X^s}\right)$, $\boldsymbol{\theta}$ and $\mathbf{C}_y$ therein by $\tilde{\mathbf{g}} \left(\tilde{\boldsymbol{\theta}}\right)$, $\tilde{\boldsymbol{\theta}}$ and $\mathbf{C}_{\mathbf{d}}$, respectively.

When there is no imperfection in the forward simulator (or when one believes so), one may choose not to introduce any correction mechanism. In this case, the parameter part $\boldsymbol{\eta}$ of $\tilde{\boldsymbol{\theta}}$ (cf. Eq. (\ref{eq:HM_joint_vector})) can be simply taken out.  Based on this observation, it is clear that adopting ensemble-based kernel approach to accounting for imperfection in the forward simulator does not significantly change our ensemble-based data assimilation algorithm. Instead, with a modified forward simulator $\tilde{\mathbf{g}} \left(\tilde{\boldsymbol{\theta}}\right)$ in Eq. (\ref{eq:new_forward_simulator}), it only requires some minor changes of the algorithm, by inserting a residual term into the original forward simulator, and then combining parameters associated with the kernel functions and the original model variables to form augmented vectors in data assimilation.. 

As will be shown later, even with a perfect forward simulator, it might be still beneficial to include a mechanism of model-error correction (i.e., the $\boldsymbol{\eta}$ term) for the improvement of data assimilation performance. The rationale behind this notion is that, similar to machine learning problems, the presence of $\boldsymbol{\eta}$ increases the dimension of $\tilde{\boldsymbol{\theta}}$, so that the assimilation algorithm would have more degrees of freedom to exploit for the search of better results.      


\section*{Numerical results in a supervised learning problem}\label{sec:results_SLP}
In this section, we investigate the performance of ensemble-based kernel learning in a toy supervised learning problem. One of our focuses here is to demonstrate a challenge arising in the toy problem, and develop a strategy that helps overcome this challenge. The insights obtained in the study will shed light on certain limitations or cautions in using the plain ensemble-based kernel learning framework, and the way for performance improvements. In turn, they will help enhance the data assimilation performance when integrating ensemble-based kernel learning into ensemble-based data assimilation.  

The supervised learning problem is designed to mimic the situation of data assimilation with an imperfect forward simulator. Specifically, we consider a forward  system
\begin{linenomath*}  
	\begin{IEEEeqnarray}{l} \label{eq:results_SLP_obs_sys}
		y^{o} = f(x) + \epsilon \, ; \\
		f(x) = \left(\vert x \vert^3 + 1 \right)^{1/2} \, ,
	\end{IEEEeqnarray}
\end{linenomath*}        
where $x \in \mathbb{R}$ is a scalar input, $y^{o} \in \mathbb{R}$ is the noisy output contaminated by Gaussian noise $\epsilon$, $f: \mathbb{R} \rightarrow \mathbb{R}$ represents the true mapping function, and  $\epsilon$ has zero mean, but its standard deviation (STD) $\sigma$ in general may depend on $f(x)$, in the form of $\sigma = \max(10^{-6},0.1 \times \vert f(x) \vert)$.

\begin{figure} 
	\centering
	\includegraphics[width=0.8\textwidth]{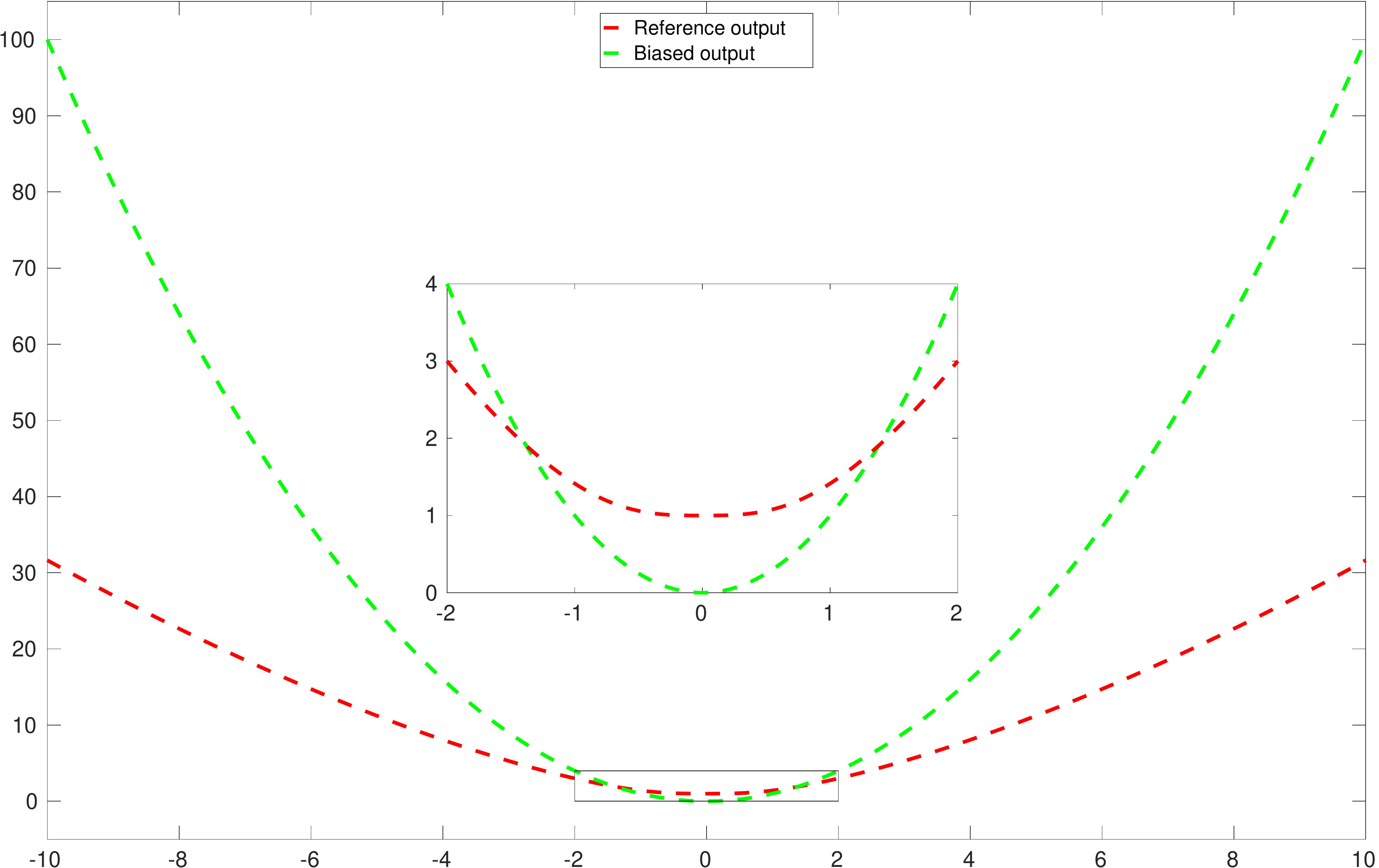}
	\caption{\label{fig:ref_bias_func} Outputs from the true (or reference) function (in red). and those from the biased one (in green). For better visualization, we re-plot the reference and biased outputs over the input interval $\left[-2, \, 2\right]$ in a separate, zoomed-in subplot.}
\end{figure}

In addition, we assume that there exists another imperfect forward simulation system
\begin{linenomath*}  
	\begin{IEEEeqnarray}{l} \label{eq:results_SLP_simObs_sys}
		y^{sim} = g(x) \, ; \\
		g(x) = x^2 \, .
	\end{IEEEeqnarray}
\end{linenomath*} 
In Figure \ref{fig:ref_bias_func}, we show the outputs of $f$ (without noise) and $g$, respectively, over the input interval $\left[-10,\, 10\right]$, while $f$ and $g$ intersect each other at $x \approx \pm 1.38$. Note that in the evaluations here (and also later), the relevant (e.g., reference, biased or prediction) functions are evaluated at the points from the set $\{-10:0.1:10\}$, which, following the MATLAB$^\copyright$ custom, represents a set of points evenly distributed over the interval $[-10, \, 10]$ with a span of $0.1$. For better visualization, we also re-plot their outputs over the input interval $\left[-2, \, 2\right]$ in a separate, zoomed-in subplot.
\renewcommand{\nScale}{0.3}

In the SLP, our objective is to learn the residual function $r(x) = f(x) - g(x)$ based on a certain set of training data. To this end, the ensemble-based approach developed in the previous section (Eqs. (\ref{eq:ensemblize_ERM}) through (\ref{eq:mean_prediction})) is adopted. In the experiments, we start with training data in which the (noisy) outputs are generated by some unimodal inputs, and then move to the more complicated situation in which the (noisy) outputs are produced using multi-modal inputs instead.  

For the purpose of comparison, in this work we adopt data mismatch and root mean squared error (RMSE) as performance measures. Following the notations in the previous sections, given the real observations $\mathbf{d}^o$, its associated observation error covariance matrix $\mathbf{C}_{d}$ and an ensemble member $\boldsymbol{\theta}_{j}$ (or $\tilde{\boldsymbol{\theta}}_j$) in SLP (or data assimilation problems), suppose that the simulated observations with respect to $\boldsymbol{\theta}_{j}$ (or $\tilde{\boldsymbol{\theta}}_j$) are $\mathbf{d}_{j}^{sim}$, then the corresponding data mismatch $\Xi_j$ is defined as  
\begin{linenomath*}    
	\begin{IEEEeqnarray}{l} \label{eq:data_mismatch}
		\Xi_j \equiv  \left( \mathbf{d}^o - \mathbf{d}_{j}^{sim} \right)^T \mathbf{C}_{d}^{-1} \left( \mathbf{d}^o - \mathbf{d}_{j}^{sim} \right), \,  j = 1, 2, \dotsb, N_e \, ,
	\end{IEEEeqnarray}
\end{linenomath*}                       
while the RMSE $e_j$ of the model $\mathbf{z}_j$ (in data assimilation problems) with respect the reference model $\mathbf{z}^{tr}$ is 
\begin{linenomath*} 
	\begin{IEEEeqnarray}{lll} \label{eq:RMSE_def}
		e_j = \dfrac{\Vert \mathbf{z}_j - \mathbf{z}^{tr} \Vert_2}{\sqrt{m_z}} \, .
	\end{IEEEeqnarray}
\end{linenomath*}

Throughout this work, we use the iES in \cite{luo2015Iterative} as the ensemble-based learning (or data assimilation) algorithm to update the relevant parameters, although in principle other iES, e.g., \cite{chen2013-levenberg,emerick2012ensemble}, may also be adopted. In the experiments, the configuration of the iES is as follows. The maximum (outer) iteration step is set to 10. If an iteration successfully reduces the average data mismatch (over the ensemble members), then the current value of the regularization parameter $\gamma$ is multiplied by a factor of $2$, aiming to further increase the step size of the next iteration. In this case, the analysis ensemble at the current iteration step will be used as the background one at the next iteration. In contrast, if the iteration leads to higher average data mismatch, then following \cite{chen2013-levenberg}, we start a trial (inner) iteration process, in which the background ensemble at the current (outer) iteration step is always used as the background ensemble in the trial process. A back-track line search strategy is adopted, in such a way that the current value of $\gamma$ is multiplied by a factor of $0.9$, and then used in a trial iteration to see if the new average data mismatch becomes lower than the original average at the current outer iteration step. The trial iteration is repeated maximum 5 times, but an earlier stop may take place if lower average data mismatch is found at a certain trial iteration step. We then use the last analysis ensemble obtained from the trial process as the background 
at the next outer iteration step. Apart from the maximum number of (outer) iteration steps, we also adopt another two stopping criteria, which become effective if (1) the change of average data mismatch values in two consecutive iterations are less than $1\%$ (for runtime control); or if (2) the average data mismatch is lower than four times the number of observations for the first time (to avoid over-fitting observations, see \citealp{luo2016sparse2d_spej}). 
For ease of comparison, localization \citep{chen2010cross,Emerick2011combining,luo2018correlation_spej,luo2018ECMOR_towards} is not adopted in the iES.  

\subsection*{Results with respect to unimodal inputs}
\begin{figure} 
	\centering
	\includegraphics[width=0.8\textwidth]{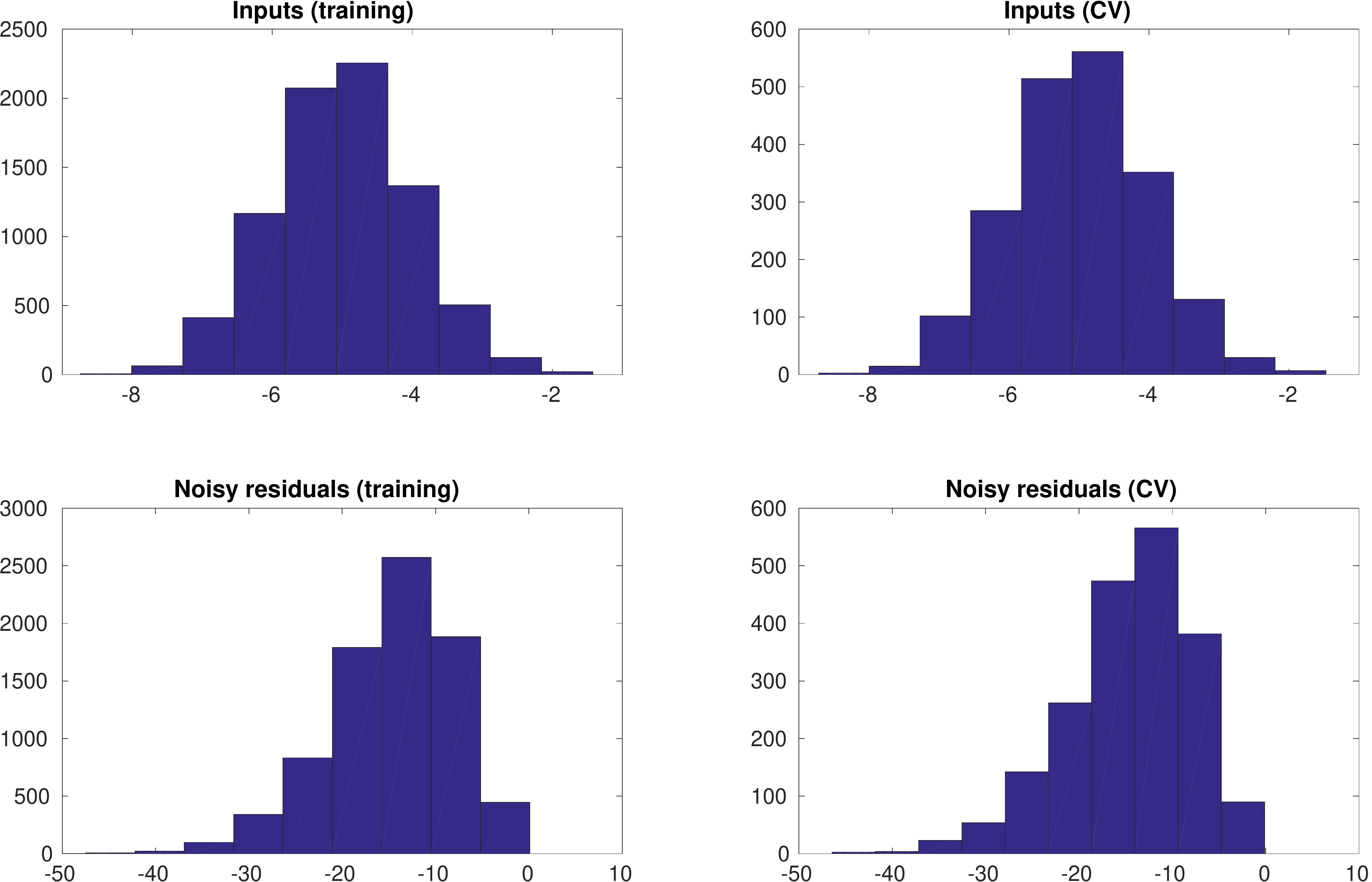}
	\caption{\label{fig:hist_unimodal} Histograms of the unimodal inputs and noisy residuals (as labels), with respect to the training and cross-validation (CV) datasets, respectively.}
\end{figure}

In the experiment, we generate a set of $10,000$ input samples drawn from the univariate Gaussian distribution $N(-5,1)$, and the corresponding noisy residuals (defined as the differences between the noisy outputs $y^{o}$ and the simulations $y^{sim}$). We randomly divide the set of input-residual pairs into two subsets: one with $8,000$ ($80\%$) of such pairs as the training dataset, whereas the rest $2,000$ ($20\%$) of such pairs as the cross-validation (CV) dataset. The training dataset is used to estimate the parameters associated with the selected kernel functions, whereas the CV dataset is not involved in learning these parameters. In a typical setting, the CV dataset can be adopted to select hyperparameter(s) in a learning algorithm. In this particular case, though, we do not have hyperparameter(s) to tune. Therefore, we simply use the CV dataset to inspect the performance of the learned parameters after the learning process is finished. Figure \ref{fig:hist_unimodal} shows the histograms of the inputs and noisy residuals in the training and CV datasets.          
\renewcommand{\nScale}{0.3}

To employ the kernel approach to approximating the residual function $r(x)$ (Eqs. (\ref{eq:rbf_approximation}) through (\ref{eq:kernel_para})), we need to specify a number of $N_{cp}$ center points $x_{k}^{cp}$ ($k = 1, 2, \dotsb, N_{cp}$). In principle, it is possible to consider both $N_{cp}$ and $x_{k}^{cp}$ as additional parameters that may be optimized through certain criteria. However, this will make the resulting learning algorithm become much more complicated. As a result, in the current work, we pre-choose $N_{cp}$ and $x_{k}^{cp}$ manually. 

Bearing this in mind, in the experiment below, we let $N_{cp} = 200$, and $x_{k}^{cp}$ be the points that evenly span the half-closed interval $\left[-6, 6\right)$. We have also tested other cases with $N_{cp} = 2000$, which turned out to lead to results similar to what we will present below. Consequently, for brevity, below we focus on the cases with $N_{cp} = 200$, with which the number of parameters (including the weights $c_k$ and the scale parameters $\beta_k$, cf. Eq. (\ref{eq:kernel_para})) is thus $2 \times 200 = 400$. 

An additional remark is that, in comparison to the histograms in Figure \ref{fig:hist_unimodal}, it is clear that the interval $\left[-6, 6\right)$ and the input ranges of both training and CV datasets do not fully cover each other. We choose such a setting to examine the impact of data coverage on the performance of the learning algorithm. 

Now we discuss how to initialize the ensembles of kernel parameters, weights $c_k$ and scale parameters $\beta_k$. For convenience of discussion, let us denote the ensembles with respect to the initial weights and the initial scale parameters by $\mathbf{C}^0 \equiv \left\{ c_{k,j}^0 \right\}_{j=1}^{N_e}$ and $\mathbf{B}^0 \equiv \left\{ \beta_{k,j}^0 \right\}_{j=1}^{N_e}$, respectively. In the current work, we let $N_e = 100$ unless otherwise stated, and $\beta_{k,j}^0$ be initialized
as follows:
\begin{linenomath*}  
	\begin{IEEEeqnarray}{l} \label{eq:init_beta}
		\beta_{k,j}^0 = \dfrac{1}{\sigma^{ti}} \times \exp(\xi_{k,j}), \text{ for } k = 1, 2, \dotsb, N_{cp}; \, j = 1, 2, \dotsb, N_e  ; \\
		\xi_{k,j} \sim N(0,1) \, ,
	\end{IEEEeqnarray}
\end{linenomath*}                 
where $\sigma^{ti}$ is the STD of the training inputs, and $\xi_{k,j}$ are random samples drawn from the normal distribution $N(0,1)$ for each center point and each ensemble member. 

For a given ensemble member (i.e., a fixed $j$ value), we then initialize $c_{k,j}^0$ ($k = 1, 2, \dotsb, N_{cp}$) as follows. We first randomly select a pair of input-label from the training dataset, denoted by $(x_{j}^{ti},\delta y_{j}^{tl})$, where $\delta y_{j}^{tl}$ is the label (noisy residual) in the training dataset that corresponds to the training input $x_{j}^{ti}$. We then insert the pair $(x_{j}^{ti},\delta y_{j}^{tl})$ into Eq. (\ref{eq:rbf_approximation}), by replacing $x$, $\hat{h}(x; \boldsymbol{\theta})$ and $\beta_k$ therein by $x_{j}^{ti}$, $\delta y_{j}^{tl}$ and $\beta_{k,j}^0$, respectively. At this stage, our goal is to find a set of weights $c_{k,j}^0$ ($k = 1, 2, \dotsb, N_{cp}$) that approximately solve the the following equation:
\begin{linenomath*}  
	\begin{IEEEeqnarray}{l} 
		\label{eq:init_c_sum_eq} \delta y_{j}^{tl} = \sum_{k=1}^{N_{cp}} c_{k,j}^0 \, K \left( x_{j}^{ti} - x_k^{cp} ; \beta_{k,j}^0\right) \, , 
	\end{IEEEeqnarray}
\end{linenomath*}         
which can be re-written as the following vector-based equation
\begin{linenomath*}  
	\begin{IEEEeqnarray}{l} 
		\label{eq:init_c_vec_eq} \delta y_{j}^{tl} = \left( \mathbf{c}_{j}^0 \right)^T \, \mathbf{K}(x_{j}^{ti}) \, ; \\
		\mathbf{c}_{j}^0 \equiv \left[c_{1,j}^0, c_{2,j}^0, \dotsb, c_{N_{cp},j}^0 \right]^T \, ; \\
		\mathbf{K}(x_{j}^{ti}) \equiv \left[ K \left( x_{j}^{ti} - x_1^{cp} ; \beta_{1,j}^0\right) , K \left( x_{j}^{ti} - x_2^{cp} ; \beta_{2,j}^0\right) , \dotsb, K \left( x_{j}^{ti} - x_{N_{cp}}^{cp} ; \beta_{{N_{cp}},j}^0\right)  \right]^T \, .
	\end{IEEEeqnarray}
\end{linenomath*}    
An approximate solution to Eq. (\ref{eq:init_c_vec_eq}) can be obtained by solving the following equation instead
\begin{linenomath*}  
	\begin{IEEEeqnarray}{l} 
		\label{eq:init_c_approx_eq} \delta y_{j}^{tl} \mathbf{K}(x_{j}^{ti}) = \left( \mathbf{K}(x_{j}^{ti}) \mathbf{K}(x_{j}^{ti})^T  + \alpha \mathbf{I}\right)  \mathbf{c}_{j}^0  \,  , 
	\end{IEEEeqnarray}
\end{linenomath*}  
where $\alpha$ is a positive scalar, and $\mathbf{I}$ is an $N_{cp} \times N_{cp}$ identity matrix. The term $\alpha \mathbf{I}$, essentially stemming from a Tikonov regularization term introduced to solve Eq. (\ref{eq:init_c_vec_eq}) as a regularized inverse problem \citep{Engl2000-regularization}, helps to improve the numerical stability of the final solution
\begin{linenomath*}  
	\begin{IEEEeqnarray}{l} 
		\label{eq:init_c_final_sol} \mathbf{c}_{j}^0 = \delta y_{j}^{tl} \mathbf{K}(x_{j}^{ti}) / \left(\alpha_{j} + \mathbf{K}(x_{j}^{ti})^T \mathbf{K}(x_{j}^{ti})\right) \, .
	\end{IEEEeqnarray}
\end{linenomath*}  
Following the implementation of the iES in \cite{luo2015Iterative}, in the current work, we let $\alpha_{j} = \exp(\xi_{j}) \, \mathbf{K}(x_{j}^{ti})^T \mathbf{K}(x_{j}^{ti})$ with $\xi_{j} \sim N(0,1)$. It is clear that the solution in Eq. (\ref{eq:init_c_final_sol}) does not solve Eq. (\ref{eq:init_c_vec_eq}) exactly. This, however, is desired, since in general the label $\delta y_{j}^{tl}$ may be noisy, and an inexact solution to Eq. (\ref{eq:init_c_vec_eq}) avoids the problem of over-fitting the training data. Applying Eq. (\ref{eq:init_c_final_sol}) to $N_e$ different pairs of $(x_{j}^{ti},\delta y_{j}^{tl})$ ($j = 1, 2, \dotsb, N_e$), we get an initial ensemble of $N_e$ different parameter vectors $\mathbf{c}_{j}^0$.  

\renewcommand{\nScale}{0.4}
\begin{figure} 
	\centering
	\begin{tabular}{cc}
		\subfloat[Training dataset]{\includegraphics[width=0.45\textwidth]{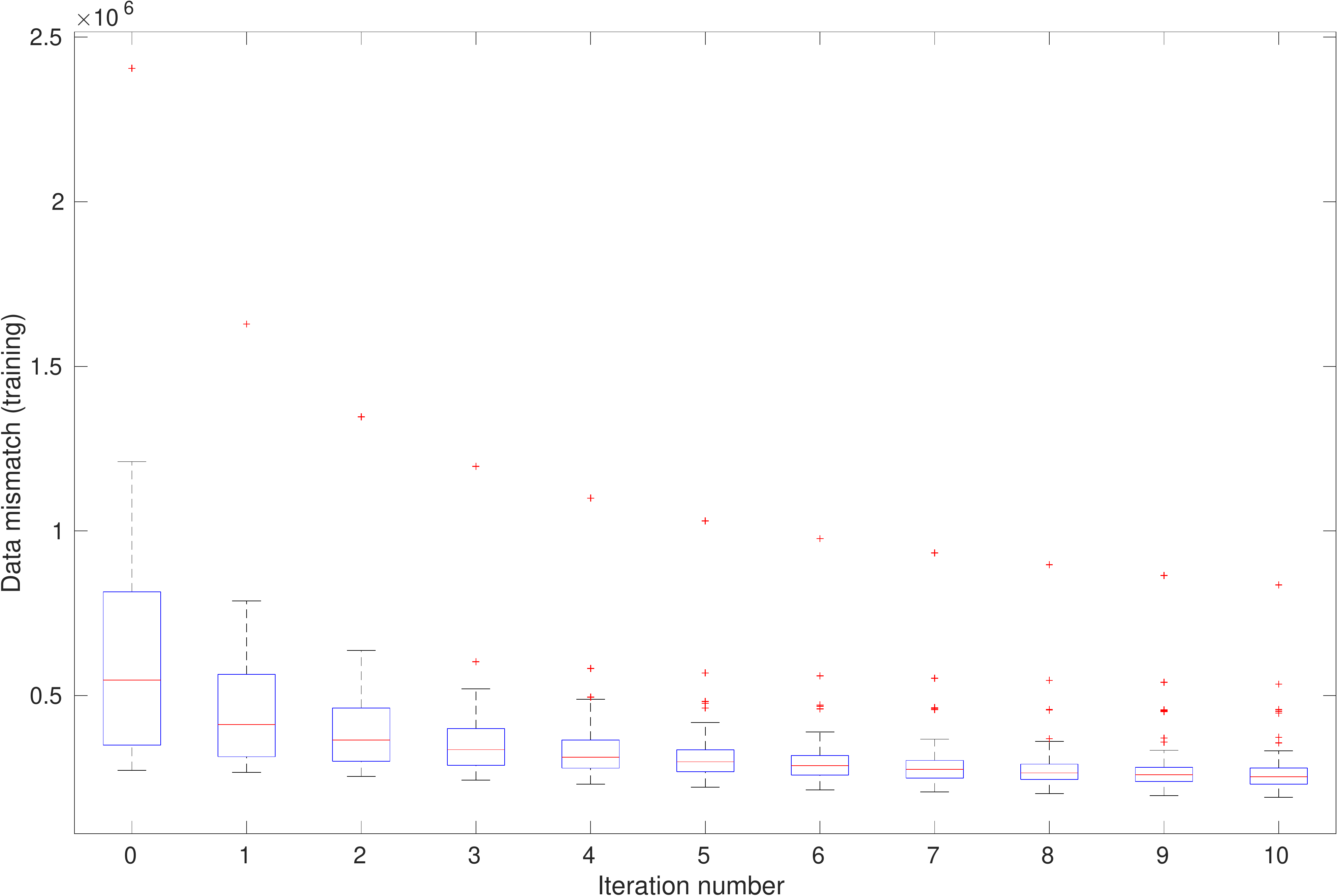}} & 
		\subfloat[CV dataset]{\includegraphics[width=0.45\textwidth]{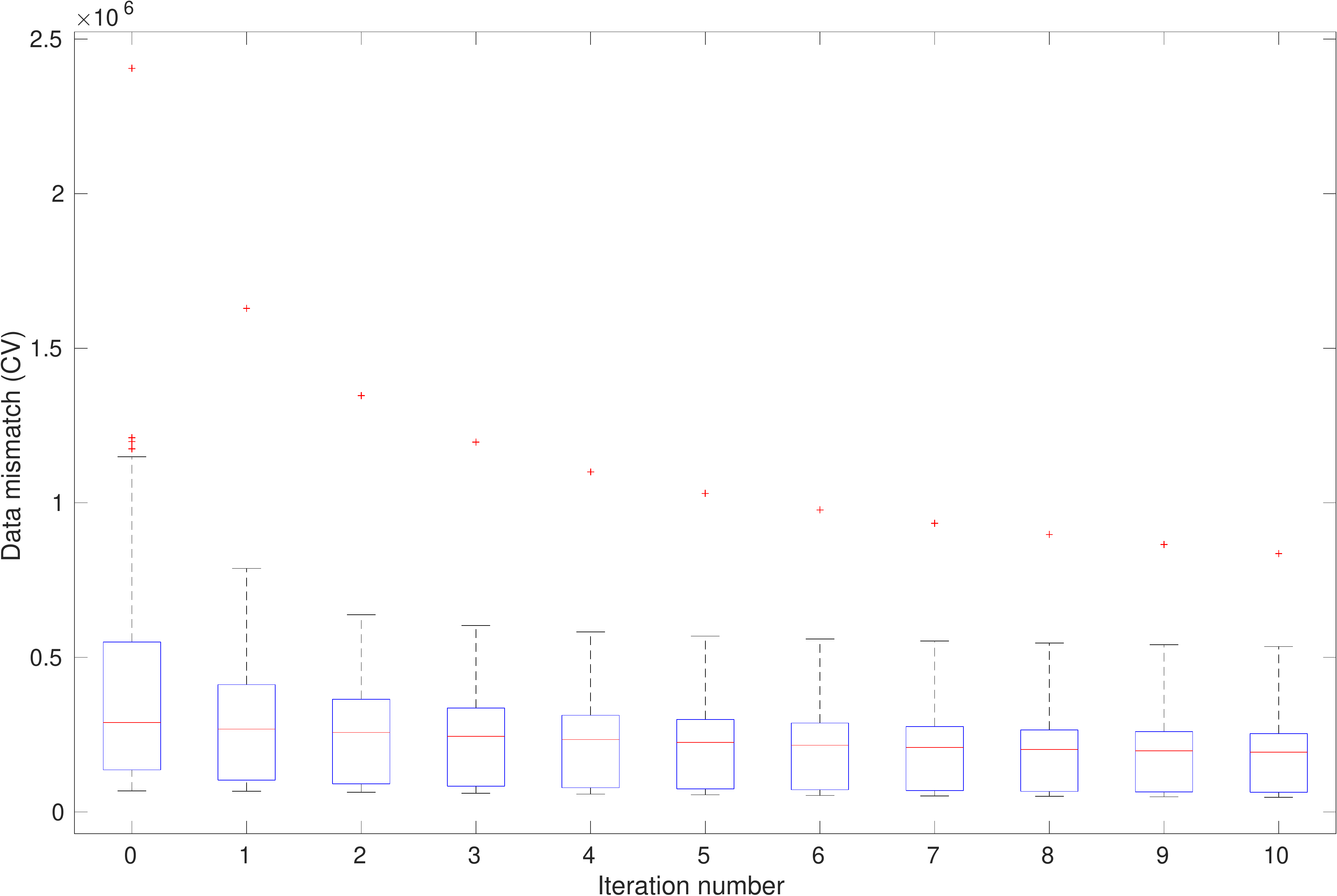}} 
	\end{tabular}
	\caption{\label{fig:boxplot_dataMismatch_cluster1_iter_unimodal} Box plots of data mismatch at different iteration steps, with respect to the (a) training and (b) CV datasets in case of unimodal inputs.}
\end{figure}

\renewcommand{\nScale}{0.2}
\begin{figure} 
	\centering
	\begin{tabular}{cc}
		\subfloat[Scale parameters ($\beta$)]{\includegraphics[width=0.45\textwidth]{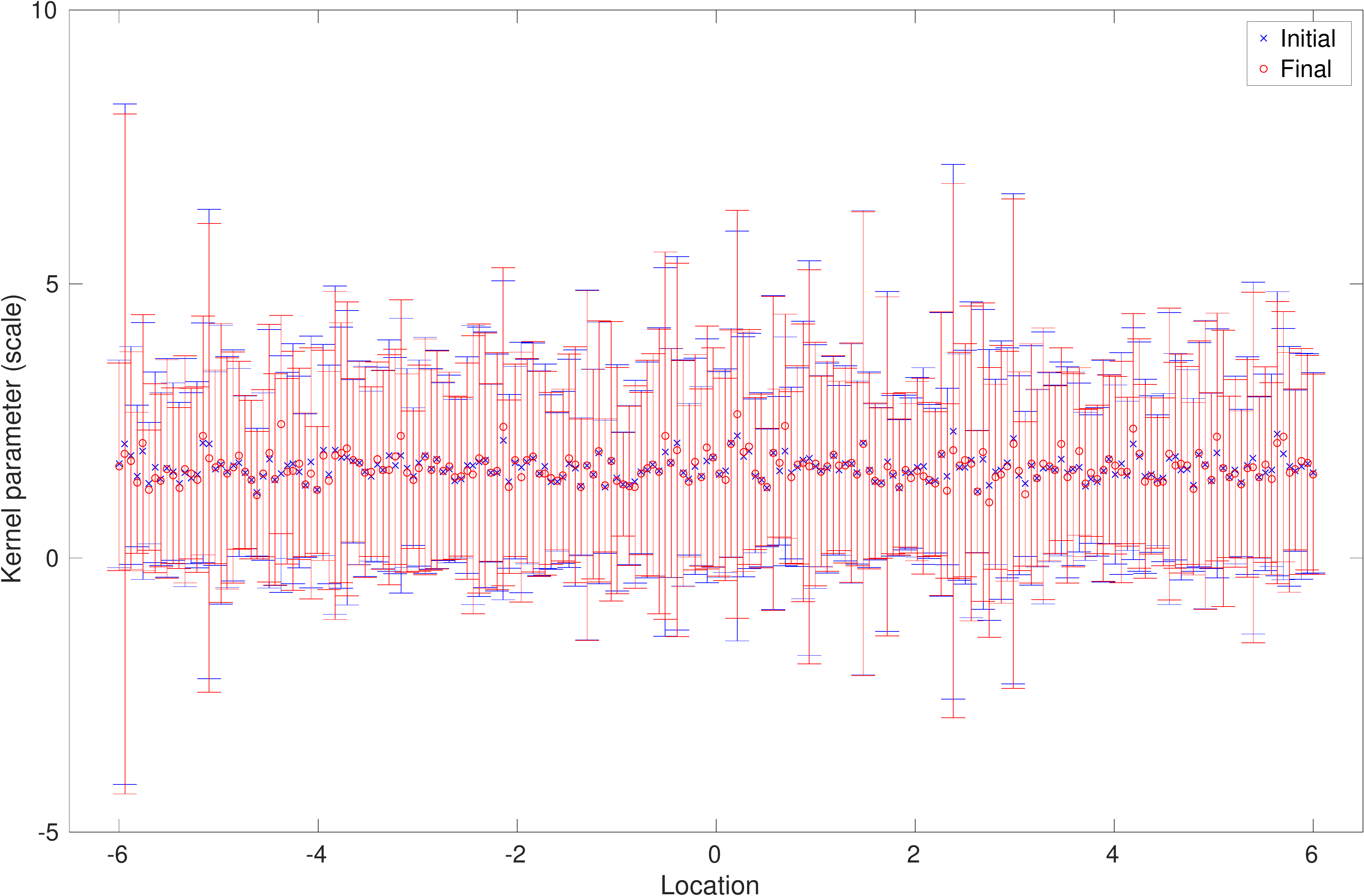}} &
		\subfloat[Weight parameters ($c$)]{\includegraphics[width=0.45\textwidth]{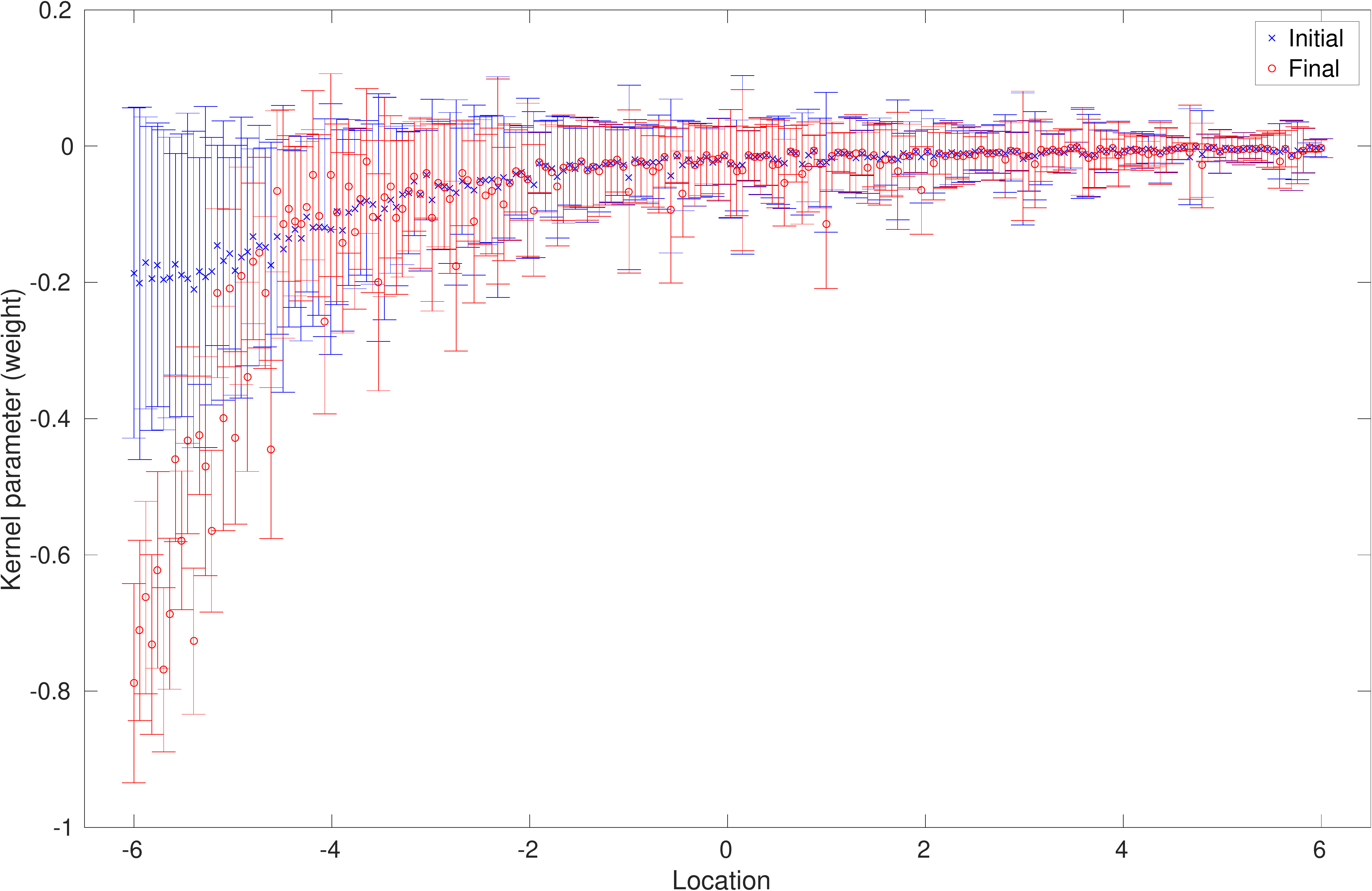}} 
	\end{tabular}
	\caption{\label{fig:errorbar_kernelPara_unimodal} Error-bar plots in case of unimodal training inputs, in the form of ensemble mean $\pm$ ensemble std, with respect to the initial (in blue) and final (in red) ensembles of scale (Panel (a)) and weight (Panel (b)) parameters, respectively, associated with $200$ center points that are evenly distributed over the interval $[-6, 6)$.}
\end{figure} 

\renewcommand{\nScale}{0.2}
\begin{figure} 
	\centering
	\begin{tabular}{cc}
		\subfloat[Initial ensemble of predictions]{\includegraphics[width=0.45\textwidth]{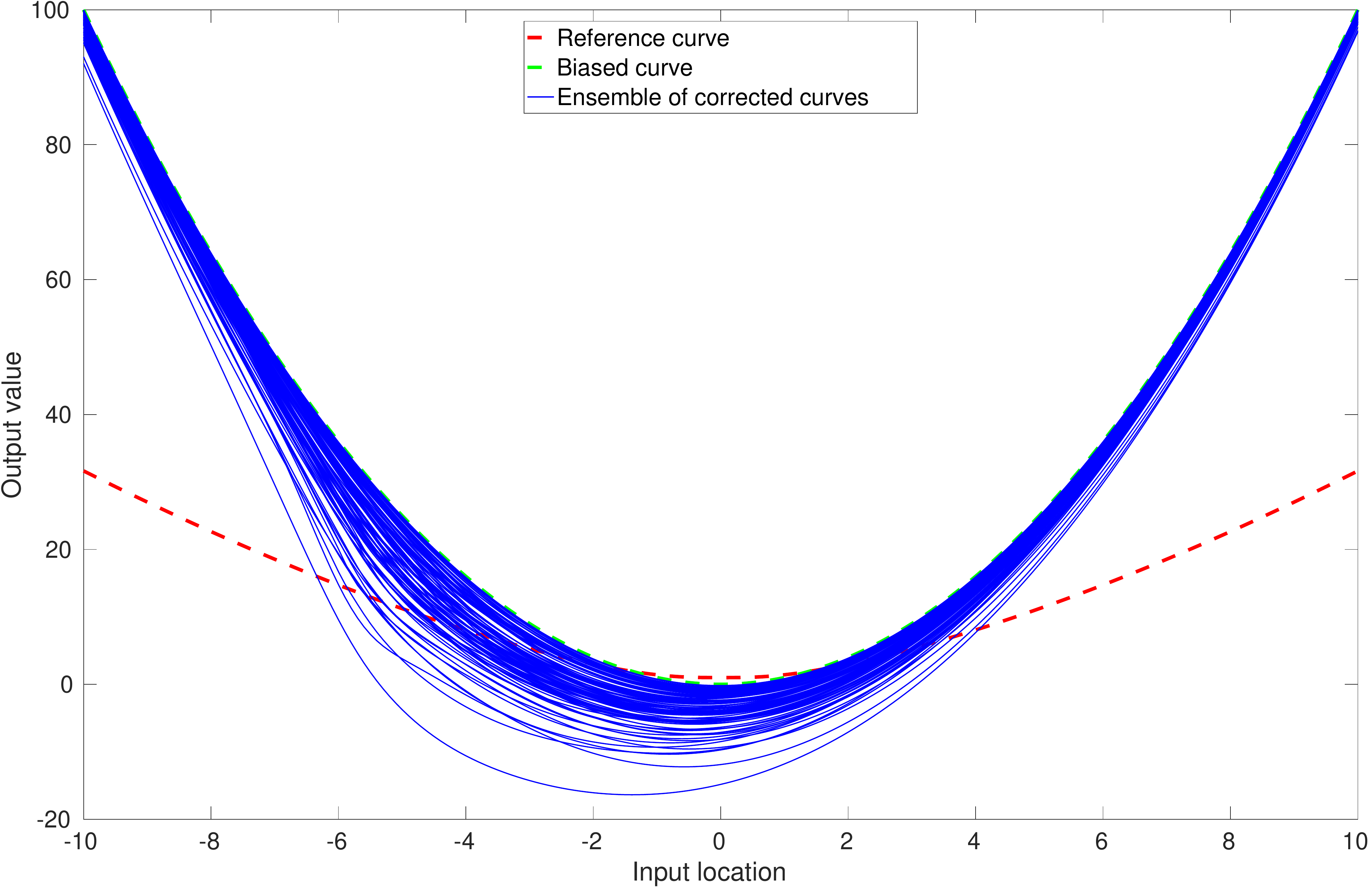}} &
		\subfloat[Initial mean predictions]{\includegraphics[width=0.45\textwidth]{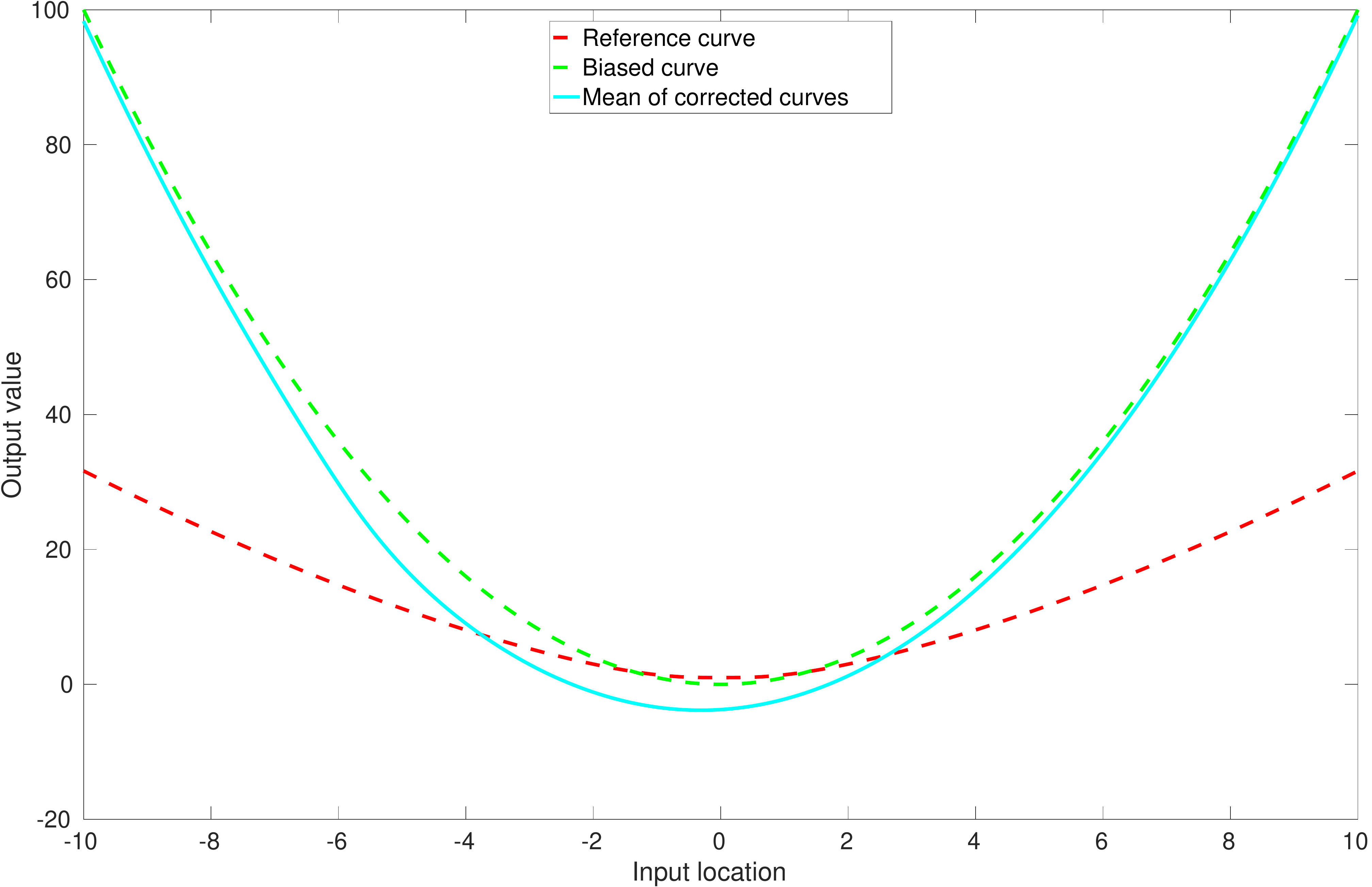}} \\
		\subfloat[Final ensemble of predictions]{\includegraphics[width=0.45\textwidth]{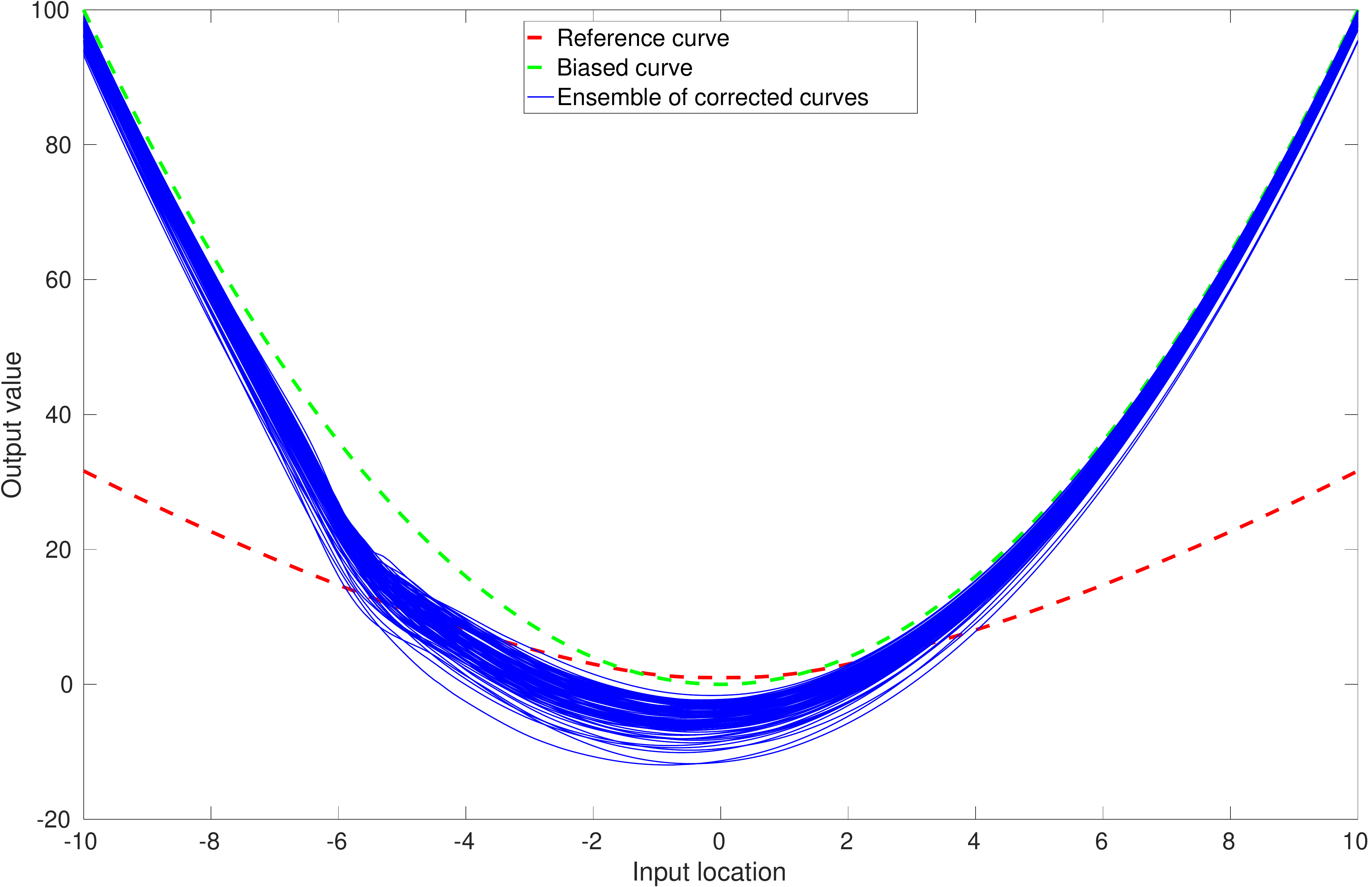}} &
		\subfloat[Final mean predictions]{\includegraphics[width=0.45\textwidth]{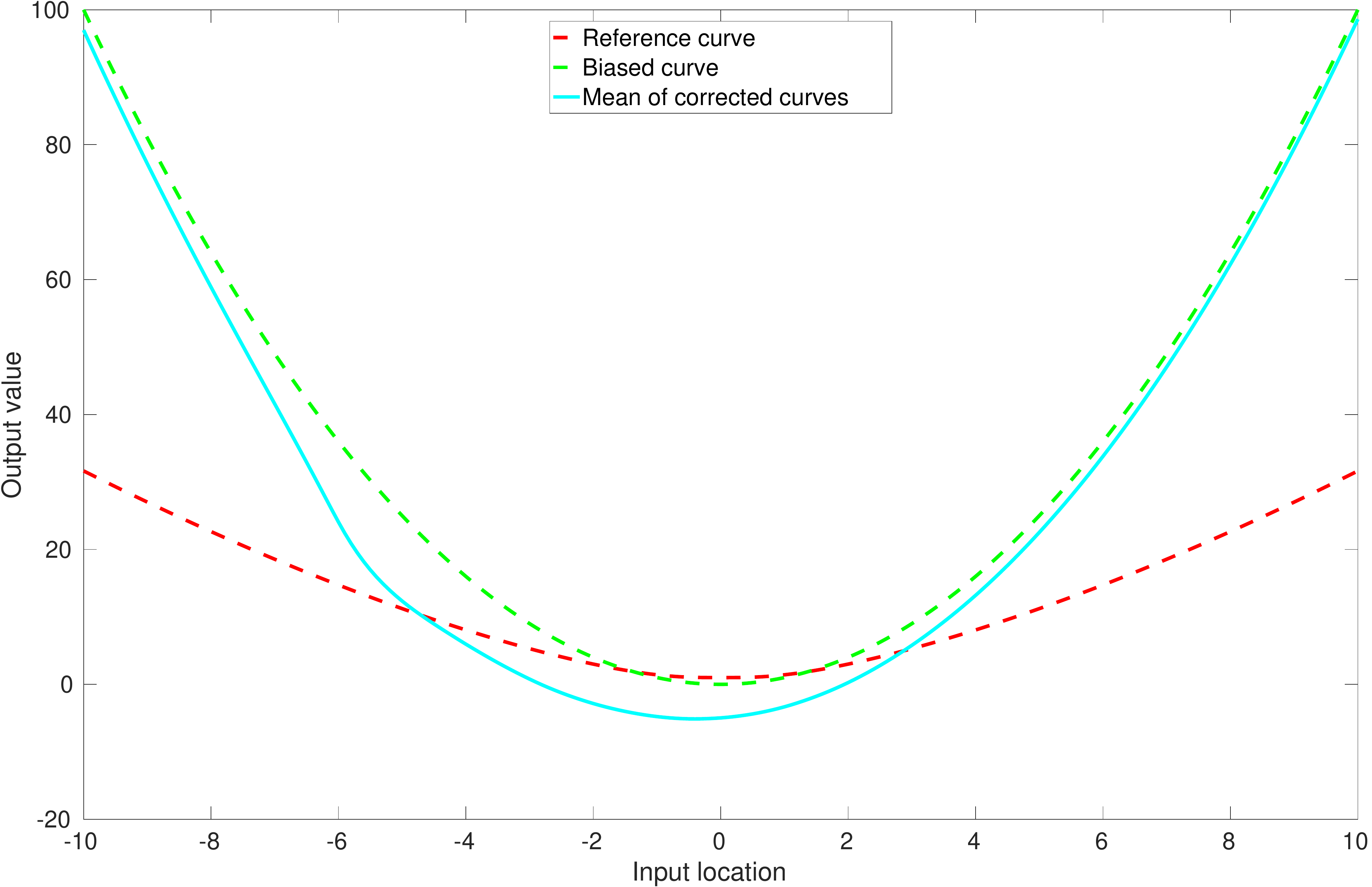}} \\ 
	\end{tabular}
	\caption{\label{fig:prediction_ensemlbe_cluster1} Red (reference curve) and green (biased curve) curves in all panels of the current figure are the same as those in Figure \ref{fig:ref_bias_func}. Panels (a) and (c) show the initial and final ensembles of predictions (with respect to the case of unimodal training inputs), obtained by adding to the biased curve the corresponding ensembles of residual terms, which are computed using Eqs. (\ref{eq:rbf_approximation}) through (\ref{eq:kernel_para}), wherein the kernel parameters correspond to the initial and final ensembles of scale and weight parameters, respectively. Panels (b) and (d) also report the means of the initial and final ensembles of predictions, respectively.}
\end{figure} 

Figure \ref{fig:boxplot_dataMismatch_cluster1_iter_unimodal} shows the box plots of data mismatch values generated by $100$ different sets of kernel (weight and scale) parameters at different iteration steps, where mismatch values are calculated using (a) training and (b) CV datasets, respectively. Note that in the course of learning, only training dataset is used to update kernel parameters. Therefore, it is not surprising to see that the reduction of training-data mismatch tends to be more significant than the reduction of CV-data mismatch. The more important observation in this case, however, is that, even though the CV dataset is not involved in training the kernel parameters, its corresponding data mismatch tends to decrease as the training (iteration) process goes on, which implies that the whole training  process appears useful and there is no need to stop the iteration earlier.    

Figure \ref{fig:errorbar_kernelPara_unimodal} depicts the error-bar plots (in terms of ensemble mean $\pm$ ensemble std) of kernel parameters, namely, scales ($\beta$) and weights ($c$), associated with $200$ center points that are evenly distributed over the interval $[-6,6)$. For scale parameters (Panel (a)), the relative changes from initial (in blue) to final (in red) values appear not so significant for all center points. In contrast, for weight parameters (Panel (b)), more substantial changes are spotted for center points located in, e.g., the interval $[-6, -2]$, whereas outside this interval, the changes tends to be less significant again. This does not appear to be surprising, if we take into account the coverage of training inputs (see the upper left panel of Figure \ref{fig:hist_unimodal}), and the fact that a Gaussian kernel function decays exponentially to zero as the distance between a training input and the center point associated with the kernel function increases. 

For further demonstration, in Panels (a) and (c) of Figure \ref{fig:prediction_ensemlbe_cluster1}, we also compare the reference curve (red) over the input interval $\left[-10,\, 10\right]$, the corresponding biased curve (green), and ensembles of corrected curves (blue), in the form of the biased curve plus ensembles of residual terms. The reference and biased curves are the same as those in Figure \ref{fig:ref_bias_func}, whereas the ensembles of residual terms are calculated using Eqs. (\ref{eq:rbf_approximation}) through (\ref{eq:kernel_para}), in which the kernel parameters correspond to either the initial or the final ensembles of scale and weight parameters, respectively. In addition, for better visualization, in panels (b) and (d), we plot the mean corrected curves (cyan).

Figures \ref{fig:prediction_ensemlbe_cluster1}(a) and (b) indicate that, compared to the biased curve, the way for us to initialize the initial ensembles of kernel parameters tends to improve the prediction accuracies over the interval (e.g., $[-6,-4]$) on which the training inputs largely concentrate (called concentration interval hereafter). In addition, the resulting ensemble of corrected predictions provides a means of conducting uncertainty quantification for the predictions. 
Recall that we adopt only $100$ (as the ensemble size) random training input-output (label) pairs to initialize the ensembles of kernel parameters. By learning (or updating) kernel parameters through more training data, it appears that the performances of both prediction and uncertainty quantification are improved over the concentration interval, as can be seen in Figures \ref{fig:prediction_ensemlbe_cluster1}(c) and (d).    

Figure \ref{fig:prediction_ensemlbe_cluster1} also shows that, for the intervals (e.g., $[-2,2]$) over which there are sparse or even no training data, the corrected predictions may be less accurate than the biased predictions themselves. In this case, a natural way to improve the performance of supervised learning is to acquire more training data over different regions. Such an investigation will be carried out in the next sub-section.      

\subsection*{Results with respect to multi-modal inputs}  
To generate more training data, here we consider a scenario with multi-modal training inputs. We will first identify a challenge for the ensemble-based learning algorithm to handle multi-modal training inputs, and then investigate a strategy that helps overcome this problem.   

\begin{figure} 
	\centering
	\includegraphics[width=0.8\textwidth]{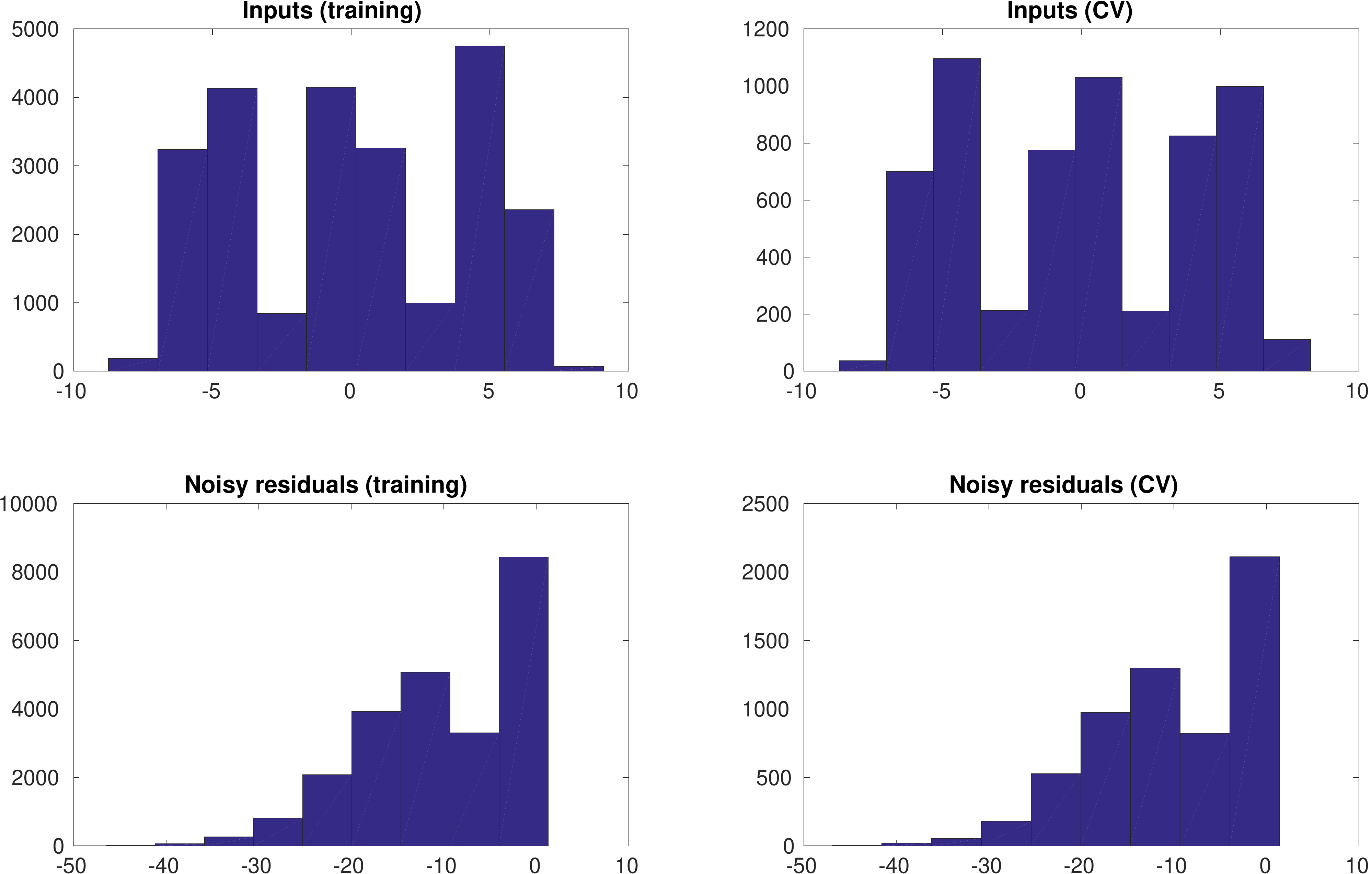}
	\caption{\label{fig:hist_multimodal} Histograms of the multi-modal inputs and noisy residuals, with respect to the training and CV datasets, respectively.}
\end{figure}

In the experiment, we generate a set of 10, 000 input samples from the distribution $N(-5,1)$, 10, 000 input samples from the distribution $N(0,1)$ and 10, 000 input samples from the distribution $N(5,1)$, and the corresponding noisy residuals. We then randomly split the resulting $30,000$ input-residual pairs into one training dataset (with $24,000$ data points) and one CV dataset (with $6,000$ data points). Figure \ref{fig:hist_multimodal} shows the histograms of the inputs and noisy residuals in the training and CV datasets.          
\renewcommand{\nScale}{0.3}

\renewcommand{\nScale}{0.2}
\begin{figure} 
	\centering
	\begin{tabular}{cc}
		\subfloat[Initial ensemble of predictions]{\includegraphics[width=0.45\textwidth]{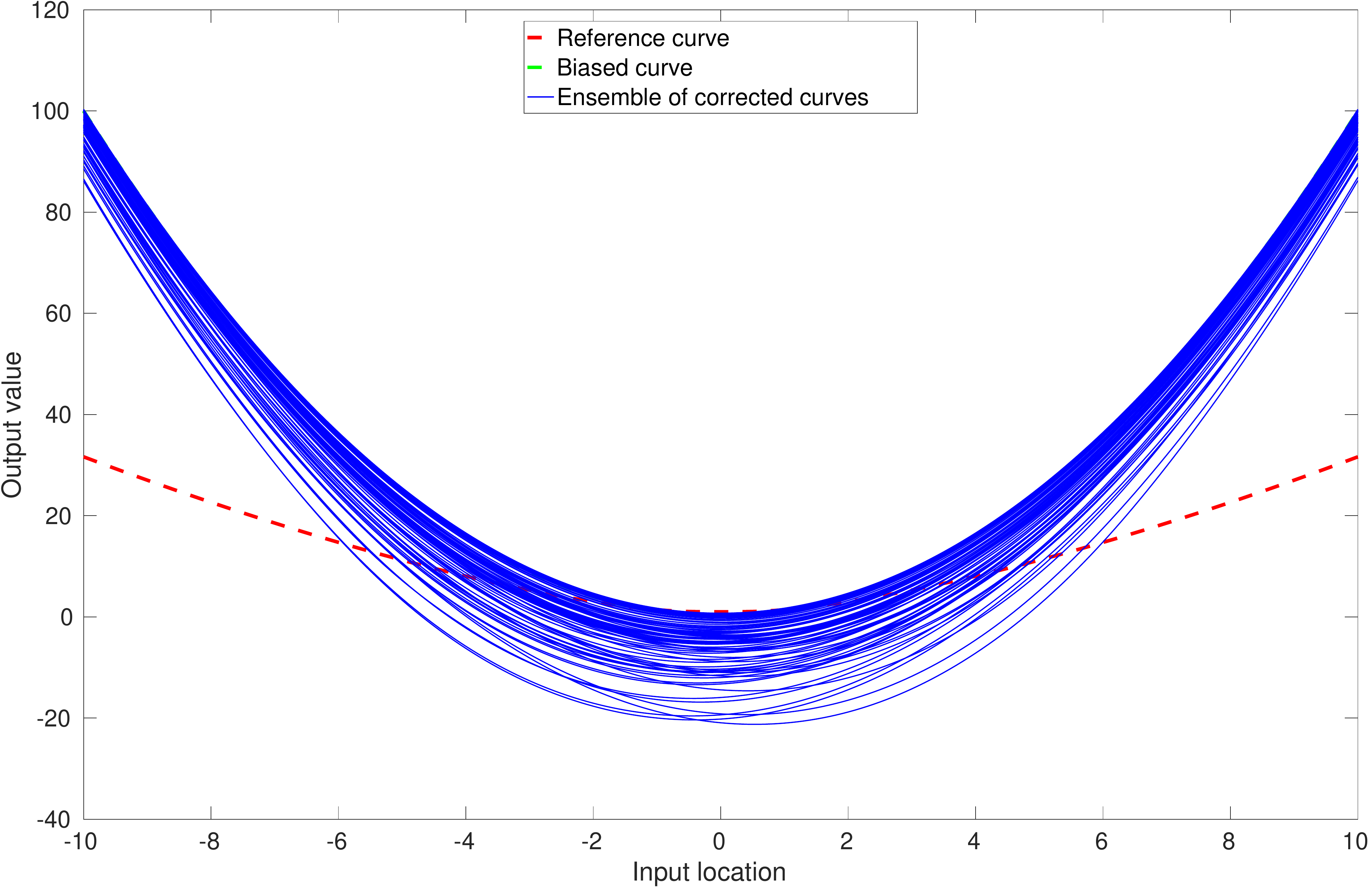}} &
		\subfloat[Initial mean predictions]{\includegraphics[width=0.45\textwidth]{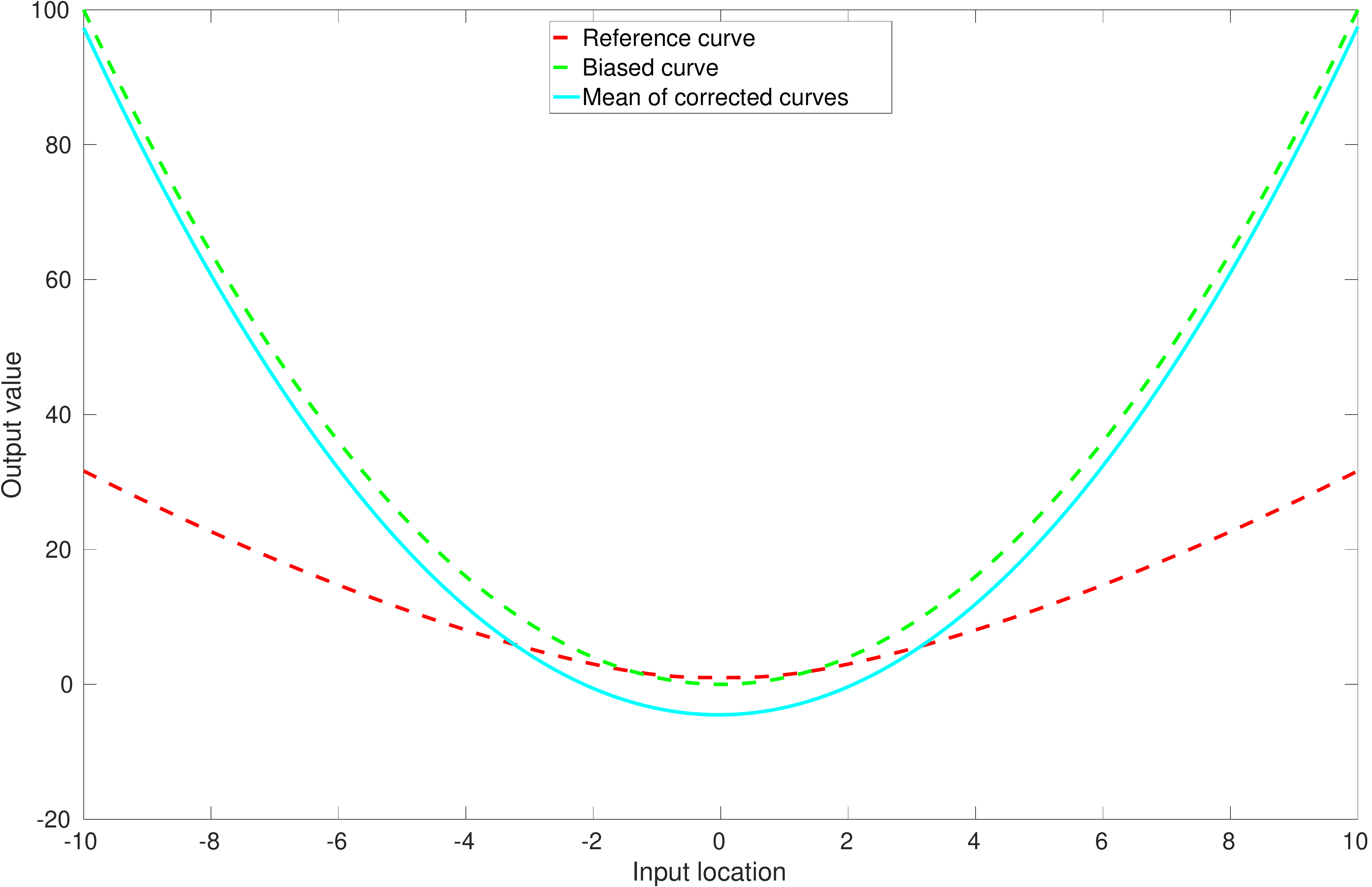}} \\
		\subfloat[Final ensemble of predictions]{\includegraphics[width=0.45\textwidth]{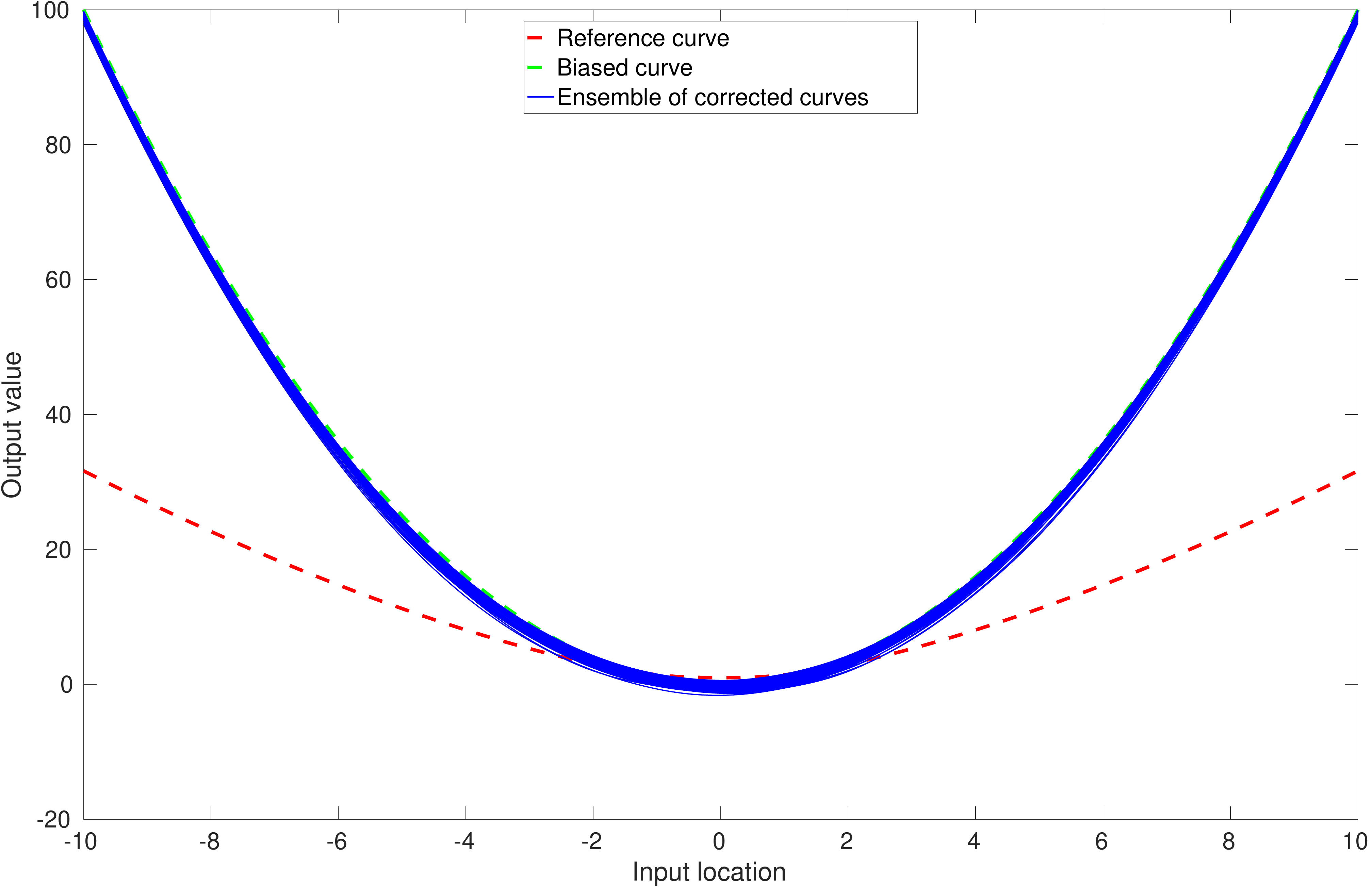}} &
		\subfloat[Final mean predictions]{\includegraphics[width=0.45\textwidth]{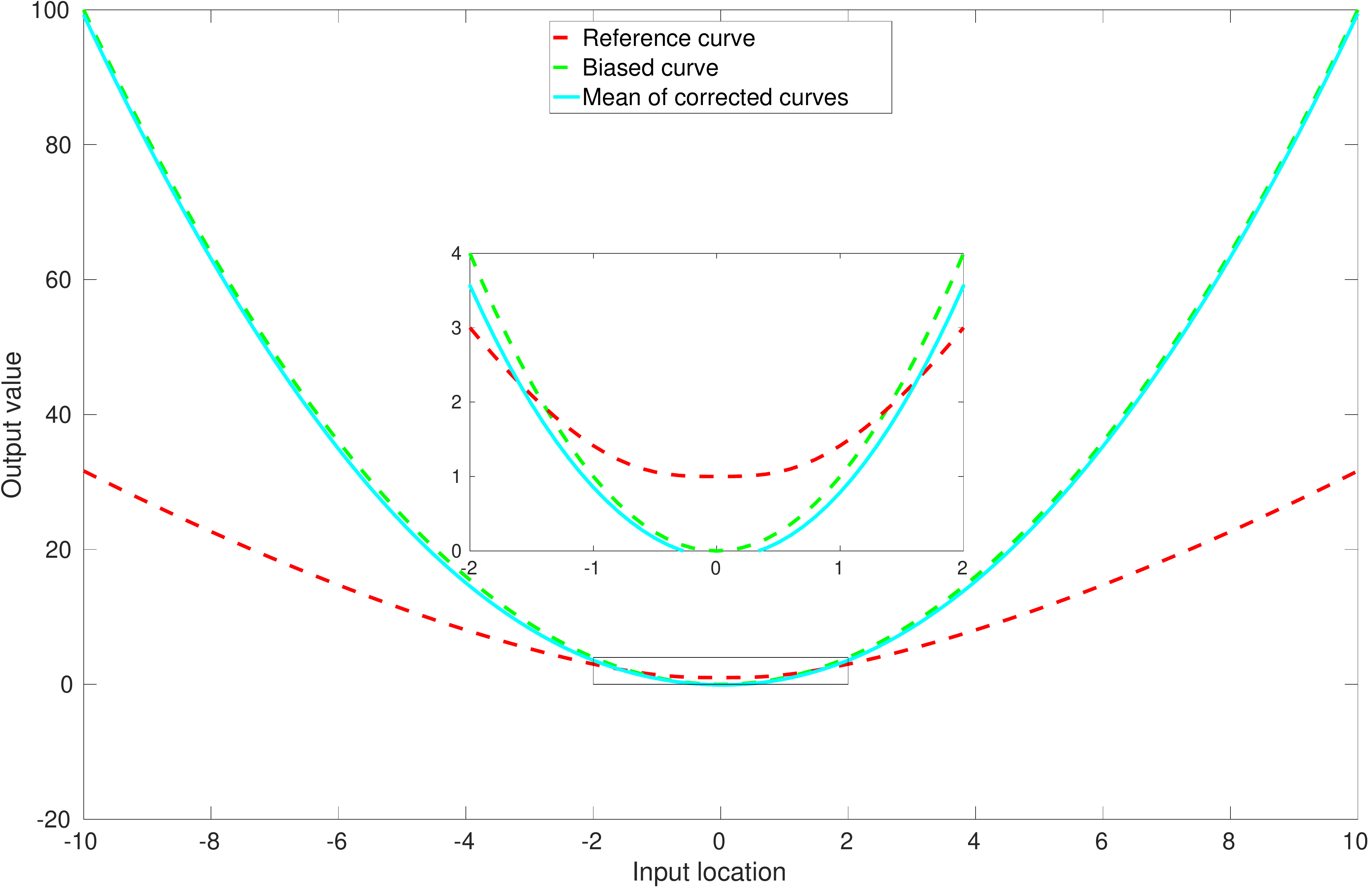}} \\ 
	\end{tabular}
	\caption{\label{fig:multimodal_prediction_ensemlbe_cluster1} As in Figure \ref{fig:prediction_ensemlbe_cluster1}, but for the case with multi-modal inputs, for which no multi-modal learning strategy (MMLS) is adopted. For better visualization, in Panel (d) we re-plot the reference, the biased and the mean corrected curves over the input interval $\left[-2, \, 2\right]$ in a separate, zoomed-in subplot.}
\end{figure} 

In the sequel, we first illustrate what will happen if one directly applies the ensemble-based learning algorithm to the training data with multi-modal inputs. In the experiment, we still adopt $200$ center points that are evenly distributed over the interval $\left[-6, 6\right)$. As in the previous sub-section, the ensemble-based learning algorithm is directly adopted to update $400$ kernel parameters, namely, the scale and weight parameters associated with each center point. However, it turns out that, in the presence of multi-modal inputs, a straightforward application of the ensemble-based learning algorithm may not achieve satisfactory performance. This point is demonstrated in Figure \ref{fig:multimodal_prediction_ensemlbe_cluster1}. With more training data than in the previous sub-section, the accuracies of corrected predictions over certain input intervals, e.g., $[-6, -4]$, are actually worsened, as is evident if one compares panels (c) and (d) of Figures \ref{fig:prediction_ensemlbe_cluster1} and \ref{fig:multimodal_prediction_ensemlbe_cluster1}.

The under-performance of plain, ensemble-based algorithms in handling multi-modal variables is also discussed in the literature, see, for example, \cite{elsheikh2013clustered,gao2017distributed,gao2018gaussian,Hoteit2012,Luo2008-spgsf1}. To deal with this problem, we equip the ensemble-based algorithm with an multi-modal learning strategy (MMLS). Concretely, similar to \cite{elsheikh2013clustered,gao2018gaussian,Hoteit2012,Luo2008-spgsf1}, we adopt a Gaussian mixture model (GMM) to fit the probability density function (pdf) of multi-modal variables, which naturally leads to a number of clustered subsets of multi-modal variables (and their corresponding noisy labels). Next, we use each cluster of training data to initialize (and then update) an ensemble of kernel parameters that are associated with the center points (note that different clusters of training data share the same set of center points). The corrected predictions are then taken as biased outputs plus certain residual terms, whereas the latter are calculated as the weighted averages of the residuals predicted using the kernel parameters in each cluster.

More specifically, suppose that the multi-modal inputs are clustered into $N_{cl}$ mutually exclusive subsets, and in each subset, the pdf of the inputs is modelled by a certain Gaussian pdf. In other words, the pdf $p(x)$ of training inputs is approximated by a GMM, in the form of
\begin{linenomath*}  
	\begin{IEEEeqnarray}{r} \label{eq:GMM}
		p(x) \approx \sum_{s=1}^{N_{cl}} w_s \, n(x;\mu_s,\sigma_s^2) \, ,
	\end{IEEEeqnarray}
\end{linenomath*}
where $w_s$ is the weight associated with the $s$th cluster, and $n(x;\mu_s,\sigma_s^2)$ is the corresponding Gaussian pdf, parametrized by the mean $\mu_s$ and STD $\sigma_s$. For each cluster, say, the $s$th one, we generate an initial ensemble $\boldsymbol{\Theta}_s^0 \equiv \{\boldsymbol{\theta}_{j,s}^0\}_{j=1}^{N_e}$ of kernel parameters for the set of center points, in the same way as in the preceding sub-subsection. The ensemble $\boldsymbol{\Theta}_s^0$ is then updated to $\boldsymbol{\Theta}_s^u \equiv \{\boldsymbol{\theta}_{j,s}^u\}_{j=1}^{N_e}$, using the training data associated with the cluster. With the above quantities, we are then able to generate corrected predictions for new inputs. For instance, given an input $x'$, we first calculate the probability $P_{s}$ of $x'$ with respect to each cluster, through
\begin{linenomath*}  
	\begin{IEEEeqnarray}{r} \label{eq:GMM_probability}
		P_{s}(x') = \dfrac{w_s n(x';\mu_s,\sigma_s^2)}{\sum_{s'=1}^{N_{cl}} w_{s'} \, n(x';\mu_{s'},\sigma_{s'}^2)} \, .
	\end{IEEEeqnarray}
\end{linenomath*}                         
Then, we can calculate an ensemble of corrected predictions, in the form of biased prediction $g(x')$ plus predicted residual $\hat{r}_j(x')$ ($j=1,2,\dotsb, N_e$), where $\hat{r}_j(x')$ is given by
\begin{linenomath*}  
	\begin{IEEEeqnarray}{r} \label{eq:predicted_residual}
		\hat{r}_j(x') = \sum_{s=1}^{N_{cl}} P_{s}(x') \hat{h}(x'; \boldsymbol{\theta}_{j,s}^u) \, , j=1,2,\dotsb, N_e \, ,
	\end{IEEEeqnarray}
\end{linenomath*}    
with $\boldsymbol{\theta}_{j,s}^u \in \boldsymbol{\Theta}_s^u$, and $\hat{h}$ a functional consisting of a set of kernel functions (cf. Eq. (\ref{eq:rbf_approximation}) or (\ref{eq:mt_rbf_approximation})) parametrized by $\boldsymbol{\theta}_{j,s}^u$. 

In terms of parametrization strategy adopted in SLP, a noticeable feature in case of multi-modal training inputs is that, each cluster of training data will have its own ensemble of kernel parameters associated with the same set of center points. Following the discussion in the text after Eq. (\ref{eq:mt_kernel_para}), for $m$-dimensional inputs, the total number of kernel parameters then becomes $(m+1) \times N_{cp} \times N_{cl}$, larger than that in case of unimodal training inputs. This may thus be considered as an additional way to improve the capacity of a learning model. 

\subsubsection*{Results with $N_{cl} = 3$}
\begin{table*} 
	\centering
	\captionsetup{width=.6\textwidth}
	\caption{\label{tab:GMM_components} Number of training data points associated with each GMM component (cluster), and the corresponding parameters estimated using the training data.}
	\begin{tabular}{||c||c||c||c||}
		\hline  
		& Cluster 1 (C1) & Cluster 2 (C2) & Cluster 3 (C3)  \\ 
		\hline		
		Number of data		& $8003$  & $8003$ & $7994$ \\
		\hline 			
		Estimated weight  & $0.3331$ & $0.3344$ & $0.3325$  \\
		\hline 			
		Estimated mean  & $5.0015$ & $-0.0074$ & $-5.0049$  \\
		\hline 			
		Estimated variance  & $0.9754$ & $1.0155$ & $0.9680$  \\
		\hline 
	\end{tabular}
\end{table*}   

%
\renewcommand{\nScale}{0.2}
\begin{figure} 
	\centering
	\begin{tabular}{cc}
		\subfloat[Training dataset (C1)]{\includegraphics[width=0.45\textwidth]{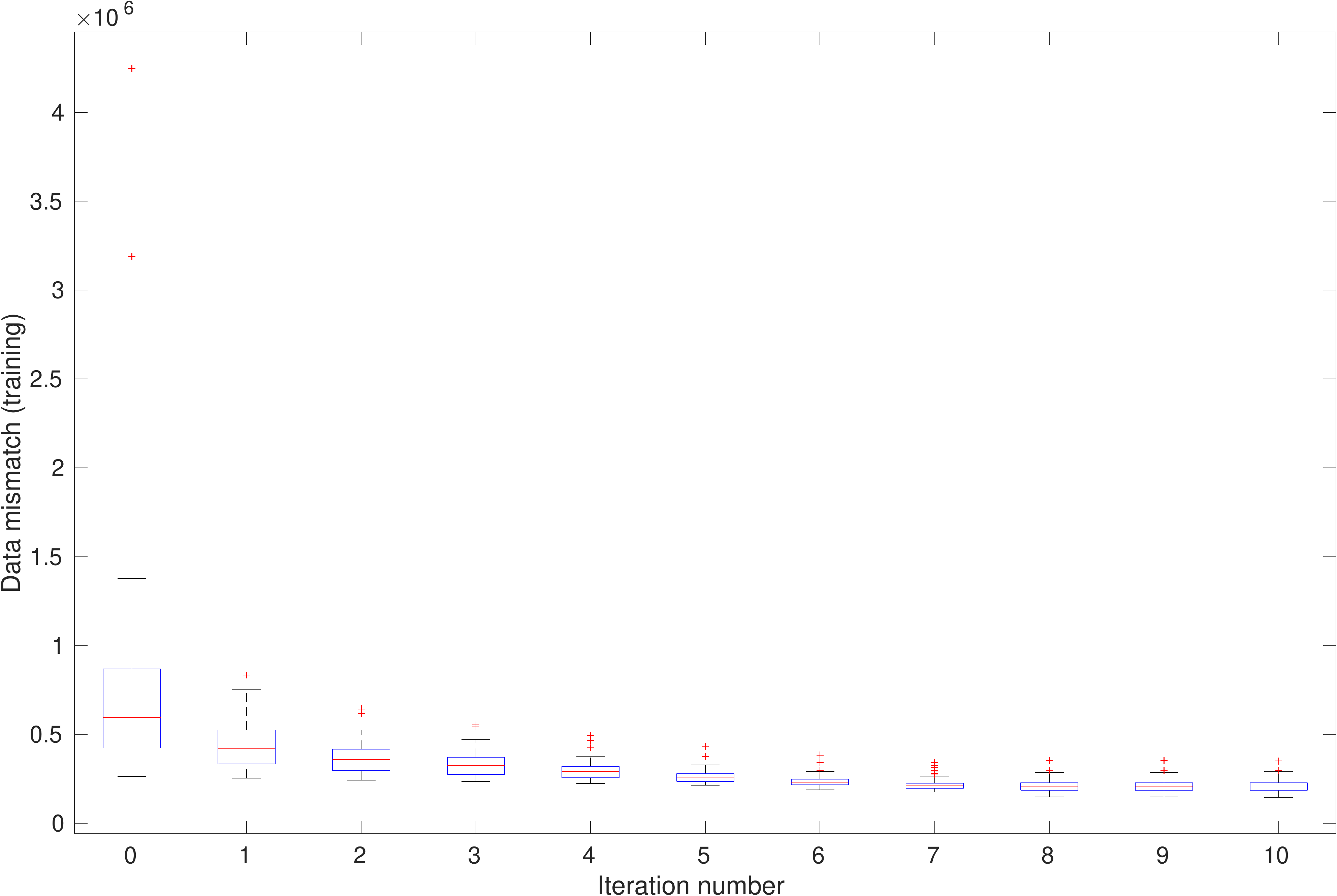}} &
		\subfloat[CV dataset (with respect to C1)]{\includegraphics[width=0.45\textwidth]{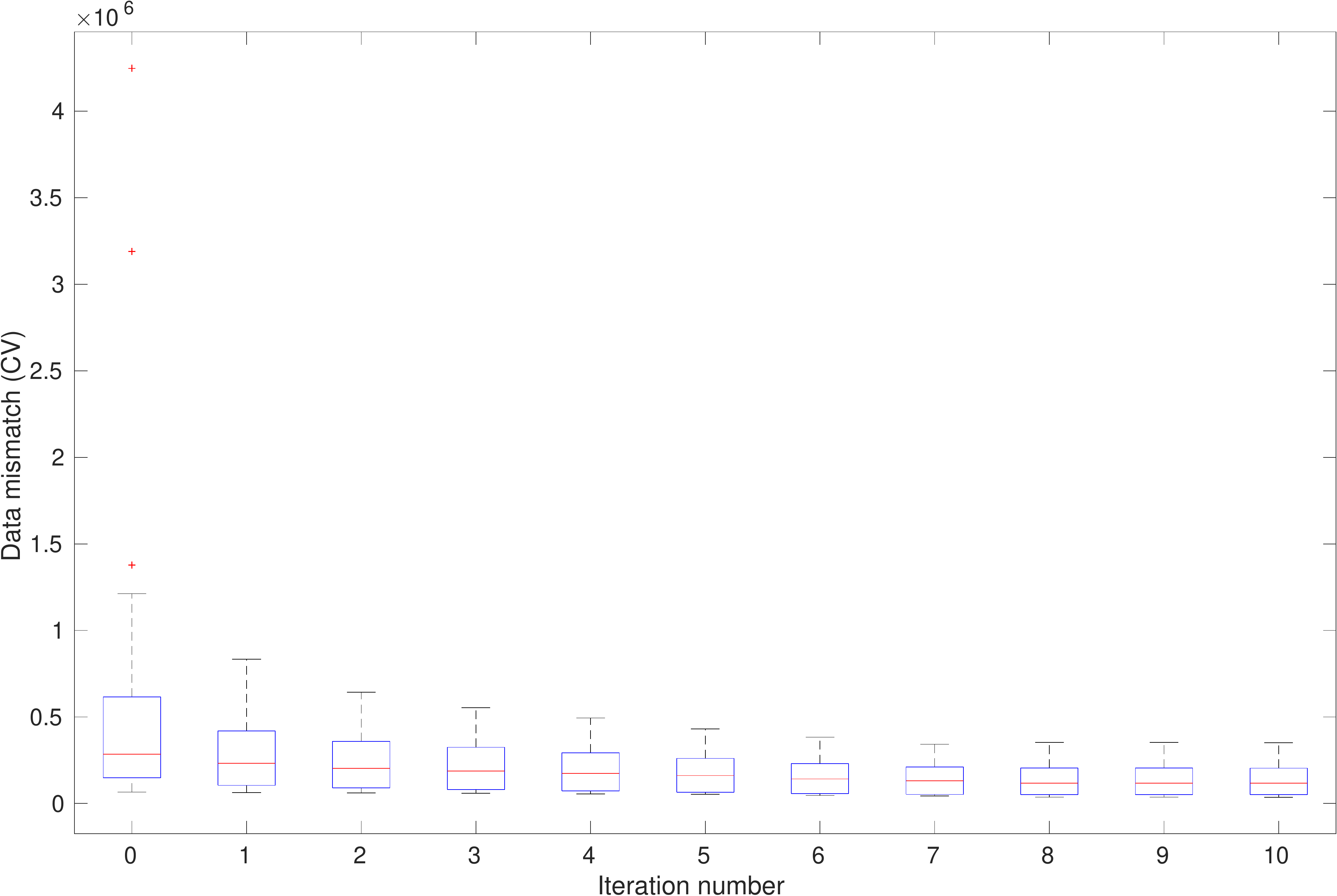}} \\
		\subfloat[Training dataset (C2)]{\includegraphics[width=0.45\textwidth]{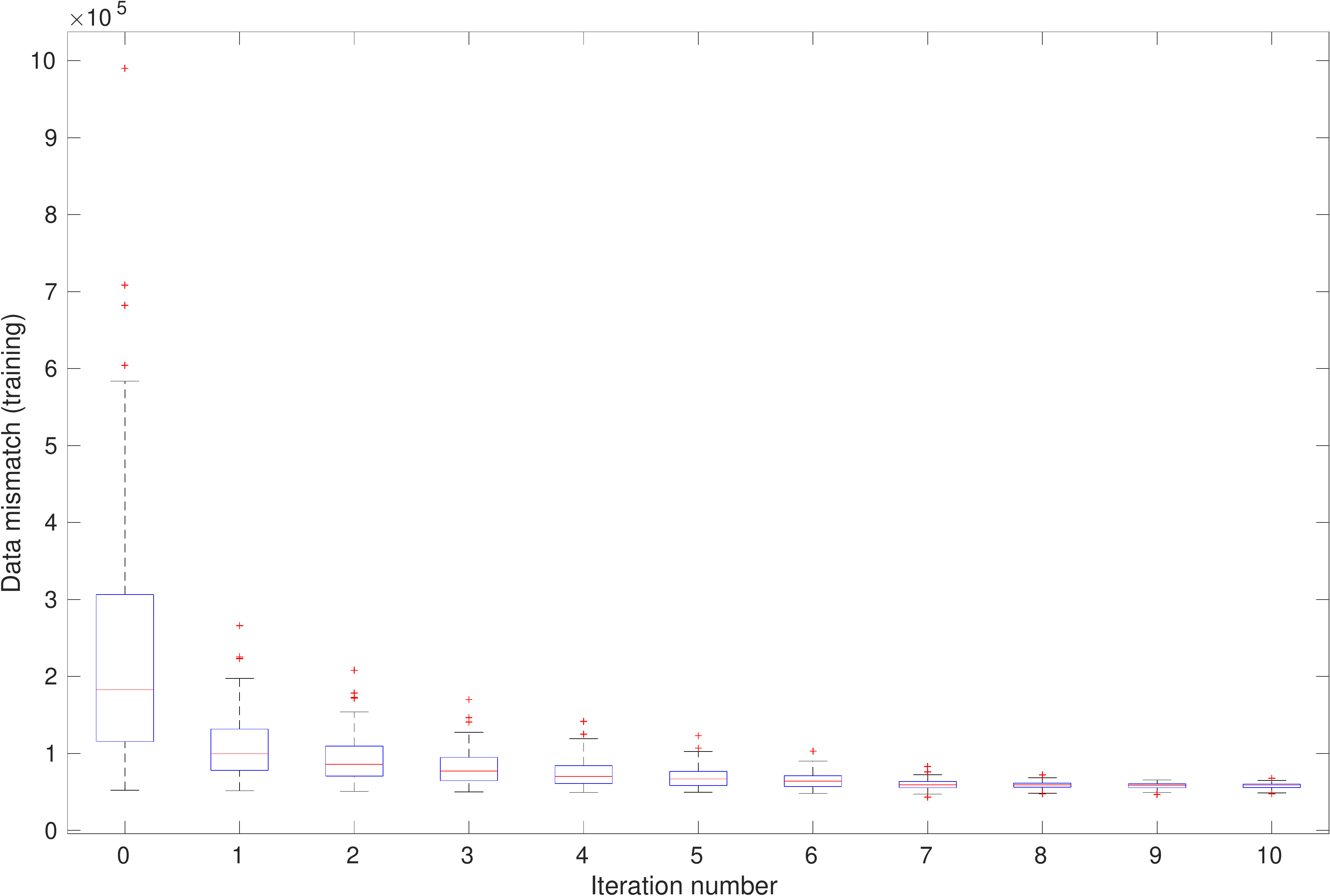}} &
		\subfloat[CV dataset (with respect to C2)]{\includegraphics[width=0.45\textwidth]{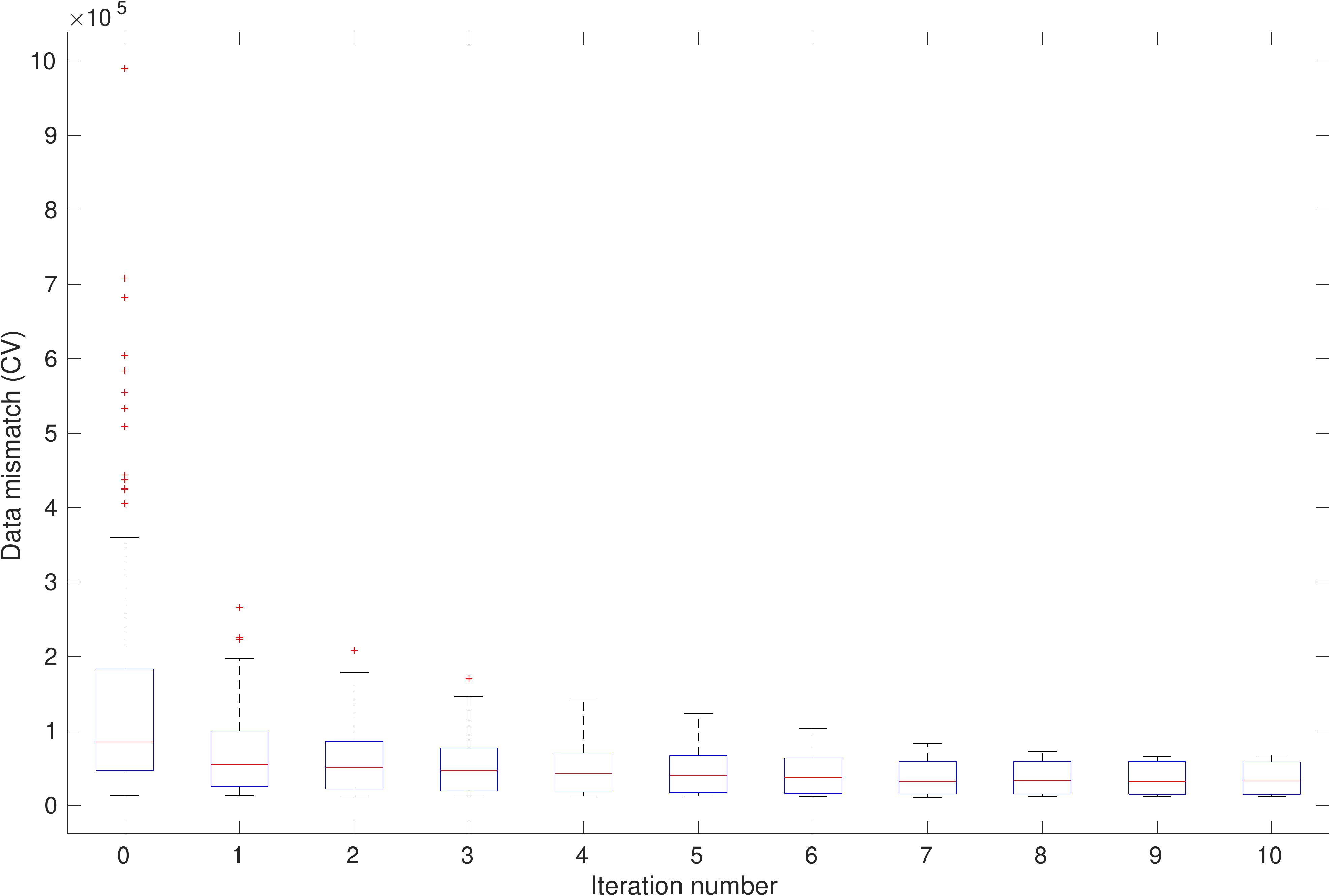}} \\
		\subfloat[Training dataset (C3)]{\includegraphics[width=0.45\textwidth]{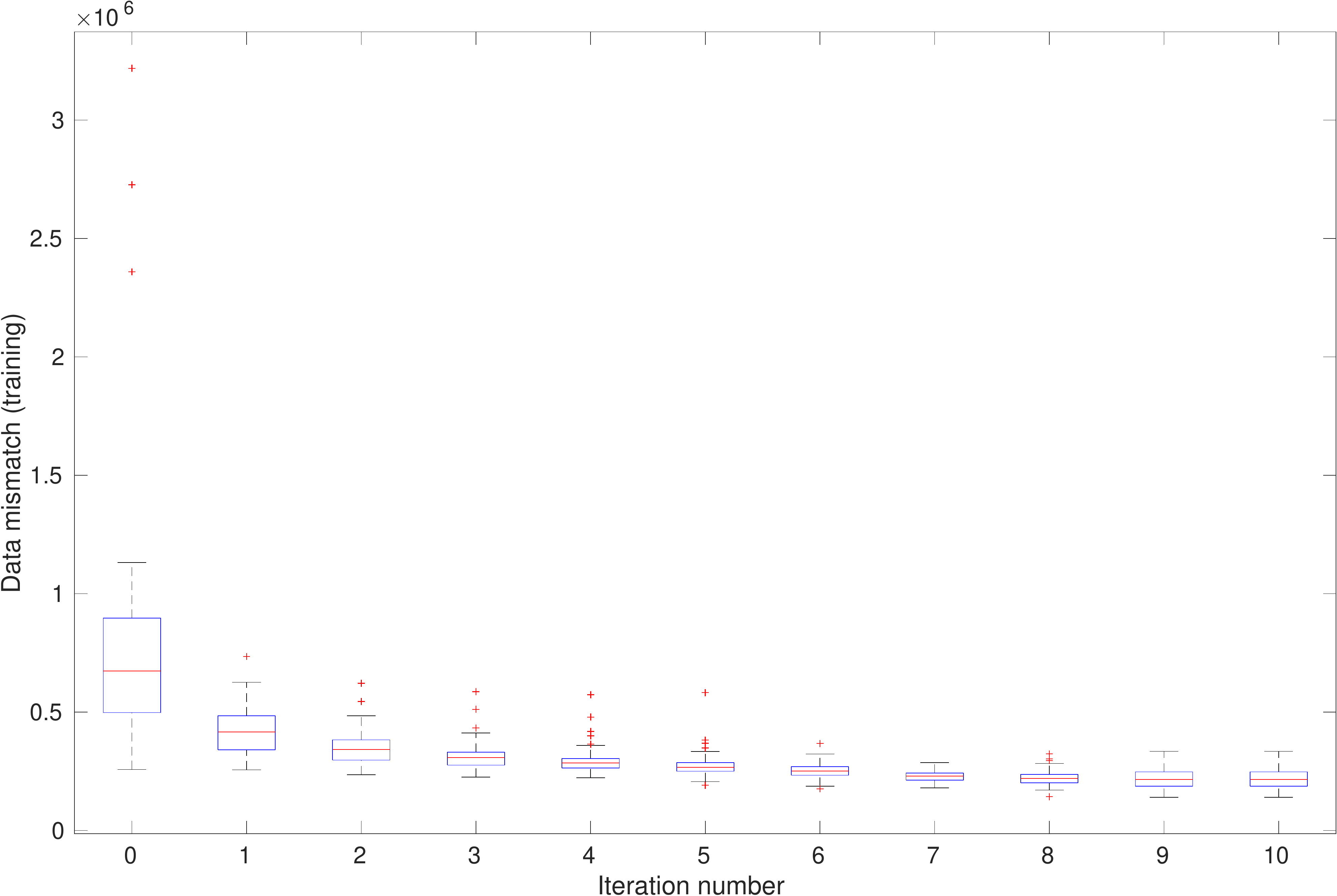}} &
		\subfloat[CV dataset (with respect to C3)]{\includegraphics[width=0.45\textwidth]{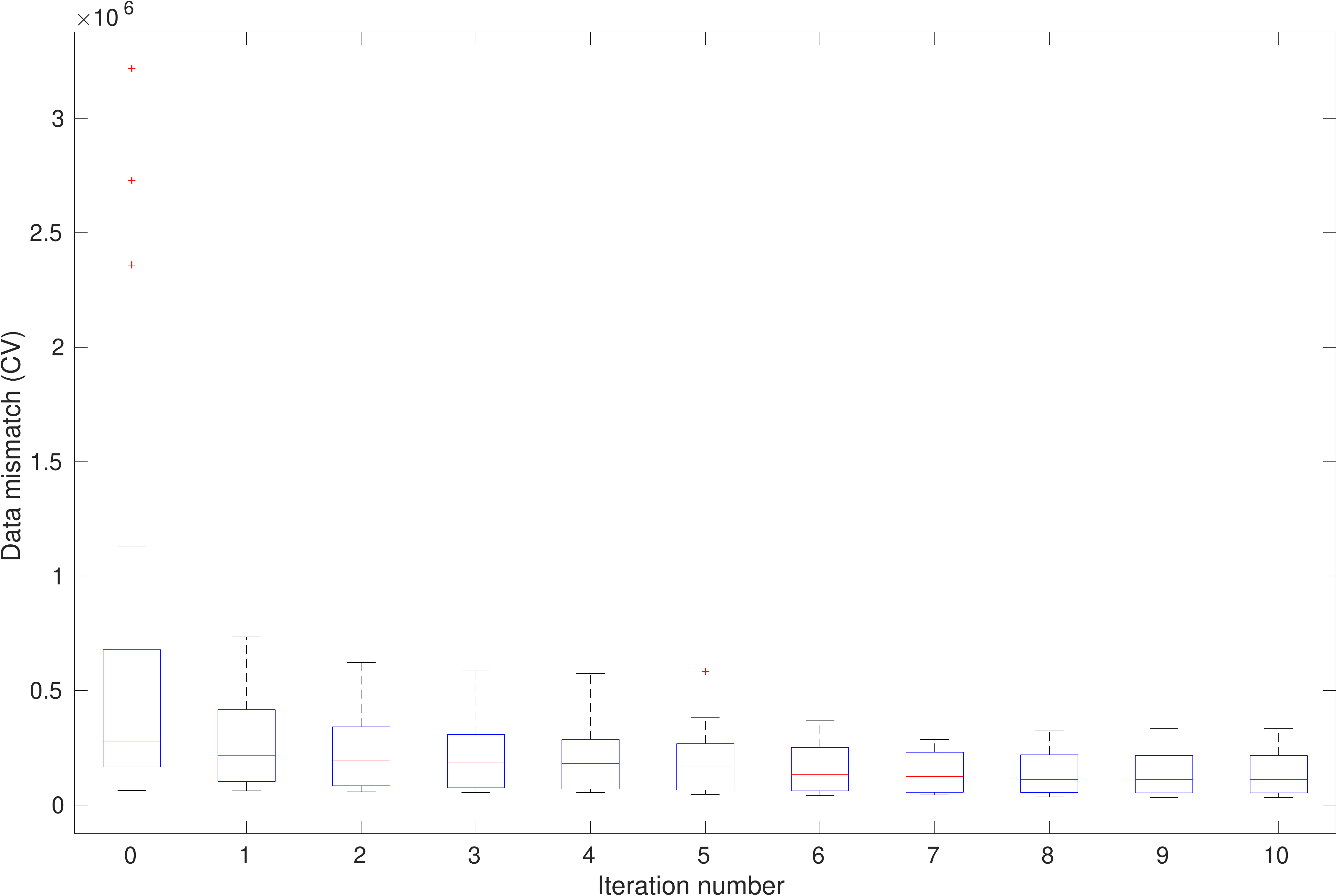}} \\
	\end{tabular}
	\caption{\label{fig:boxplot_dataMismatch_multimodal} Box plots of data mismatch at different iteration steps, with respect to the (a) training and (b) CV datasets in case of multi-modal inputs.}
\end{figure}

\renewcommand{\nScale}{0.2}
\begin{figure} 
	\centering
	\begin{tabular}{cc}
		\subfloat[Scale parameters ($\beta$) associated with C1]{\includegraphics[width=0.45\textwidth]{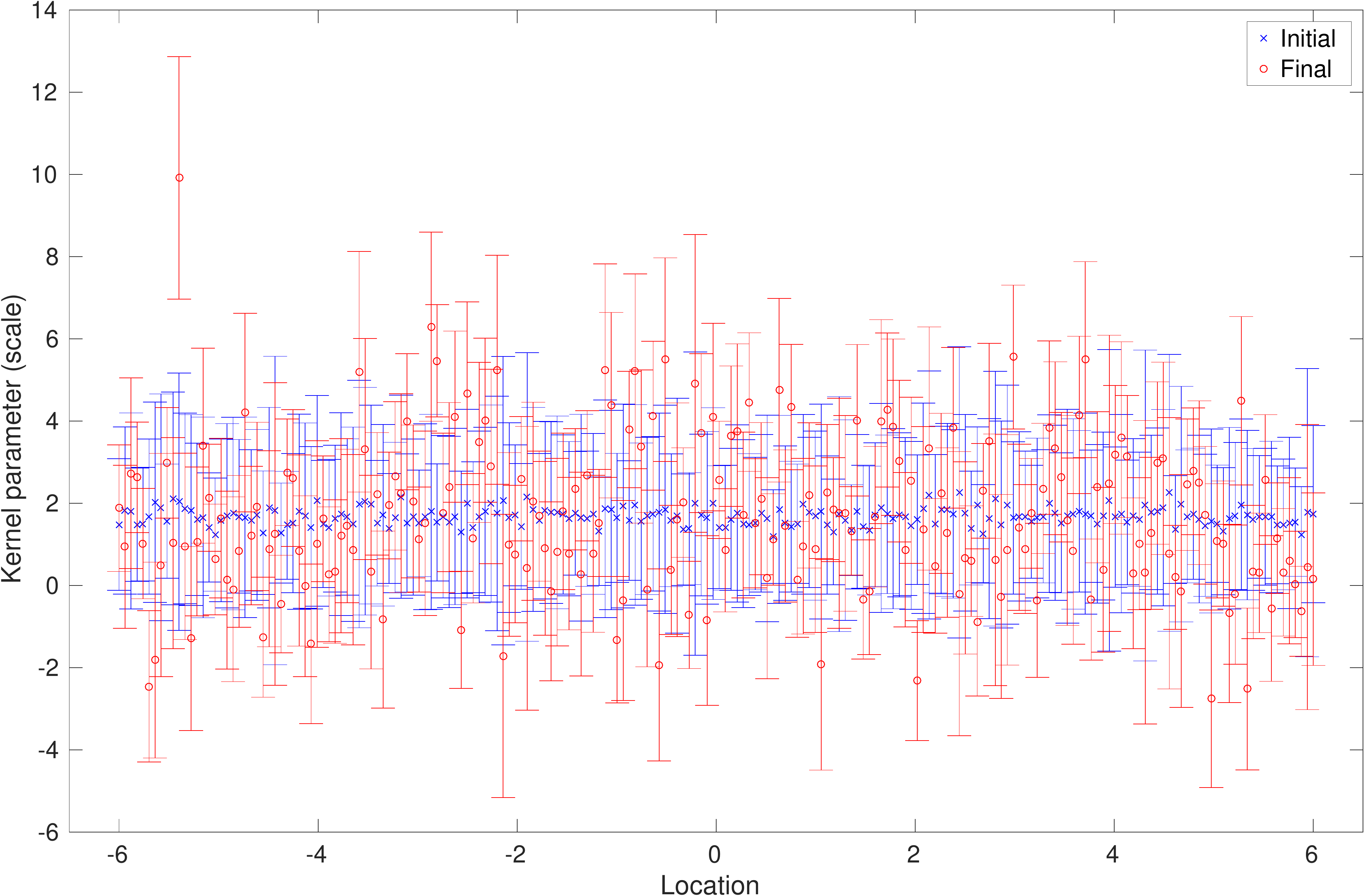}} &
		\subfloat[Weight parameters ($c$) associated with C1]{\includegraphics[width=0.45\textwidth]{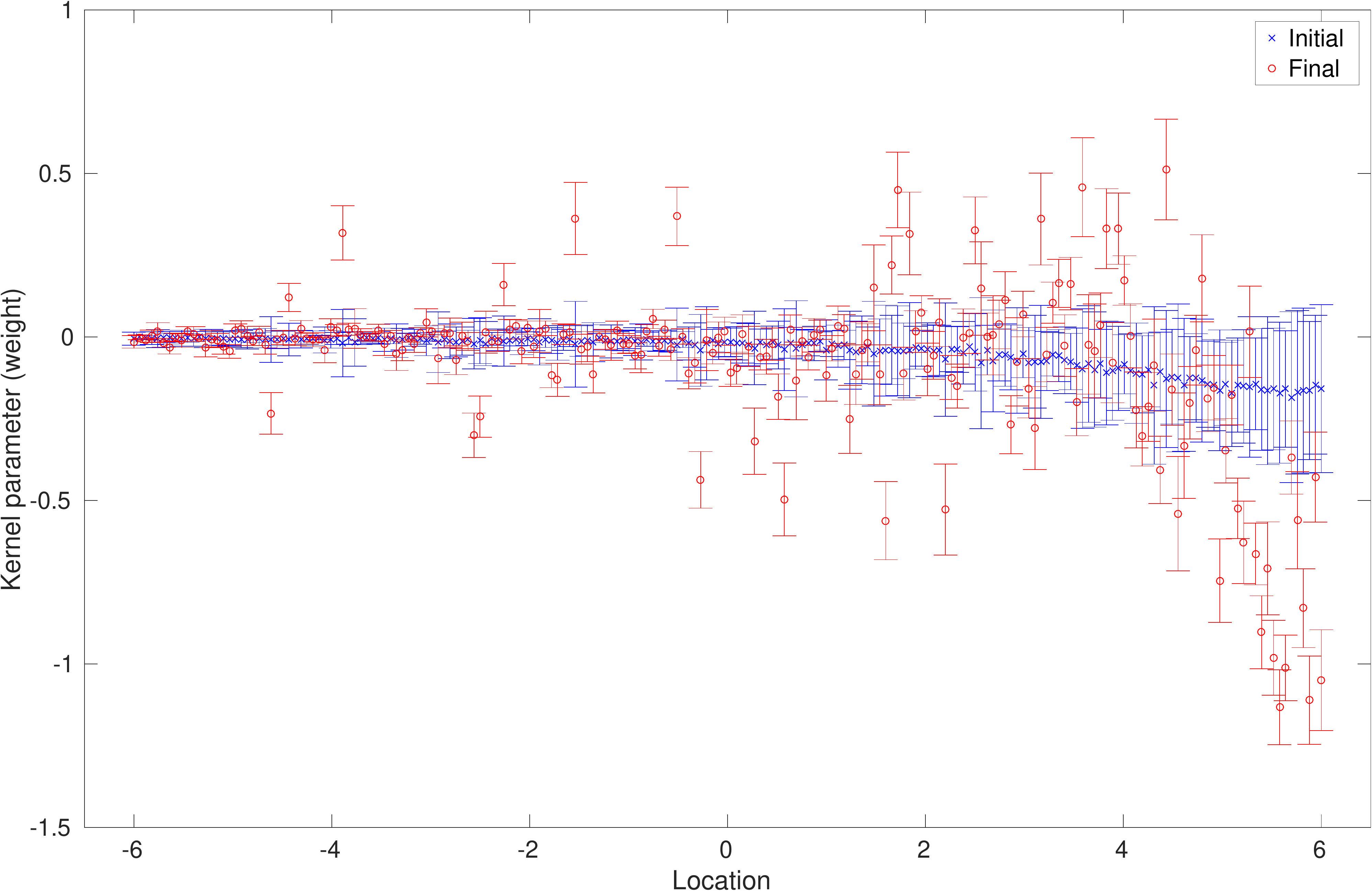}} \\
		\subfloat[Scale parameters ($\beta$) associated with C2]{\includegraphics[width=0.45\textwidth]{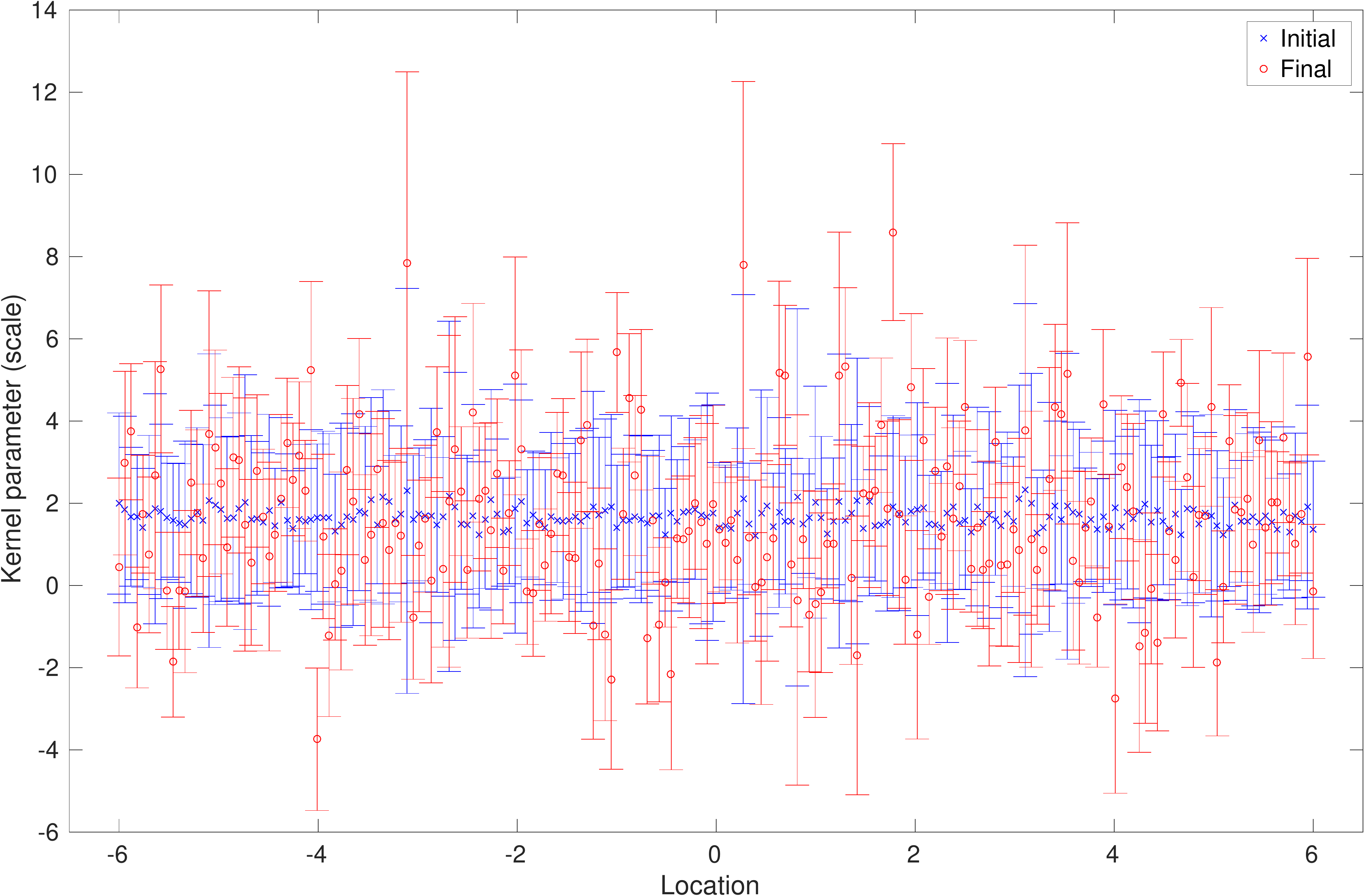}} &
		\subfloat[Weight parameters ($c$) associated with C2]{\includegraphics[width=0.45\textwidth]{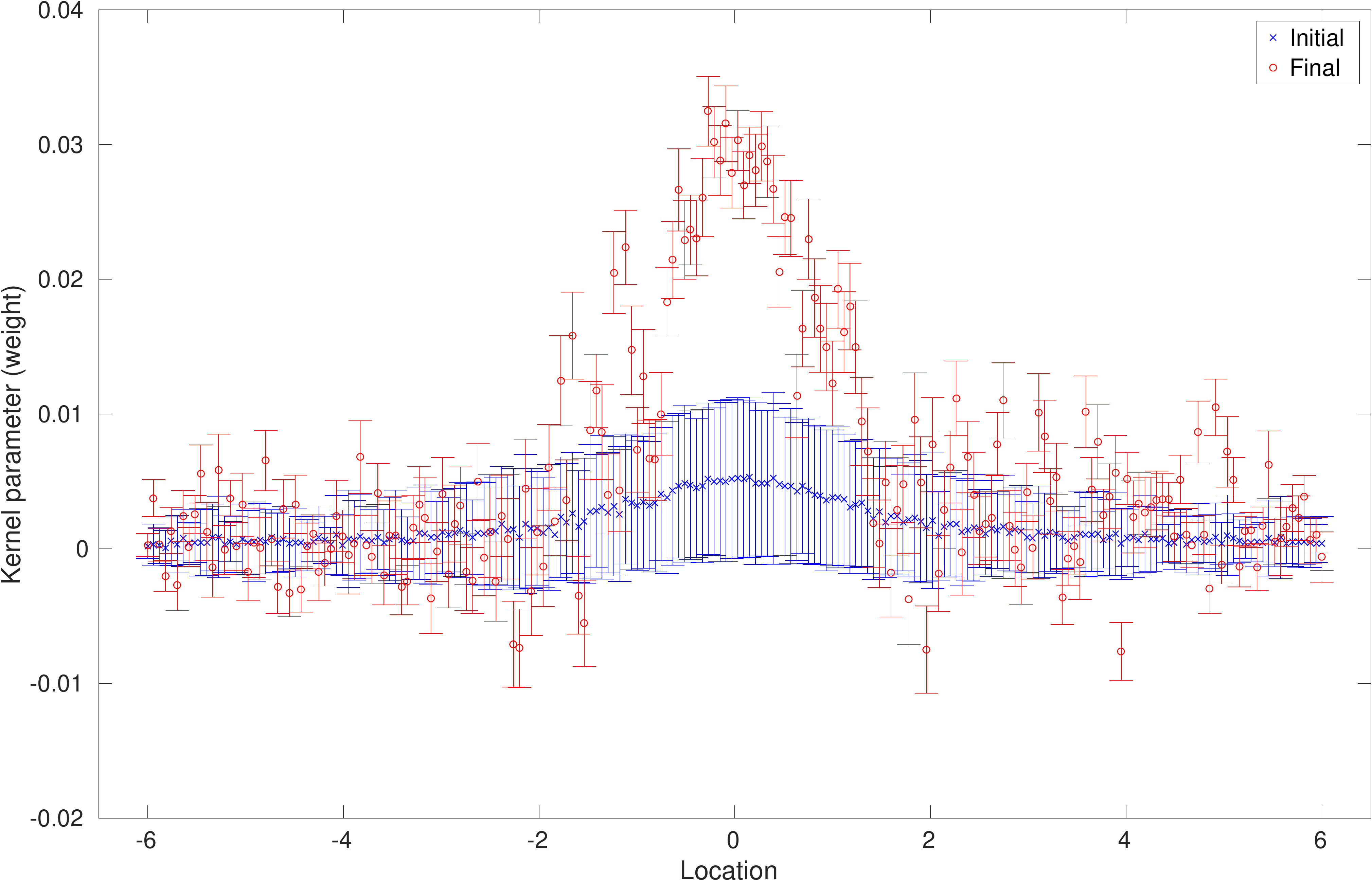}} \\
		\subfloat[Scale parameters ($\beta$) associated with C3]{\includegraphics[width=0.45\textwidth]{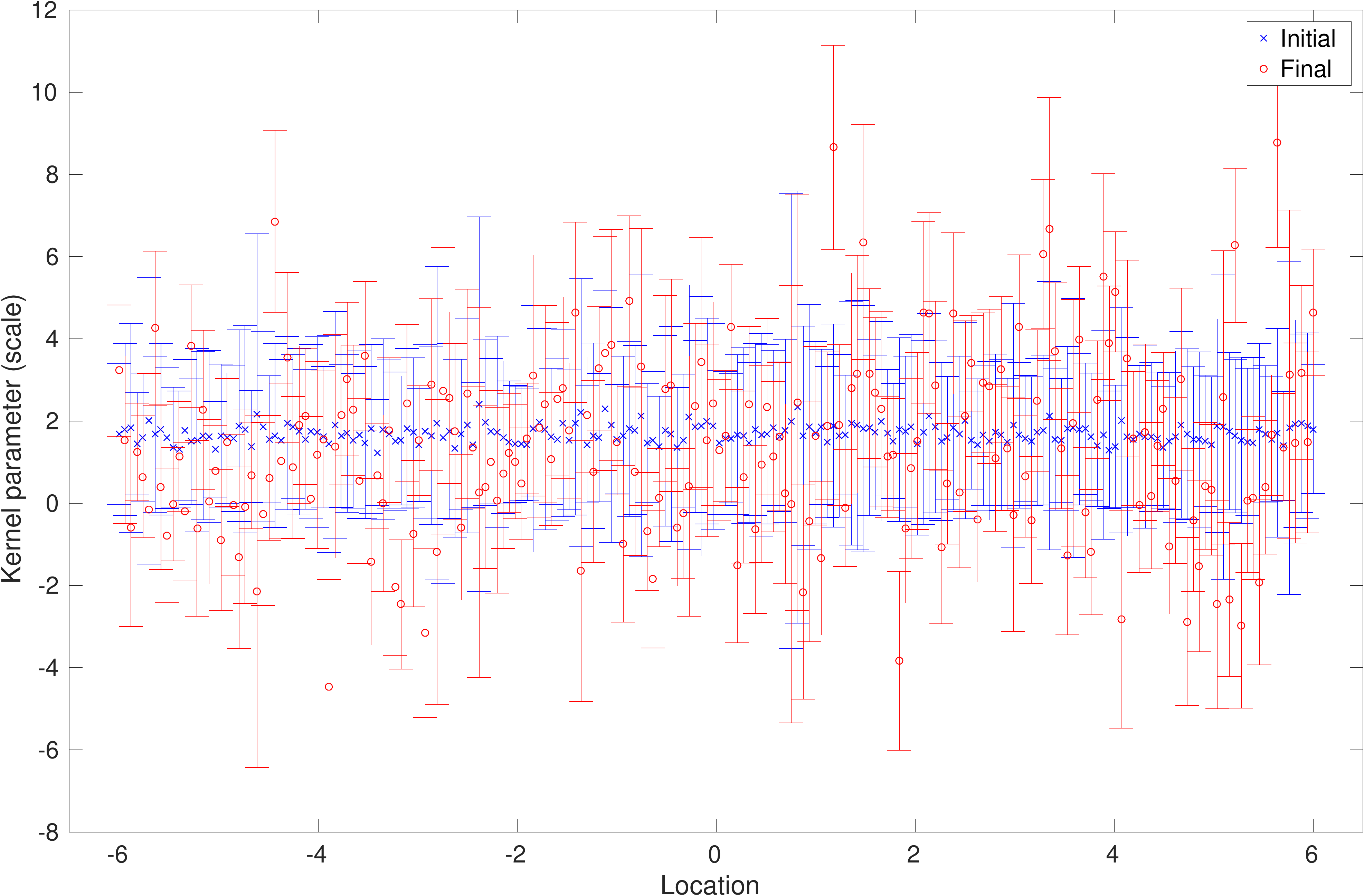}} &
		\subfloat[Weight parameters ($c$) associated with C3]{\includegraphics[width=0.45\textwidth]{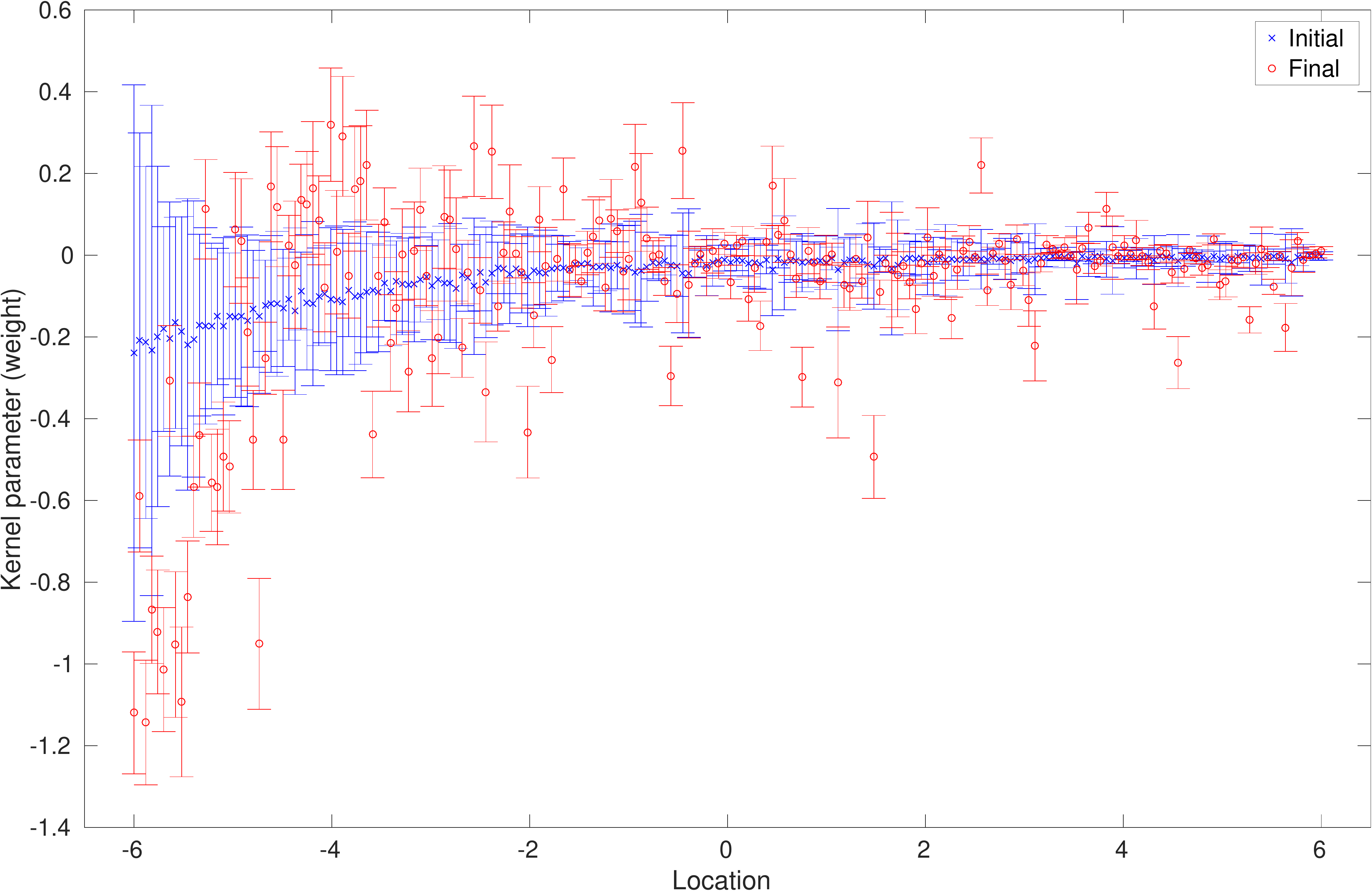}} \\
	\end{tabular}
	\caption{\label{fig:errorbar_kernelPara_multimodal} Similar to Figure \ref{fig:errorbar_kernelPara_unimodal}, but for multi-modal training inputs. For visualization, we plot scale (left column) and weight (right column) parameters associated with different clusters separately.}
\end{figure}

\renewcommand{\nScale}{0.2}
\begin{figure} 
	\centering
	\begin{tabular}{cc}
		\subfloat[Initial ensemble of predictions]{\includegraphics[width=0.45\textwidth]{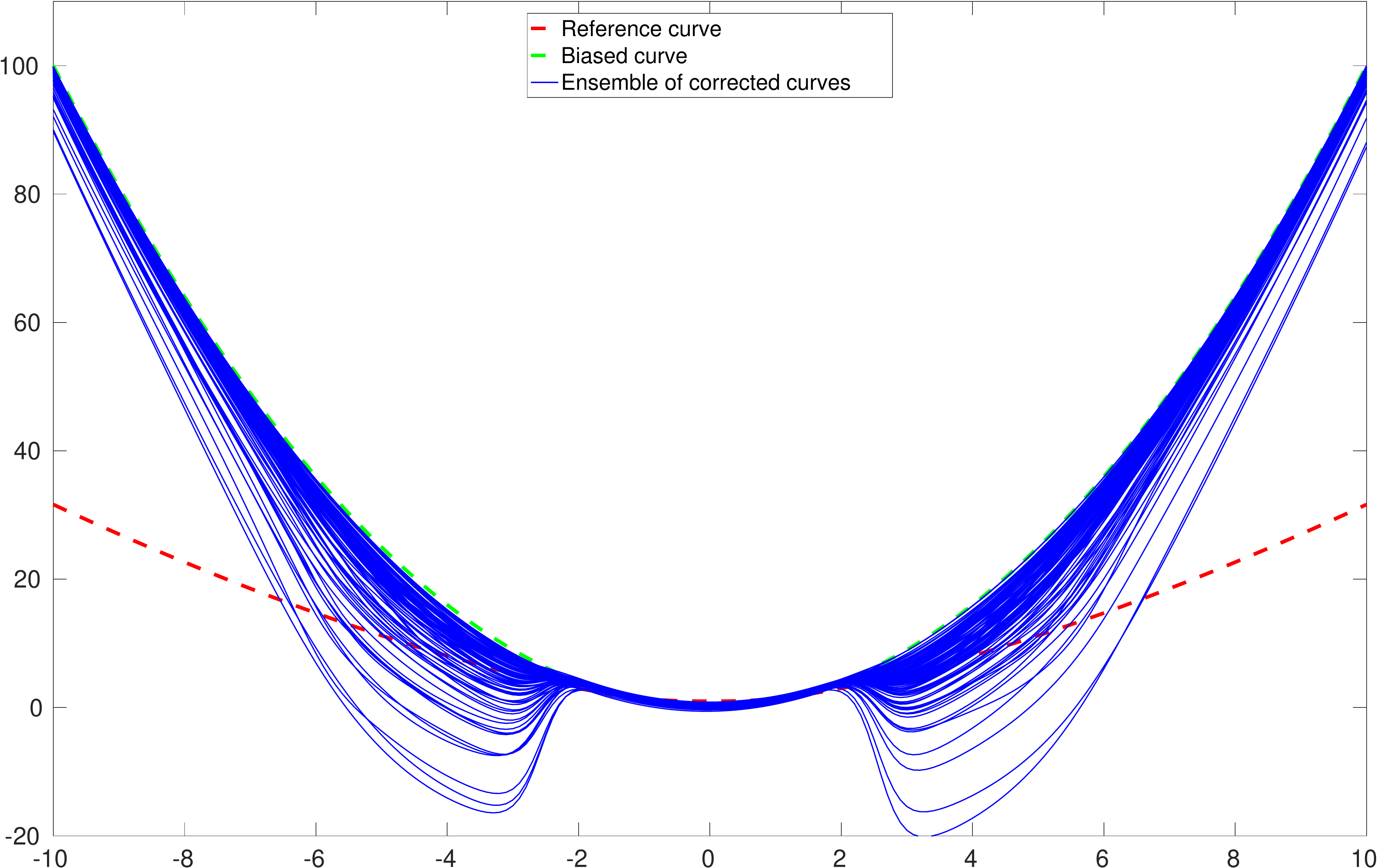}} &
		\subfloat[Initial mean predictions]{\includegraphics[width=0.45\textwidth]{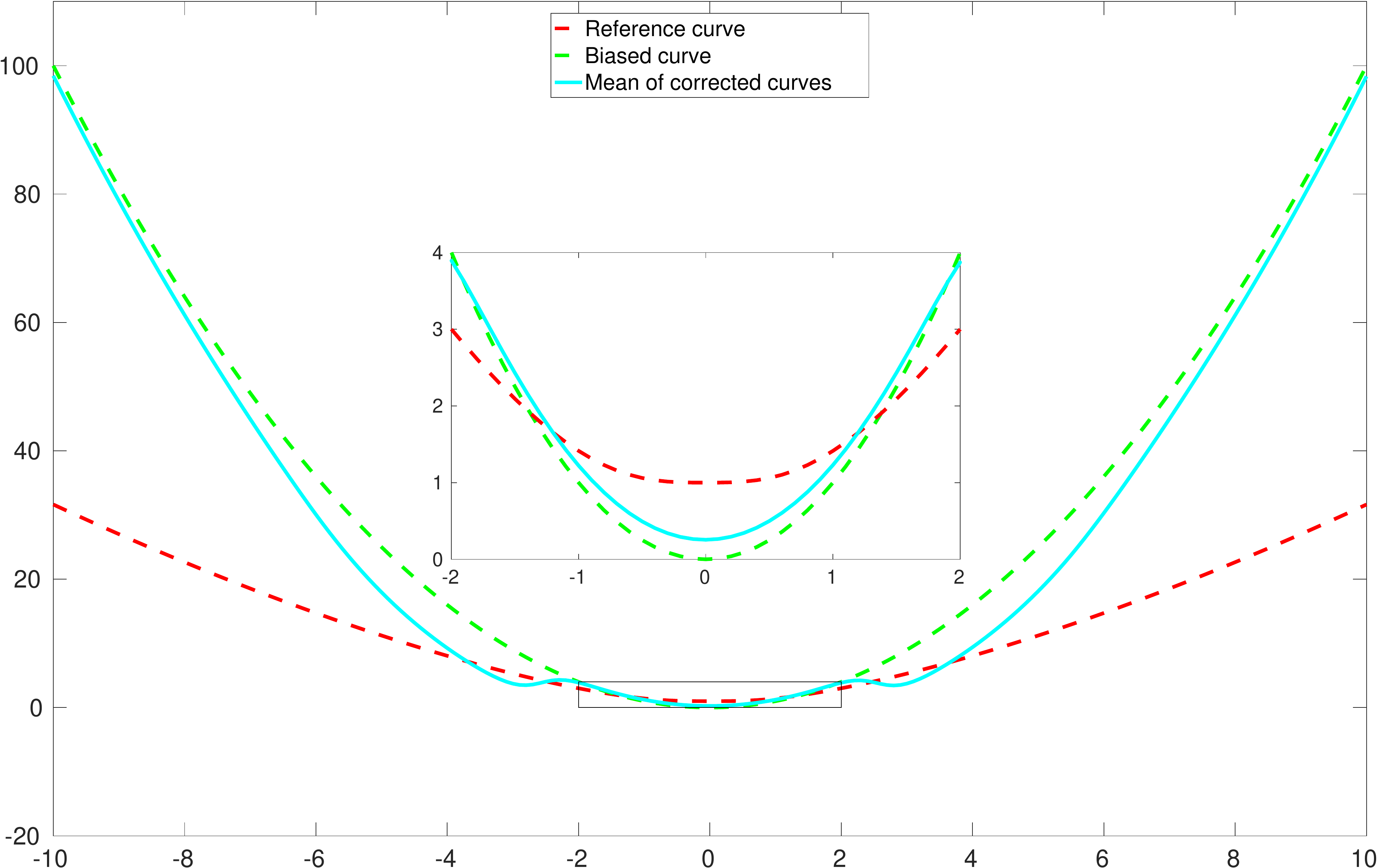}} \\
		\subfloat[Ensemble of predictions (after learning C1)]{\includegraphics[width=0.45\textwidth]{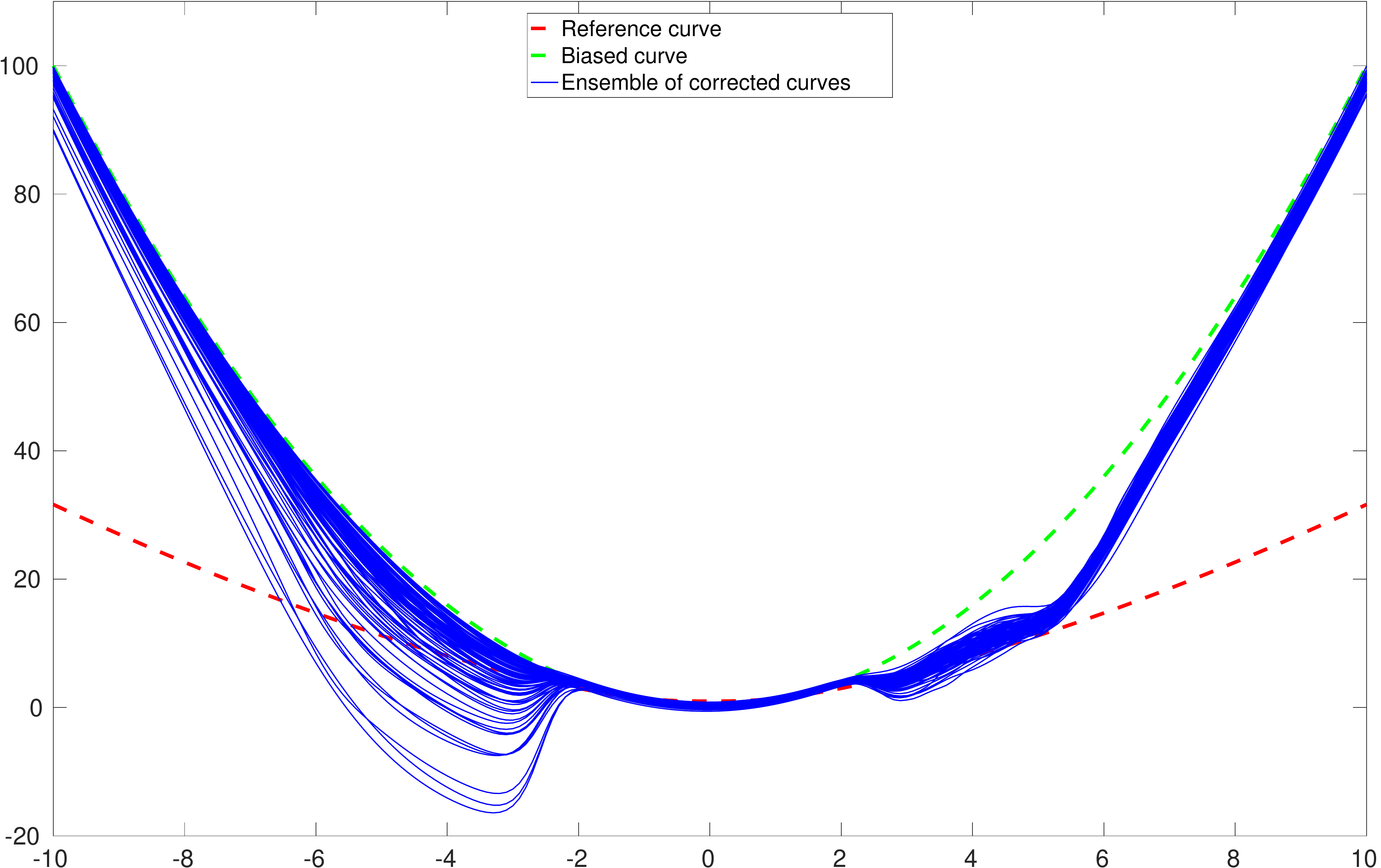}} &
		\subfloat[Mean predictions (after learning C1)]{\includegraphics[width=0.45\textwidth]{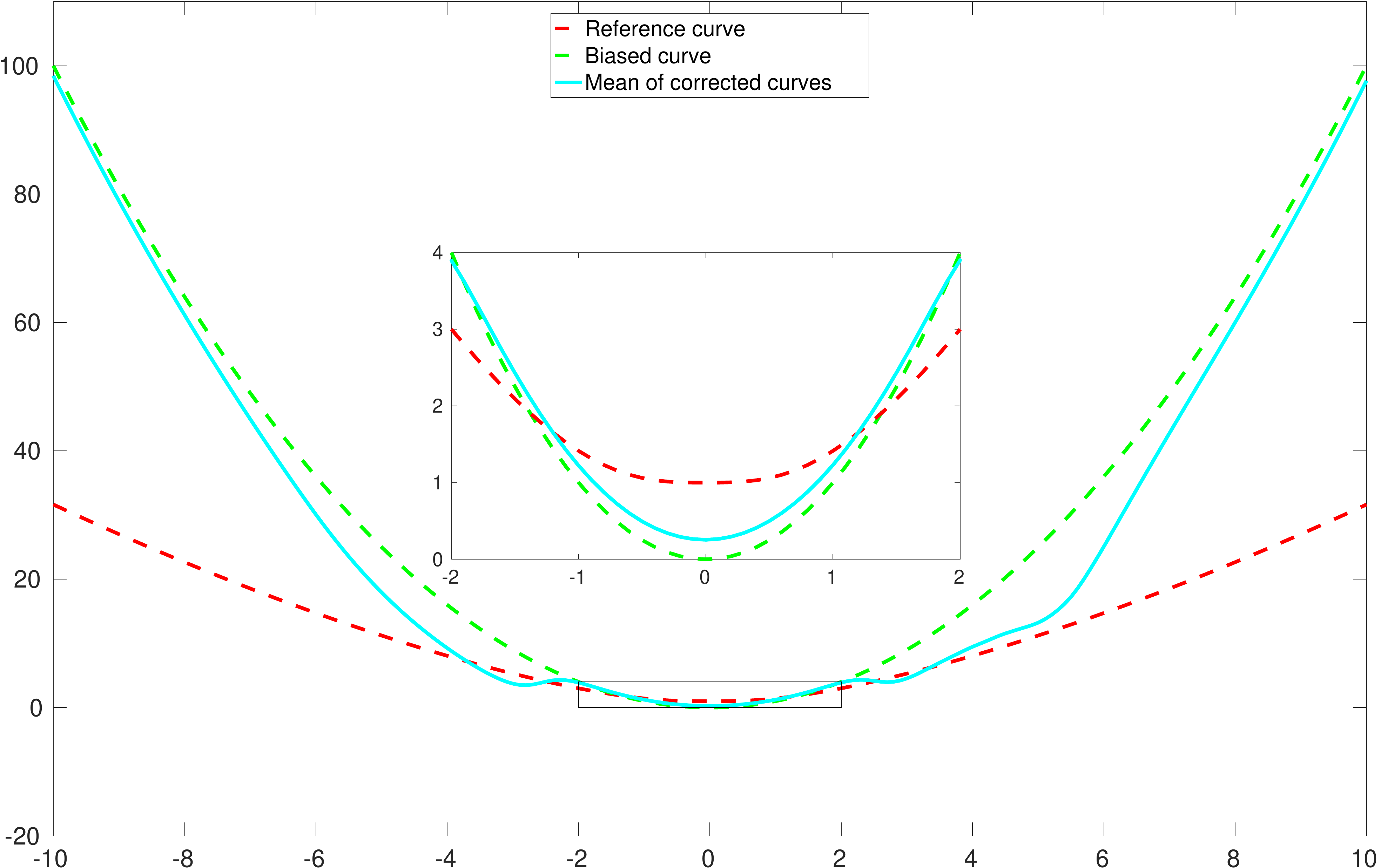}} \\
		\subfloat[Ensemble of predictions (after learning C1 and C2)]{\includegraphics[width=0.45\textwidth]{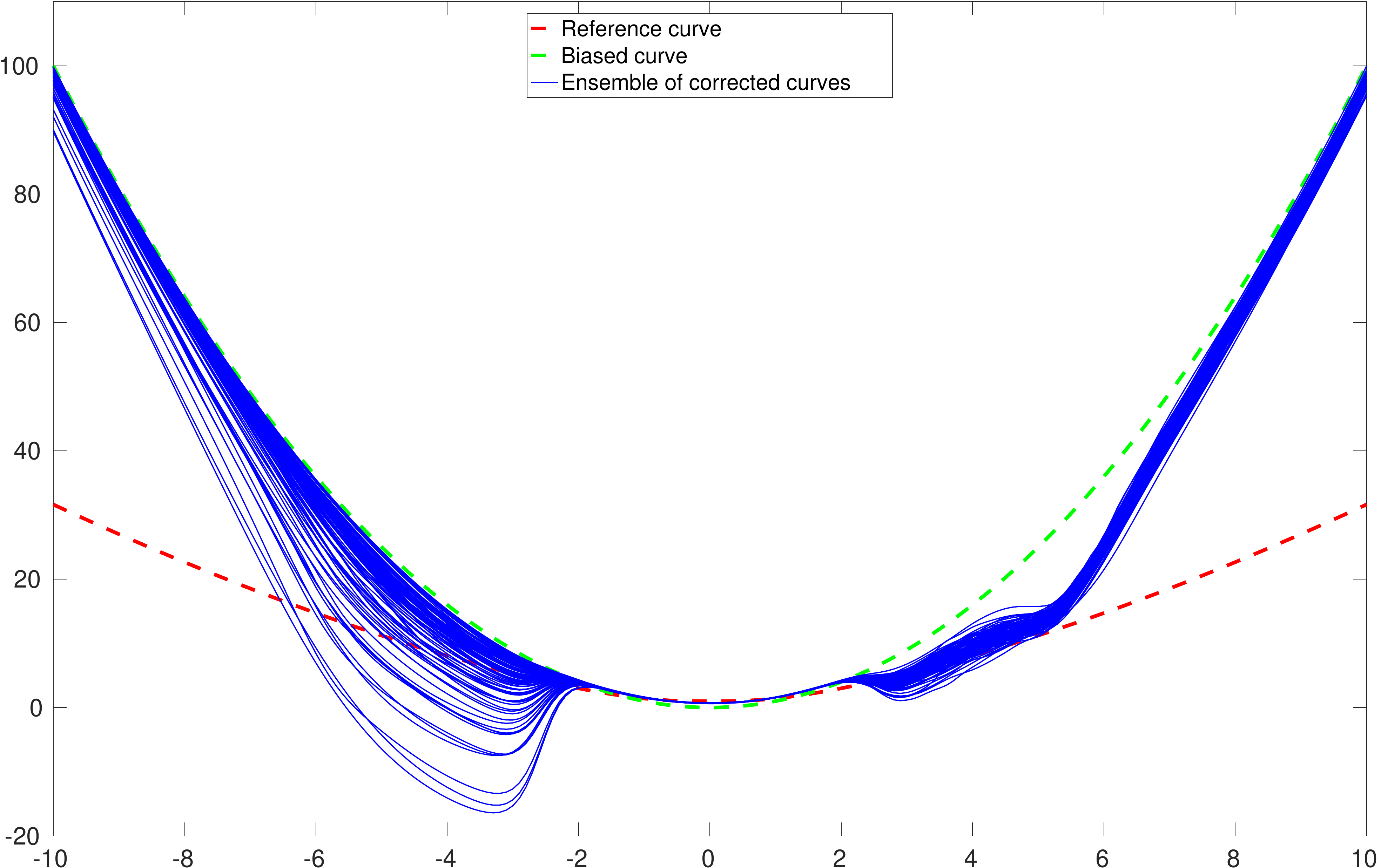}} &
		\subfloat[Mean predictions (after learning C1 and C2)]{\includegraphics[width=0.45\textwidth]{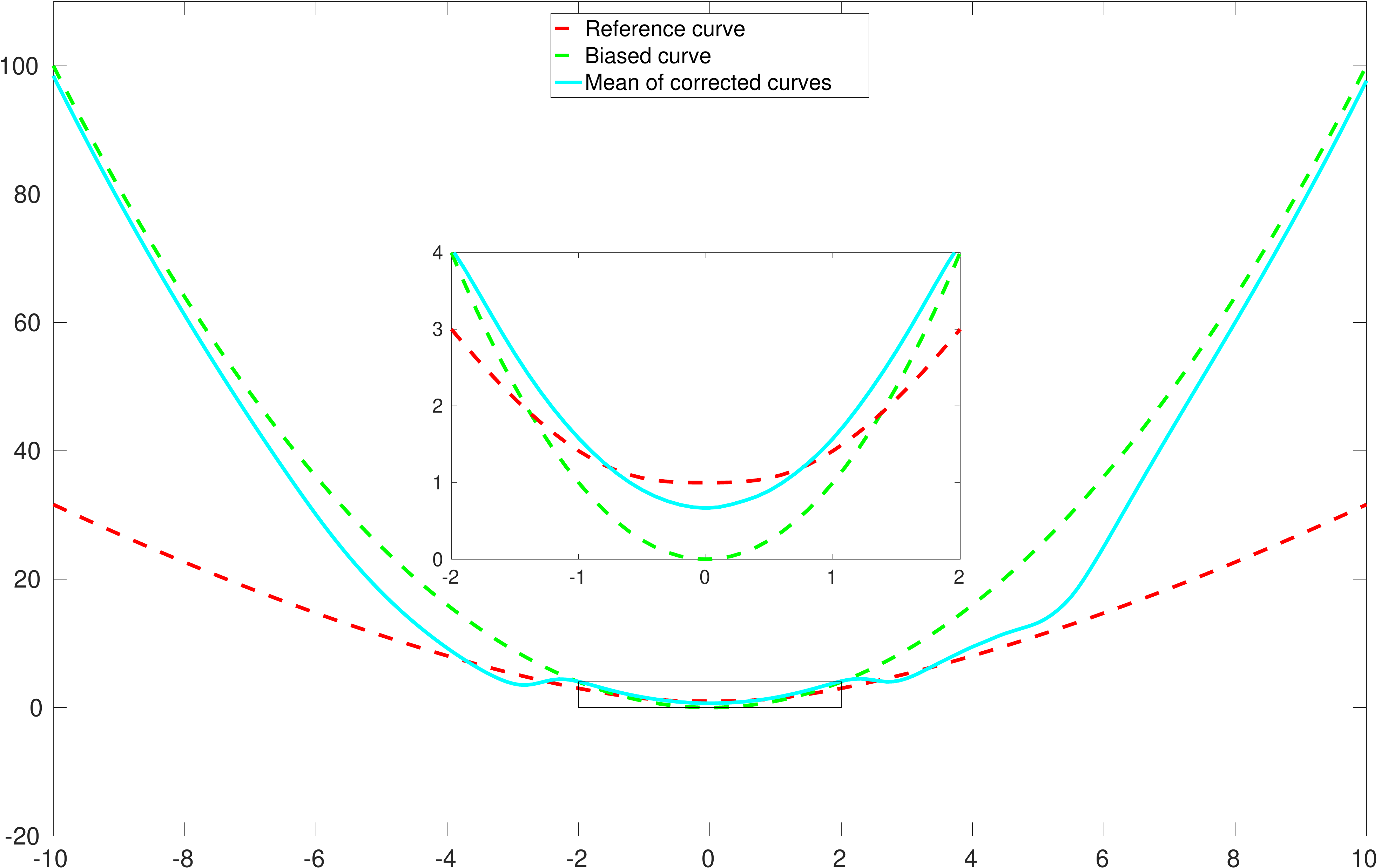}} \\
		\subfloat[Ensemble of predictions (after learning C1, C2 and C3)]{\includegraphics[width=0.45\textwidth]{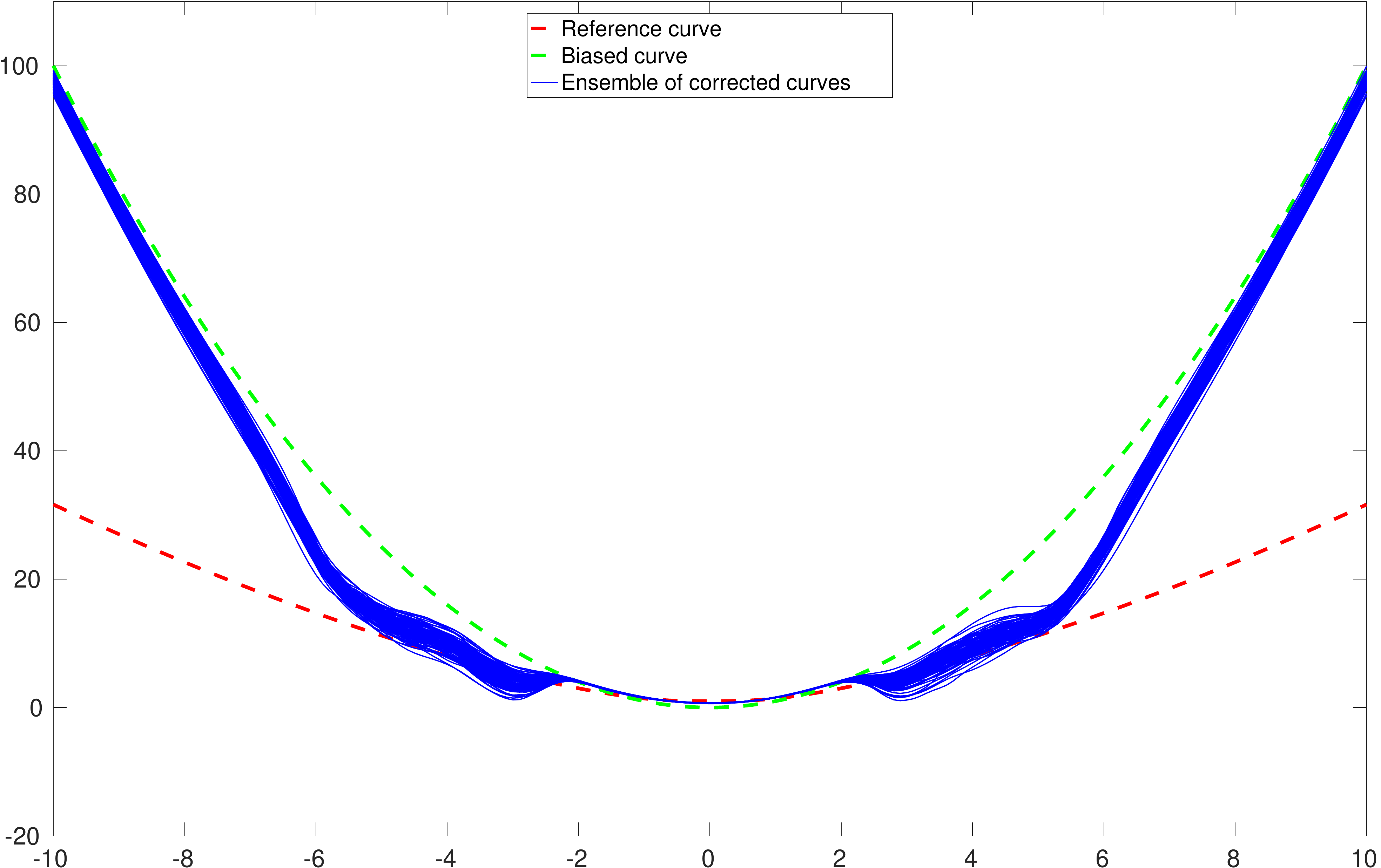}} &
		\subfloat[Mean predictions (after learning C1, C2 and C3)]{\includegraphics[width=0.45\textwidth]{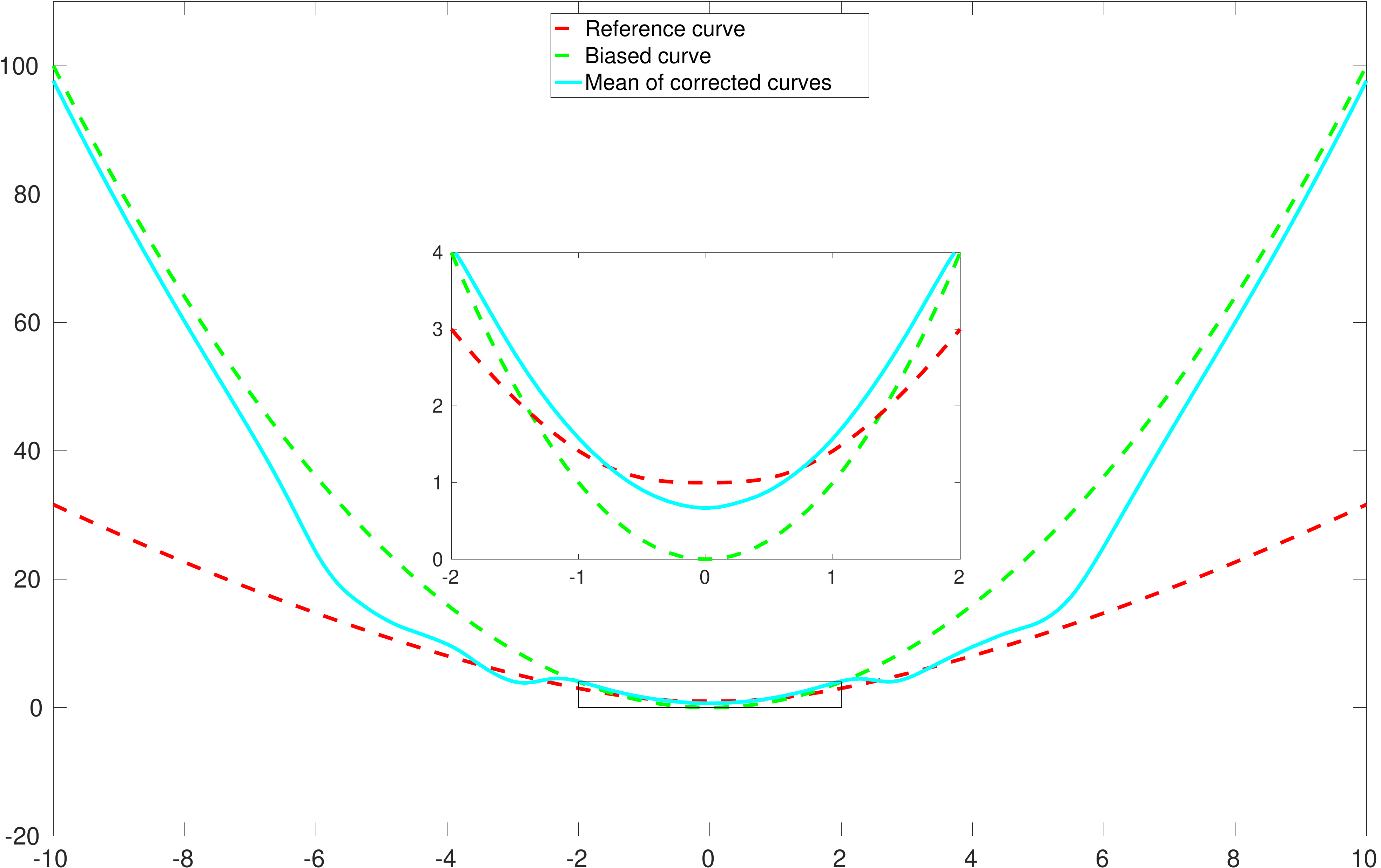}} \\
	\end{tabular}
	\caption{\label{fig:clustering_prediction_ensemlbe} Similar to Figure \ref{fig:multimodal_prediction_ensemlbe_cluster1}, but for the case in which the multi-modal learning strategy (MMLS) is adopted. In the experiment, the number $N_{cl}$ of clusters is $3$, the same as the number of modes in the training inputs. Note that the learning process is carried out cluster by cluster.}
\end{figure} 

%
In the first experiment, we investigate the case where the number of clusters is the same as the number of modes in the training inputs, i.e., $N_{cl} = 3$. We use the MATLAB$^\copyright$ function ``fitgmdist'' to estimate the parameters like weight ($w$), mean ($\mu$) and variance ($\sigma^2$) (cf. Eq. (\ref{eq:GMM})) associated with each Gaussian component. Table \ref{tab:GMM_components} summarizes the number of training data points, as well as the values of the aforementioned parameters, associated with each component (cluster). This indicates that the GMM is fitted quite well, in light of how these data are generated. Table \ref{tab:GMM_components} also provides each Gaussian component a label (e.g., ``C1''), which will be adopted in the discussions below.

With the aforementioned settings, in principle one can update the kernel parameters associated with each cluster in parallel, although in the current work, such updates are conducted in a sequential manner, namely, C1 $\rightarrow$ C2 $\rightarrow$ C3. Figure \ref{fig:boxplot_dataMismatch_multimodal} shows the box plots of data mismatch, with respect to training (left column) and CV (right column) datasets, respectively. For the training dataset, we report data mismatch at different iteration steps cluster by cluster. For instance, in Panel (a) of Figure \ref{fig:boxplot_dataMismatch_multimodal}, data mismatch is calculated using the differences between the training outputs in C1, and the predicted outputs with respect to the training inputs in C1. Under this setting, one can see that the ensemble-based learning algorithm progressively reduces data mismatch within each cluster.        

For the CV dataset, we do not pre-cluster the data points into different clusters. To compute data mismatch with respect to the CV dataset, we use all the CV data points ($6000$ in total). For better comprehension, Eqs. (\ref{eq:GMM_probability}) and (\ref{eq:predicted_residual}) are referred in our discussion below. Given a CV input $x^{cv}$ and an ensemble of kernel parameters $\boldsymbol{\theta}_{j,s}$ ($j=1,2,\dotsb, N_e$) for a certain cluster $s$, we compute an ensemble of predicted outputs $\hat{h}(x^{cv}; \boldsymbol{\theta}_{j,s})$ (cf. Eq. (\ref{eq:predicted_residual})), as well as a probability $P_{s}(x^{cv})$ with respect to the GMM (cf. Eq. (\ref{eq:GMM_probability})). Data mismatch with respect to the cluster $s$ is then calculated using the differences between the CV output weighted by $P_{s}(x^{cv})$, and the predicted output $\hat{h}(x^{cv}; \boldsymbol{\theta}_{j,s})$ that is also weighted by $P_{s}(x^{cv})$. In this way, we are able to cross-validate the impacts of supervised learning within individual clusters. As reported in Panels (b), (d) and (f) of Figure \ref{fig:boxplot_dataMismatch_multimodal}, data mismatch of the CV dataset with respect to all clusters tends to decrease through the iterations, indicating that the learning process goes reasonably well.   

Similar to Figure \ref{fig:errorbar_kernelPara_unimodal}, in Figure \ref{fig:errorbar_kernelPara_multimodal} we also plot the initial (in blue) and final (in red) ensembles of scale (left column) and weight (right column) parameters associated with different clusters. For scale parameters, compared to the case with unimodal training inputs (cf. Figure \ref{fig:errorbar_kernelPara_unimodal}(a)), there appear to be more substantial differences between initial and final values in all three clusters (cf. Figures \ref{fig:errorbar_kernelPara_multimodal}(a), \ref{fig:errorbar_kernelPara_multimodal}(c) and \ref{fig:errorbar_kernelPara_multimodal}(e)). For weight parameters, similar to the case with unimodal training inputs (cf. Figure \ref{fig:errorbar_kernelPara_unimodal}(b)), significant changes from initial to final values can also be spotted in the areas surrounding the mode of each Gaussian component. 

Similar to Figure \ref{fig:multimodal_prediction_ensemlbe_cluster1}, Figure \ref{fig:clustering_prediction_ensemlbe} shows the results after the MMLS is adopted to train kernel parameters cluster by cluster. As one can see, with the MMLS, the initial ensemble in Figure \ref{fig:clustering_prediction_ensemlbe}(a) exhibits multimodality, which is not the case for the initial ensemble in Figure \ref{fig:multimodal_prediction_ensemlbe_cluster1}(a), where the MMLS is not employed. On top of the multi-modal initial ensemble, the ensemble-based training algorithm in general tends to improve the predictions, by updating the kernel parameters sequentially through the use of training data in individual clusters.

\subsubsection*{The impact of the number $N_{cl}$ of clusters}

%
\renewcommand{\nScale}{0.2}
\begin{figure} 
	\centering
	\begin{tabular}{cc}
		\subfloat[Final ensemble of predictions with $N_{cl} = 2$]{\includegraphics[width=0.45\textwidth]{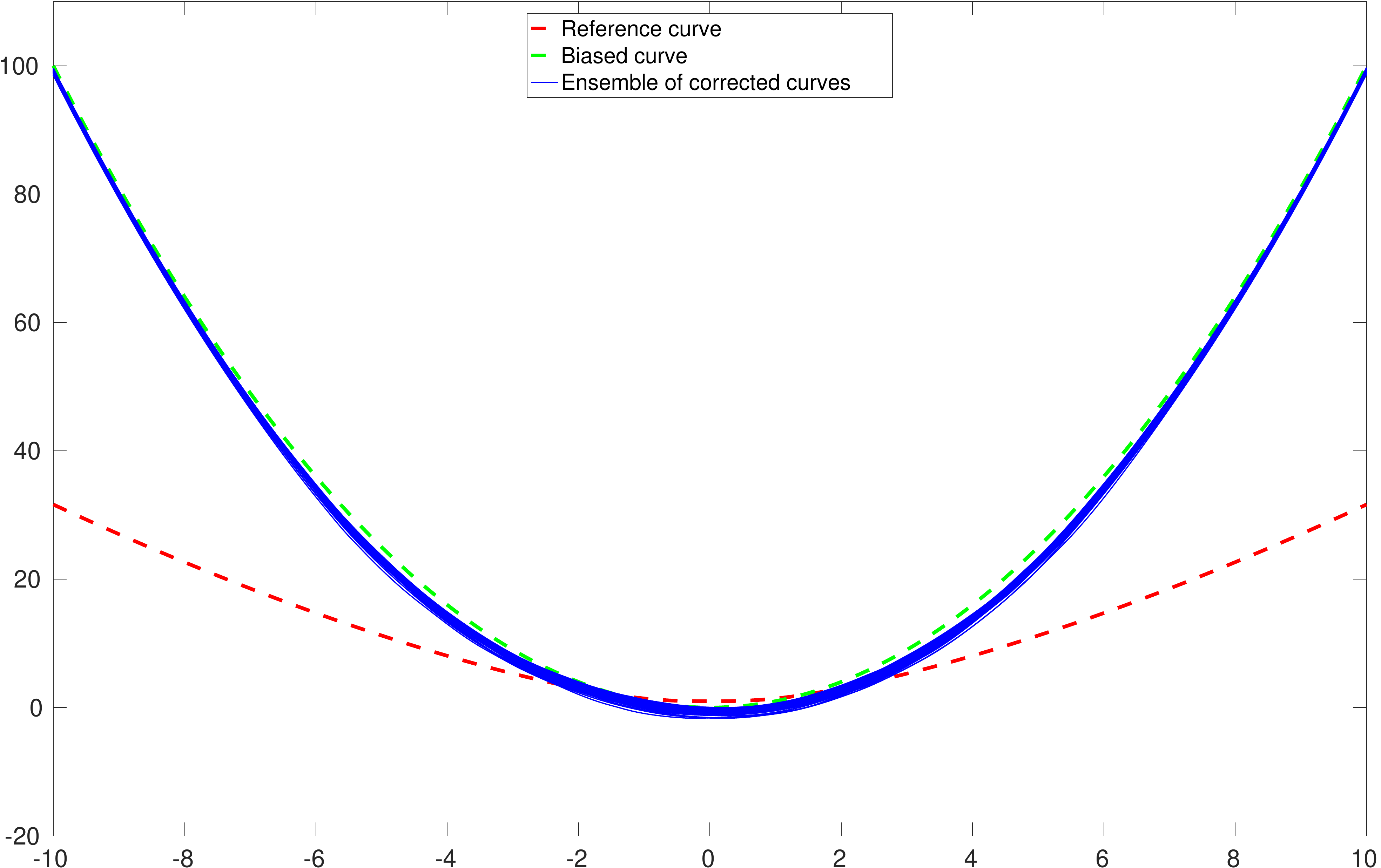}} &
		\subfloat[Final mean predictions  with $N_{cl} = 2$]{\includegraphics[width=0.45\textwidth]{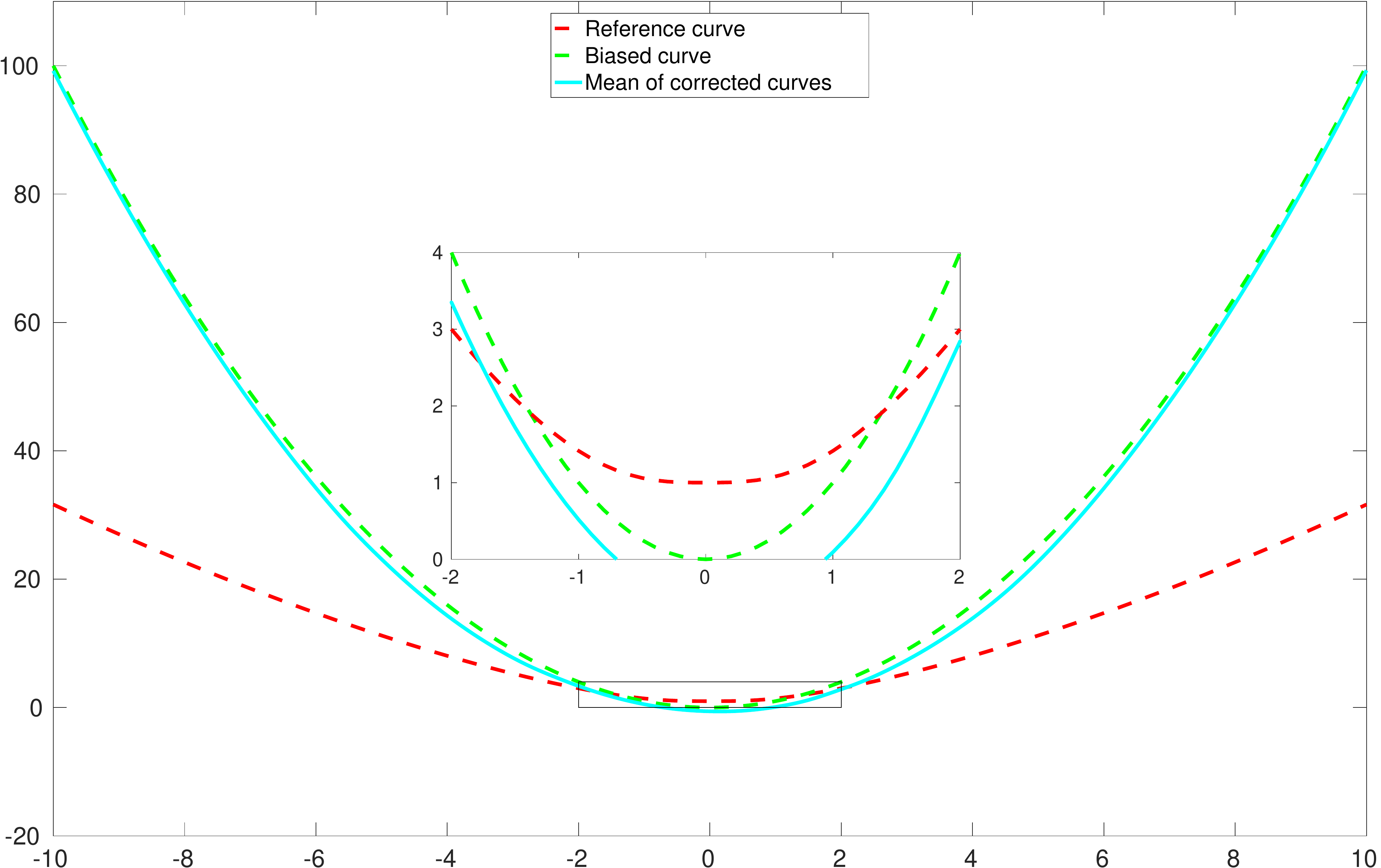}} \\
		\subfloat[Final ensemble of predictions with $N_{cl} = 4$]{\includegraphics[width=0.45\textwidth]{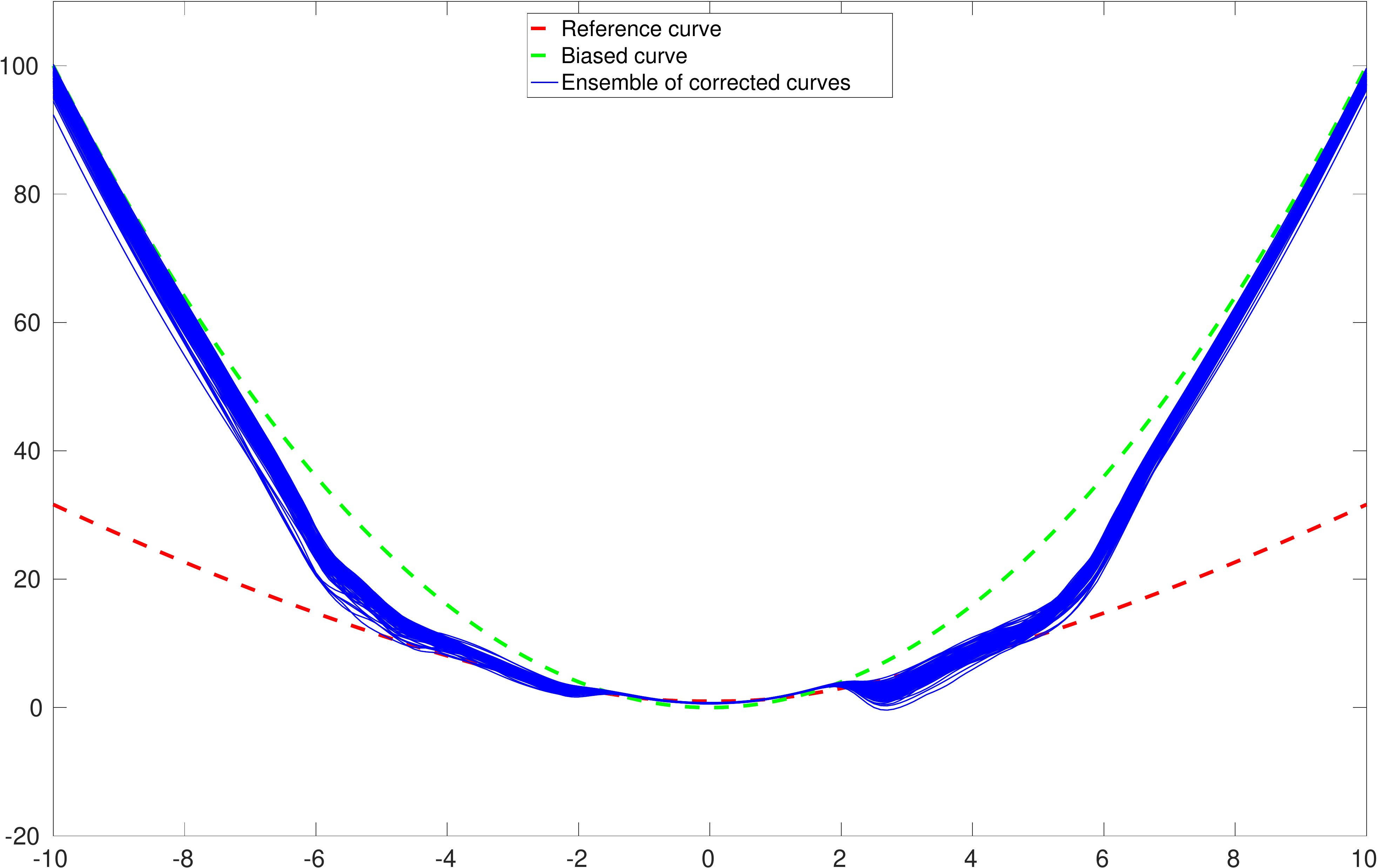}} &
		\subfloat[Final mean predictions  with $N_{cl} = 4$]{\includegraphics[width=0.45\textwidth]{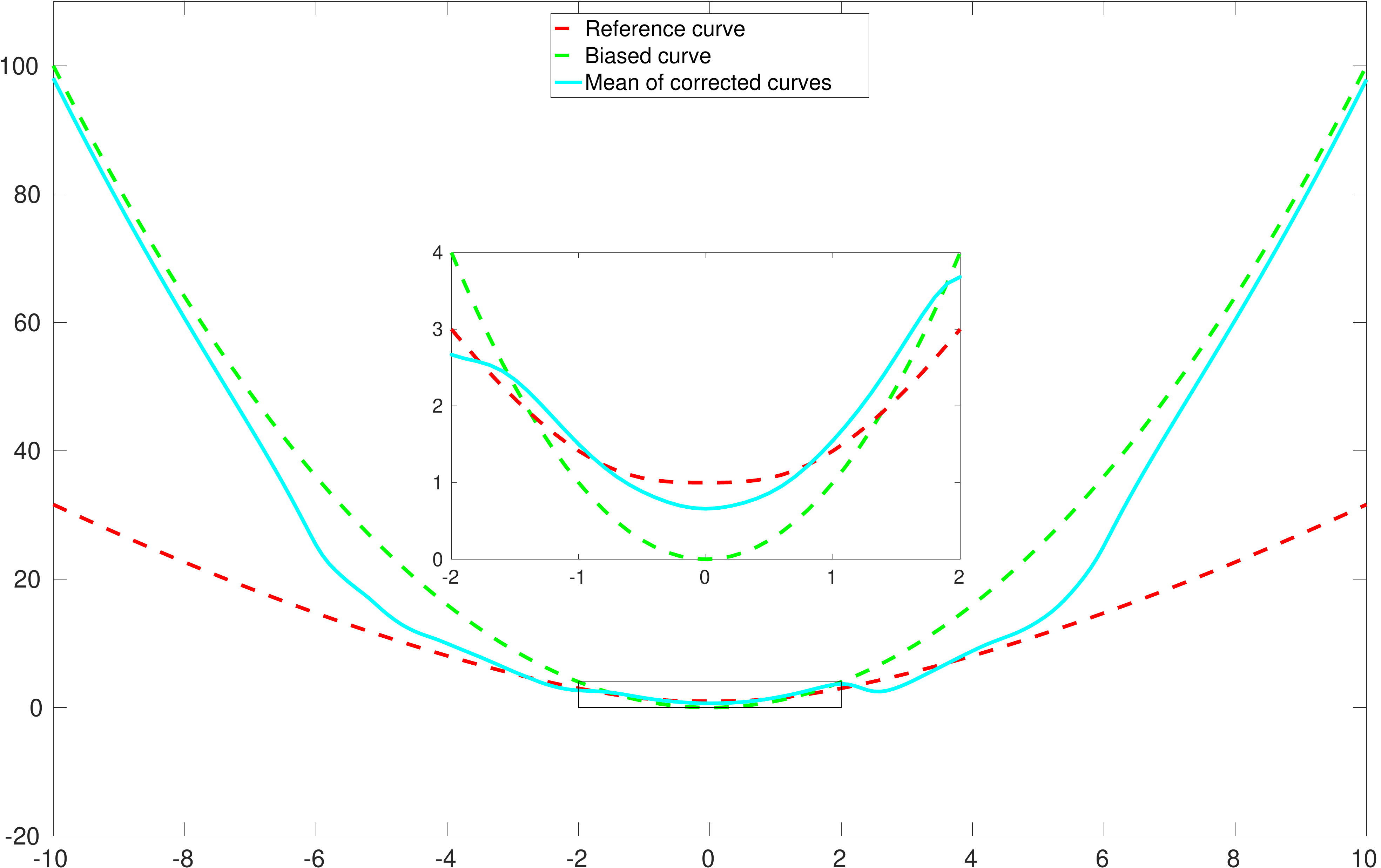}} \\
		\subfloat[Final ensemble of predictions with $N_{cl} = 6$]{\includegraphics[width=0.45\textwidth]{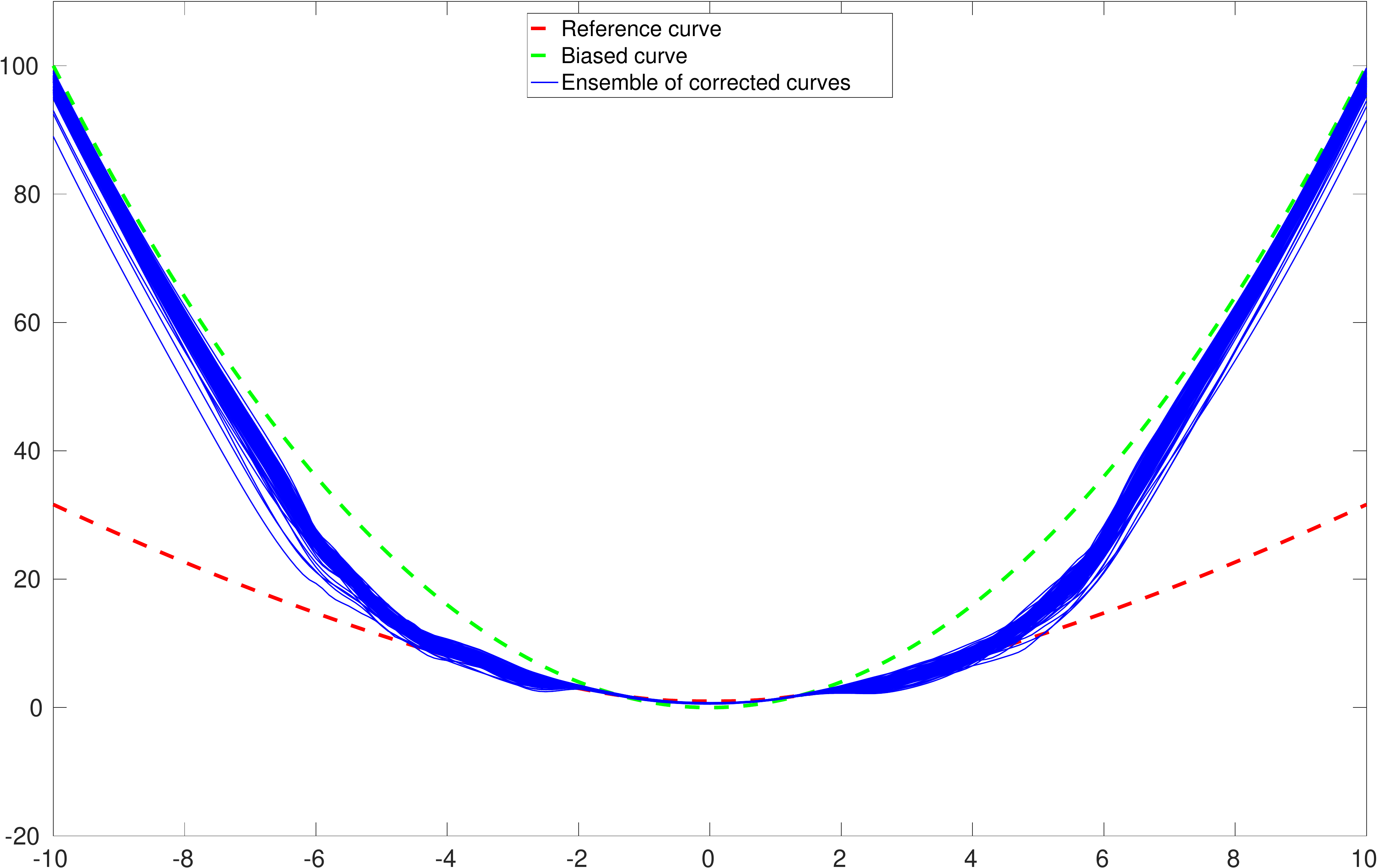}} &
		\subfloat[Final mean predictions  with $N_{cl} = 6$]{\includegraphics[width=0.45\textwidth]{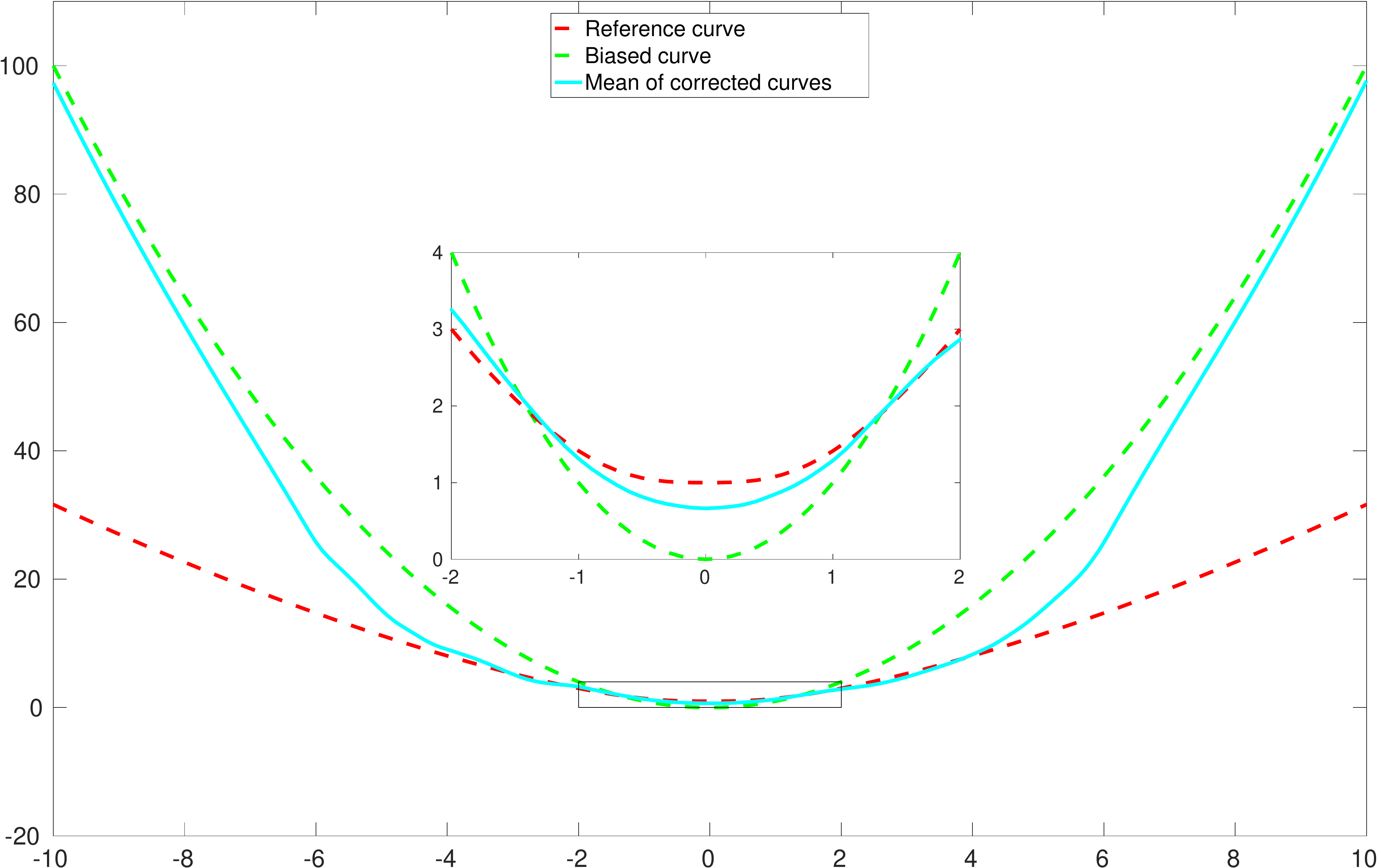}} \\
		\subfloat[Final ensemble of predictions with $N_{cl} = 8$]{\includegraphics[width=0.45\textwidth]{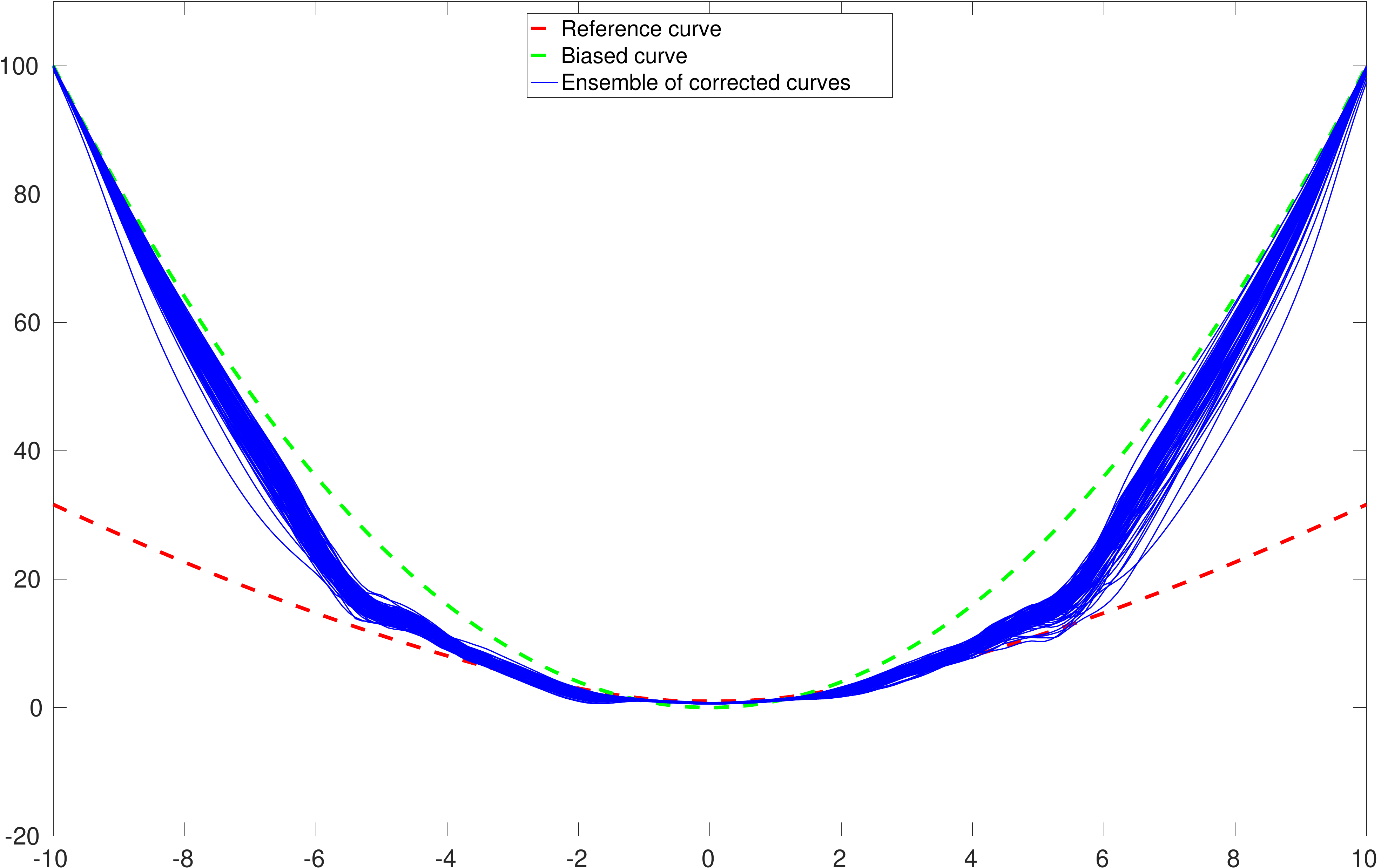}} &
		\subfloat[Final mean predictions  with $N_{cl} = 8$]{\includegraphics[width=0.45\textwidth]{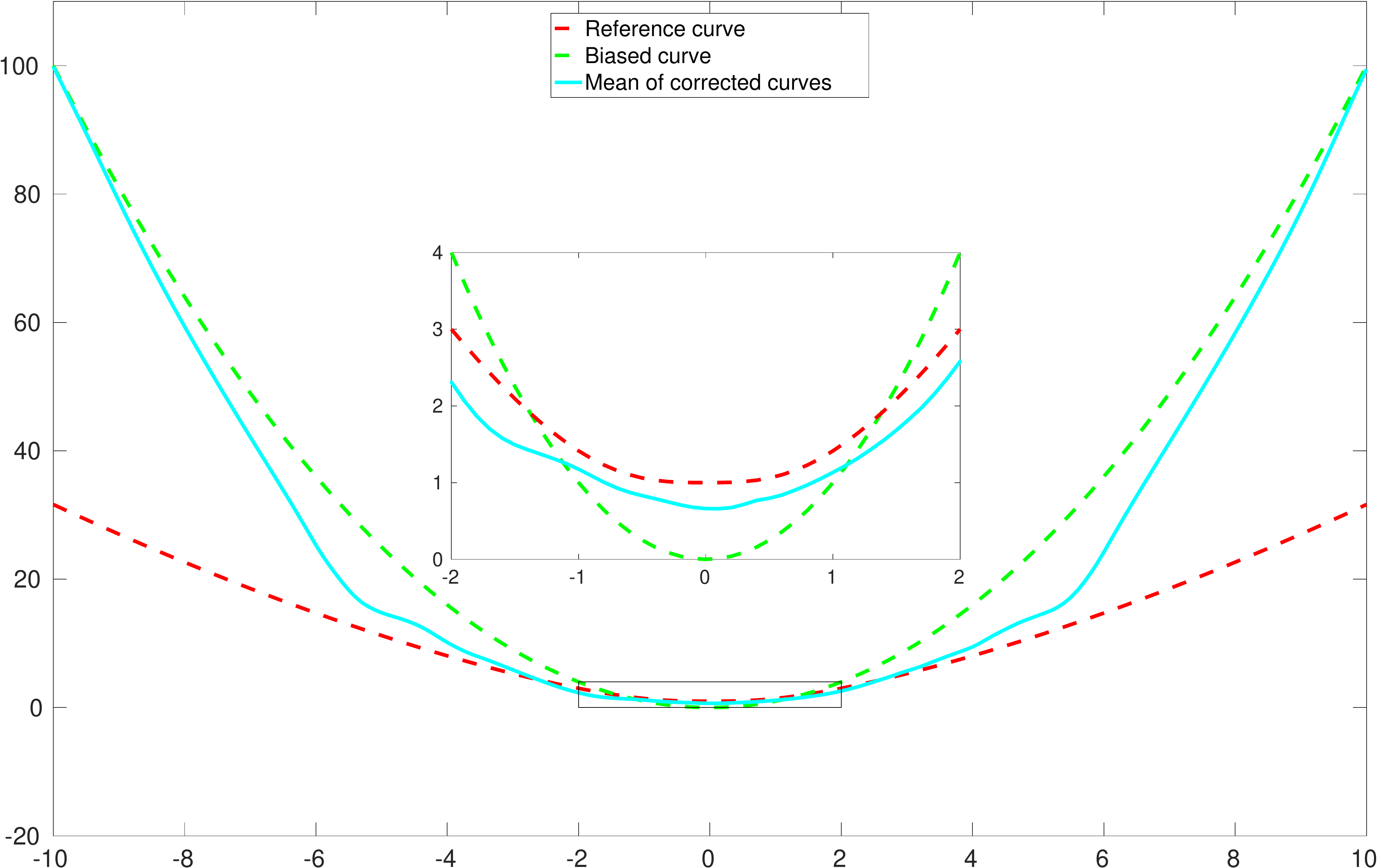}} \\
	\end{tabular}
	\caption{\label{fig:prediction_ensemlbe_multicluster} Similar to Figure \ref{fig:clustering_prediction_ensemlbe}, but for the final prediction results after all the training data in different clusters are used to learn kernel parameters. Presented here are the results with respect to of the choices of using $2$, $4$, $6$ and $8$ clusters to fit the GMM (from top to bottom), respectively.}
\end{figure}

%

The previous results indicate that, when the MMLS is adopted and the number of clusters is the same as the number of modes in the training inputs, one can improve the performance of predictions using the learned kernel parameters. Here, we also examine what will happen, when the MMLS is adopted, but the number of clusters is not necessarily the same as the number of modes in the training inputs.

Figure \ref{fig:prediction_ensemlbe_multicluster} presents some of the final prediction results (after using all training data to learn kernel parameters) from an experiment, in which we adopt different numbers $N_{cl}$ of clusters (e.g., $2$, $4$, $6$ and $8$) to fit the GMM using the same training inputs (with $3$ modes) as in the previous experiment. Combining the results in Figures \ref{fig:multimodal_prediction_ensemlbe_cluster1}, \ref{fig:clustering_prediction_ensemlbe} and \ref{fig:prediction_ensemlbe_multicluster}, it appears that, if the number of the clusters is less than the number of modes in the training inputs, then the learned ensemble of models tends to have insufficient capacities to perform relatively well in the prediction tests. However, when the number of the clusters becomes no less than the number of modes, then the capacities of the learned models tend to improve. In this particular case, it seems that, if $N_{cl}$ is slightly larger than the number of modes (e.g., $N_{cl} = 6$), then one might actually achieve better prediction accuracies over certain intervals, in comparison to the choice of $N_{cl} = 3$. Of course, given a fixed number of training data, on average the number of training data per cluster will reduce as $N_{cl}$ increases. Therefore, if $N_{cl}$ becomes too large (e.g., $N_{cl} = 8$), the prediction accuracies may be instead worsened as the number of training data within each cluster decreases. This insight will be useful for us to handle data assimilation problems in the presence of forward-simulator imperfection, yielding improved flexibility and assimilation performance, as will be shown in the next section.              

\section*{Numerical results in a data assimilation problem with an imperfect forward simulator}\label{sec:results_DA} 
The preceding section indicates that, when combined with the MMLS, the ensemble-based kernel learning algorithm performs reasonably well in the presented SLP. As discussed previously, the idea of kernel-based functional approximation can also be extended to handle data assimilation problems with imperfect forward simulators. As a proof-of-concept study, in what follows, we illustrate the performance of the integrated data assimilation (history matching) framework, Eqs. (\ref{eq:HM_with_SLP}) through (\ref{eq:new_forward_simulator}), in a synthetic 2D problem. In the experiment, we have a reference model in the dimension of $100 \times 120$ (cf. Figure \ref{fig:DA_models_perfect}(a)). The corresponding (noisy) observations (cf. Figure \ref{fig:DA_obs_perfect}(a)) are generated by first applying a function $f(z) = \left(\vert z \vert^3 + 1 \right)^{1/2}$ to each gridblock of the reference model, and then adding $10\%$ Gaussian white noise (relative to magnitudes) to the simulation outputs. As a result, in data assimilation, we have a set of observations distributing over the same gridblocks as in the reference model. 

The reference model in Figure \ref{fig:DA_models_perfect}(a) is generated through a fast Gaussian simulation method \citep{lorentzen2017history,Luo2018CorrLoc_Norne}, as a realization of a 2D Gaussian random field with zero mean, and an anisotropic covariance model whose STD is 2, and whose length scales along x and y directions are 15 and 25 gridblocks, respectively. The initial ensemble (with $100$ members) is generated in a similar way, but using a slightly different covariance model, whose STD is 2.2, and whose length scales along x and y directions are 17 and 23 gridblocks, respectively. Figure \ref{fig:DA_models_perfect}(b) shows the mean of the initial ensemble.    

As mentioned earlier, to use kernel-based functional approximation in Eq. (\ref{eq:residual_rbf_approximation}), we need to specify a set
$\mathbf{Z}^{cp} \equiv \{ z_k^{cp} \}_{k=1}^{N_{cp}}$ of center points, and a corresponding set $\mathbf{D}^{o,cp} \equiv \{ d_k^{o,cp} \}_{k=1}^{N_{cp}}$ of observations associated with $\mathbf{Z}^{cp}$. In the experiments, we do not assume to have hard data to condition on. Instead, we construct $\mathbf{Z}^{cp}$ and $\mathbf{D}^{o,cp}$ as follows. We set $N_{cp} = 200$, and take $z_k^{cp}$ as the points that evenly span an interval $[z_l,z_u)$, where $z_l = z_{min} - 0.1|z_{min}|$ and $z_u = z_{max} + 0.1|z_{max}|$, with $z_{min}$ and $z_{max}$ being the minimum and maximum values of the initial ensemble of model variables, respectively. To choose $d_k^{o,cp}$, we first compute the mean $\hat{\mathbf{z}}_0$ of the initial ensemble, and treat $\hat{\mathbf{z}}_0$ as if it were the ground truth that generates the real observations $\mathbf{d}^{o}$. With this treatment, for each given $z_k^{cp}$, we find 20 variables in $\hat{\mathbf{z}}_0$ that are closest to $z_k^{cp}$. We then use the locations of these $20$ nearest neighbors to identify the corresponding $20$ data points in $\mathbf{d}^{o}$, and take $ d_k^{o,cp}$ as the (equally weighted) mean of these $20$ data points. Of course, in general, $\hat{\mathbf{z}}_0$ and $\mathbf{d}^{o}$ may not be ``consistent''. This inconsistency, however, is partially taken into account by including $d_k^{o,cp}$ as a part of the inputs to the kernel function (cf. Eq. (\ref{eq:residual_rbf_approximation})), and assigning additional scale parameters ($\beta$) to adjust its influence in the course of data assimilation.  

In the experiments, we consider two scenarios. In the first one, we study the case in which there is no imperfection in the forward simulator $g(z)$, i.e., $g(z) = f(z) = \left(\vert z \vert^3 + 1 \right)^{1/2}$. Our objective here is to inspect the impact of kernel-based model-error correction (MEC) mechanism on the performance of data assimilation, when there is no imperfection in the forward simulator, but MEC is still adopted. For reference later, we call this perfect (simulator) scenario (PS). In the second scenario, we investigate the case in which imperfection indeed exists in the forward simulator, with $g(z) = z^2$. We examine how the performance of data assimilation may change in the presence of simulator imperfection. Likewise, we call this imperfect scenario (IS). 

\subsection*{Results in the perfect scenario (PS)} 
%
%

\begin{table*} 
	\centering
	\caption{\label{tab:rmse_dm_ps} Means and STDs of data mismatch and RMSE with respect to the initial ensemble, and the final ensembles with or without model-error correction (MEC), in the perfect scenario.}
	\begin{tabular}{||c||c||c||c||}
		\hline  
		& Initial ensemble & Final ensemble (no MEC) & Final ensemble (with MEC)  \\ 
		\hline 			
		Data mismatch (mean $\pm$ STD)   & $1.0694 \pm 0.5361 (\times 10^7)$ & $3.7326 \pm 0.0223 (\times 10^4)$ &  $6.2645 \pm 1.7551 (\times 10^4)$  \\
		\hline		
		RMSE (mean $\pm$ STD)		& $2.5240 \pm 0.3070$  & $1.0889 \pm 0.0025$ & $0.8836 \pm 0.0133$ \\
		\hline 
	\end{tabular}
\end{table*}   

\renewcommand{\nScale}{0.2}
\begin{figure} 
	\centering
	\begin{tabular}{cc}
		\subfloat[Reference model]{\includegraphics[width=0.45\textwidth]{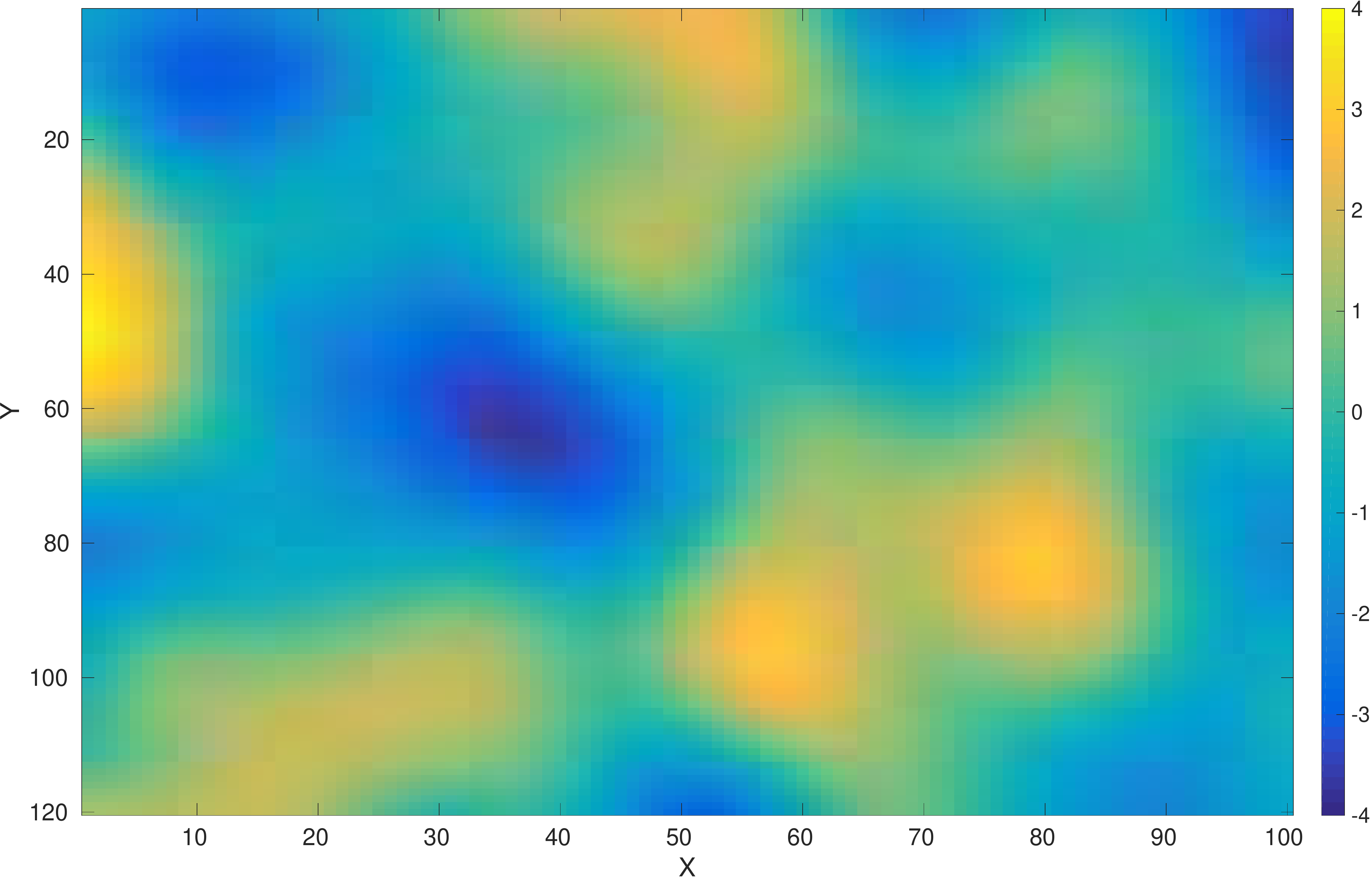}} &
		\subfloat[Mean of the initial ensemble]{\includegraphics[width=0.45\textwidth]{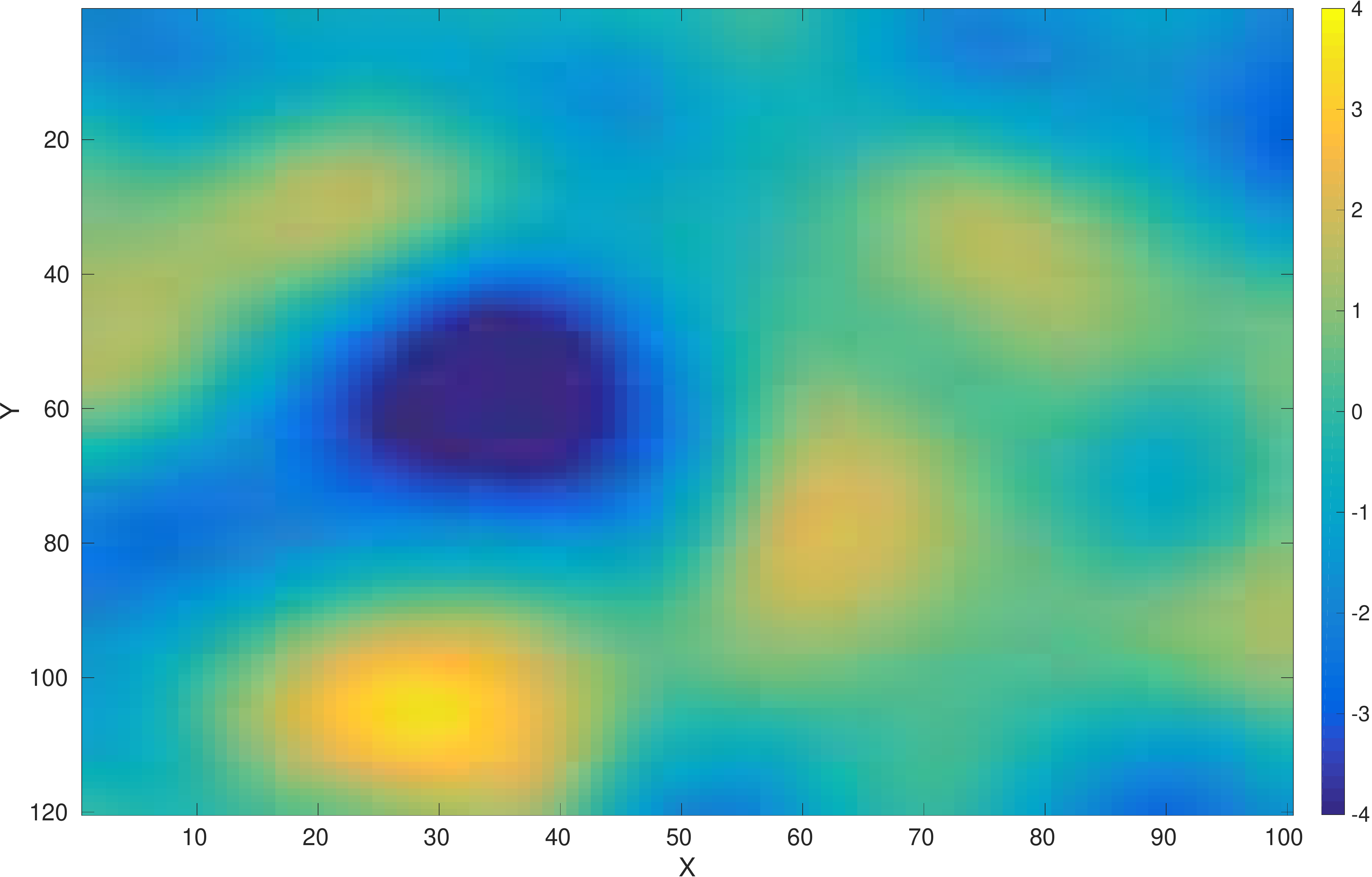}} \\
		\subfloat[Mean of the final ensemble \textit{without} MEC]{\includegraphics[width=0.45\textwidth]{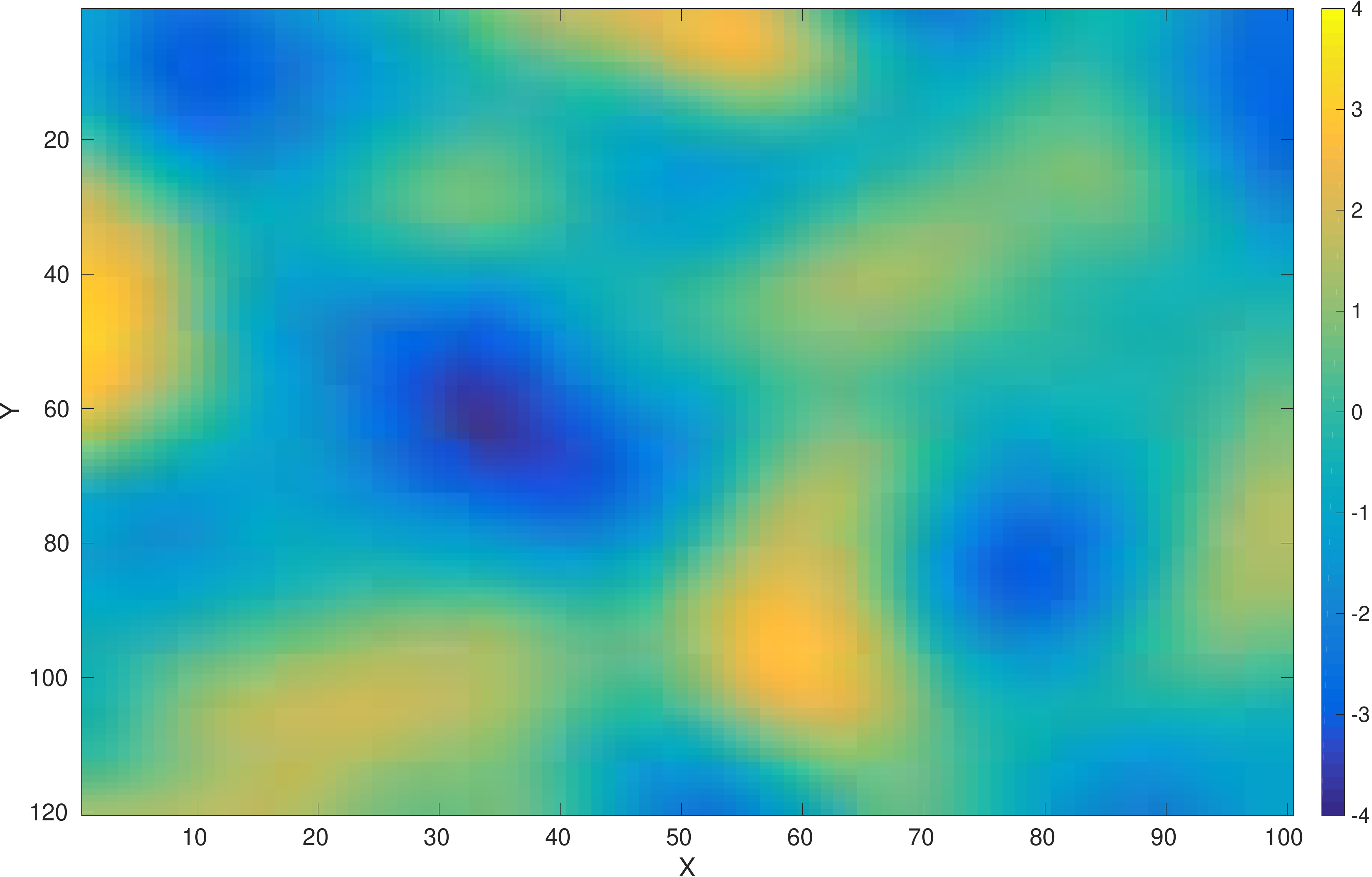}} &
		\subfloat[Mean of the final ensemble \textit{with} MEC]{\includegraphics[width=0.45\textwidth]{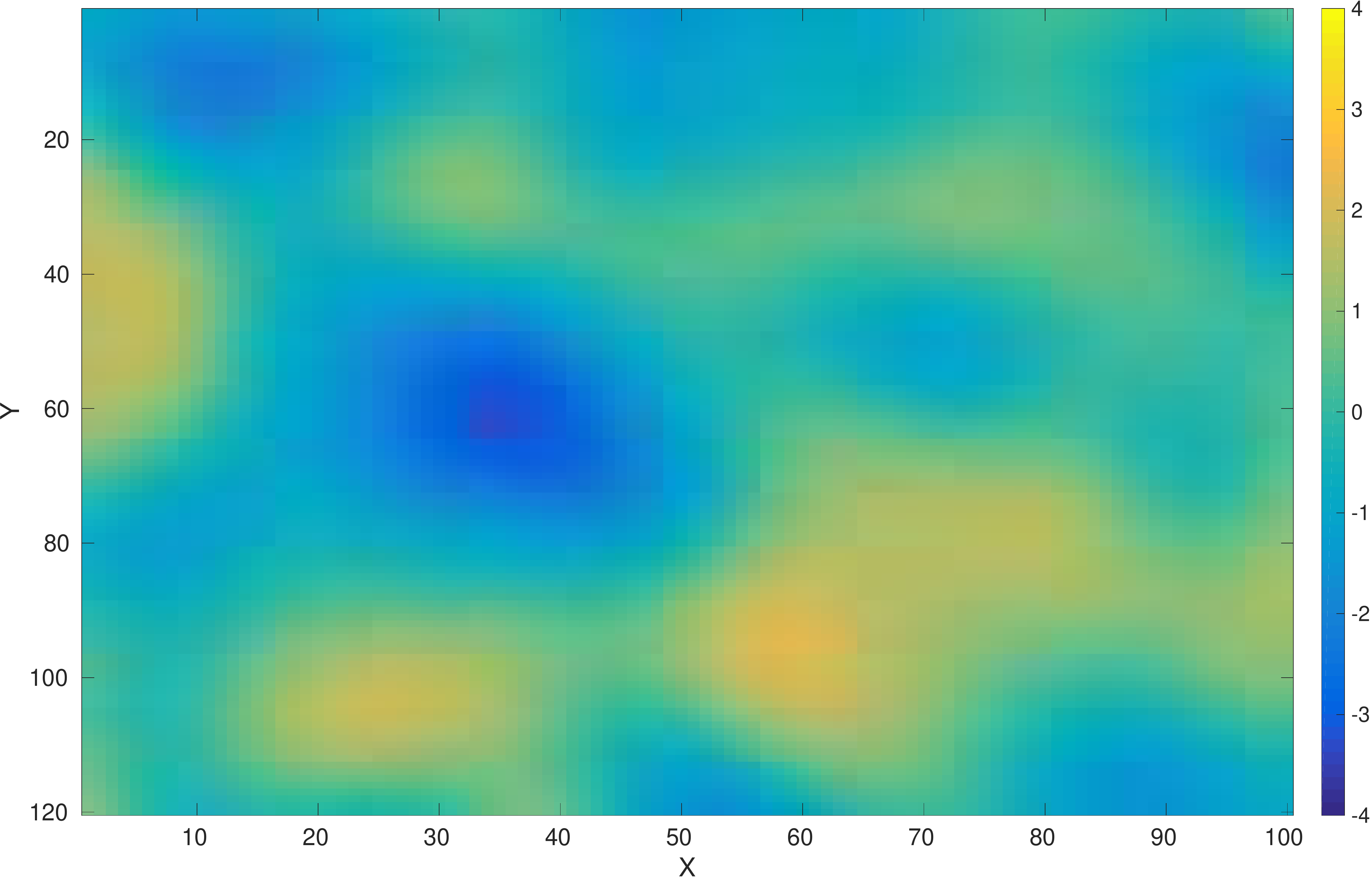}} \\
	\end{tabular}
	\caption{\label{fig:DA_models_perfect} Reference and mean models in the perfect scenario. Top row: Reference model (Panel (a)) used to generate observations (cf. Figure \ref{fig:DA_obs_perfect}(a)), and the mean model (Panel (b)) of the initial ensemble. Bottom row: mean of the final ensemble obtained through data assimilation without any model-error correction (MEC) (Panel (c)), and the corresponding mean when MEC is still adopted (Panel (d)) even though the forward simulator is perfect.}
\end{figure} 

\renewcommand{\nScale}{0.2}
\begin{figure} 
	\centering
	\begin{tabular}{cc}
		\subfloat[Real observations]{\includegraphics[width=0.45\textwidth]{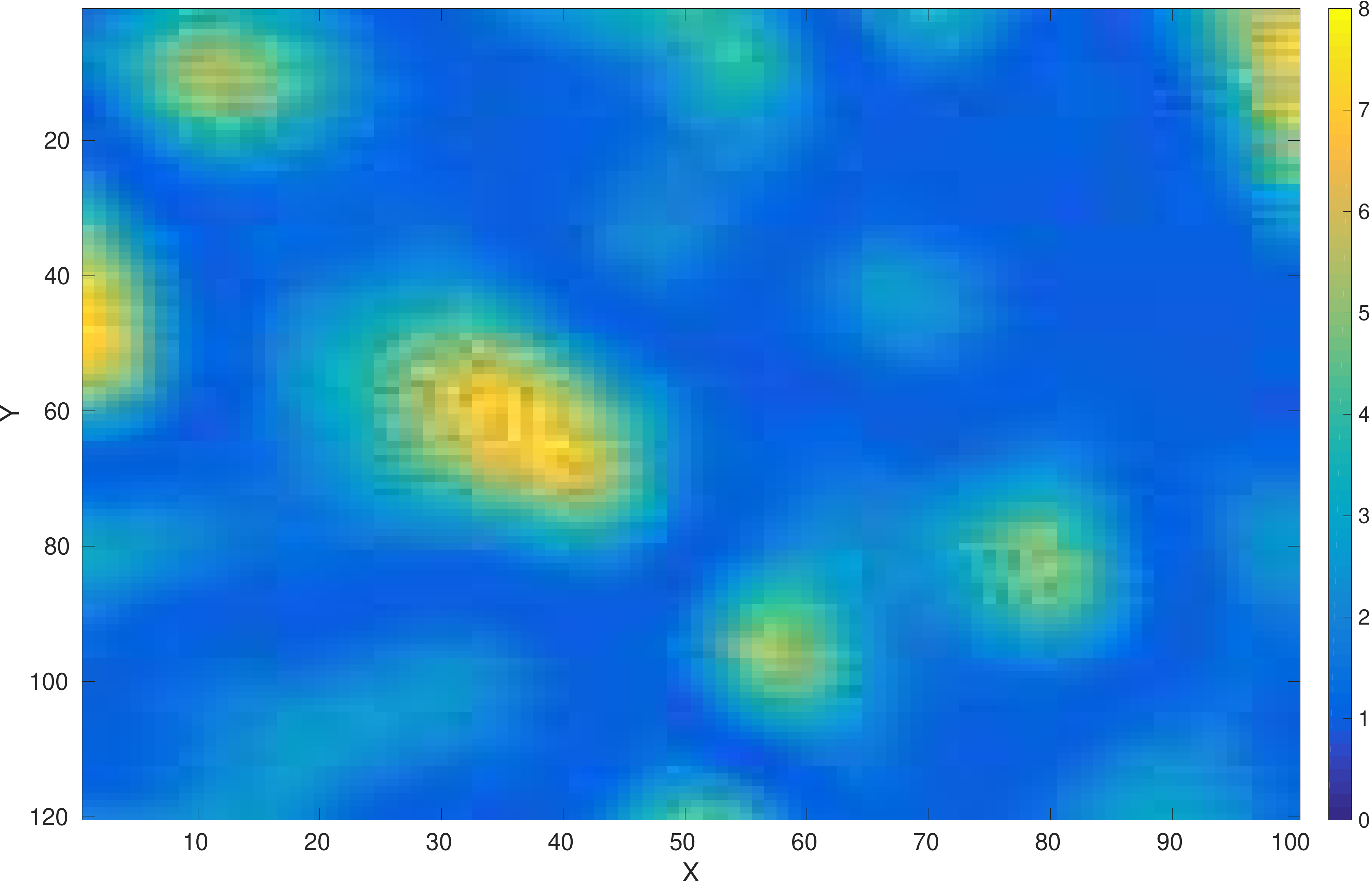}} &
		\subfloat[Mean of initial simulated observations]{\includegraphics[width=0.45\textwidth]{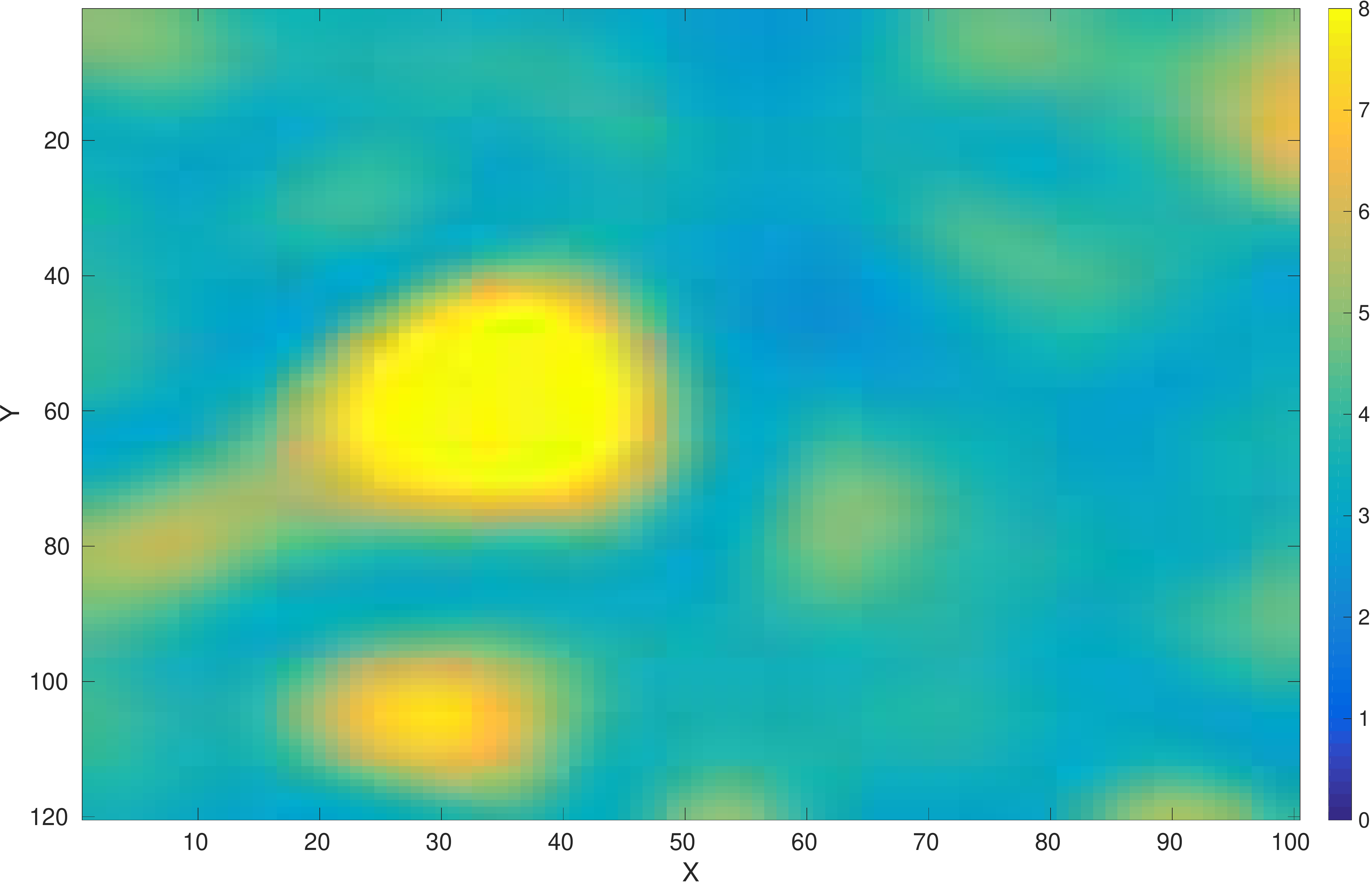}} \\
		\subfloat[Mean of final simulated observations \textit{without} MEC]{\includegraphics[width=0.45\textwidth]{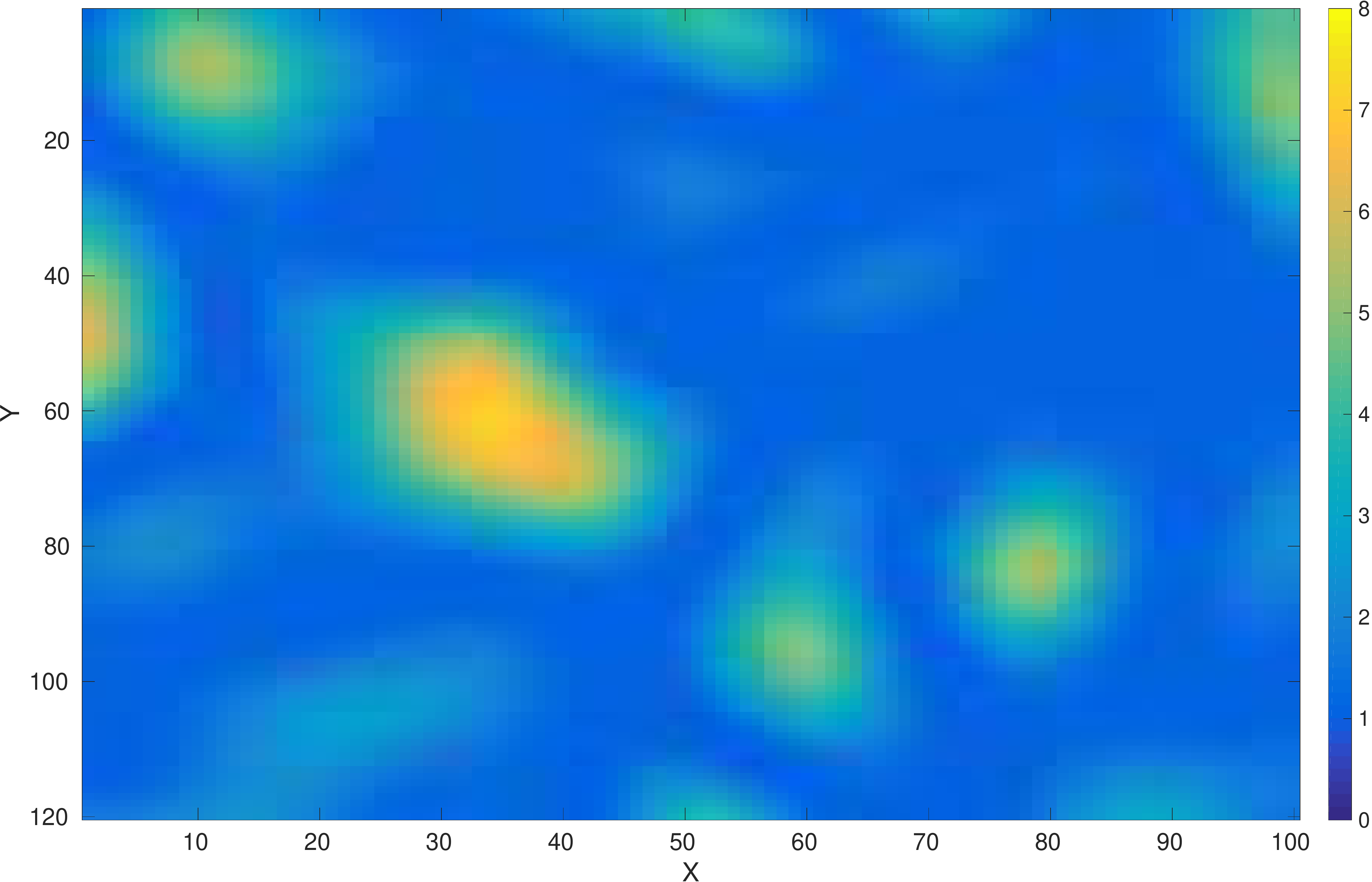}} &
		\subfloat[Mean of final simulated observations \textit{with} MEC]{\includegraphics[width=0.45\textwidth]{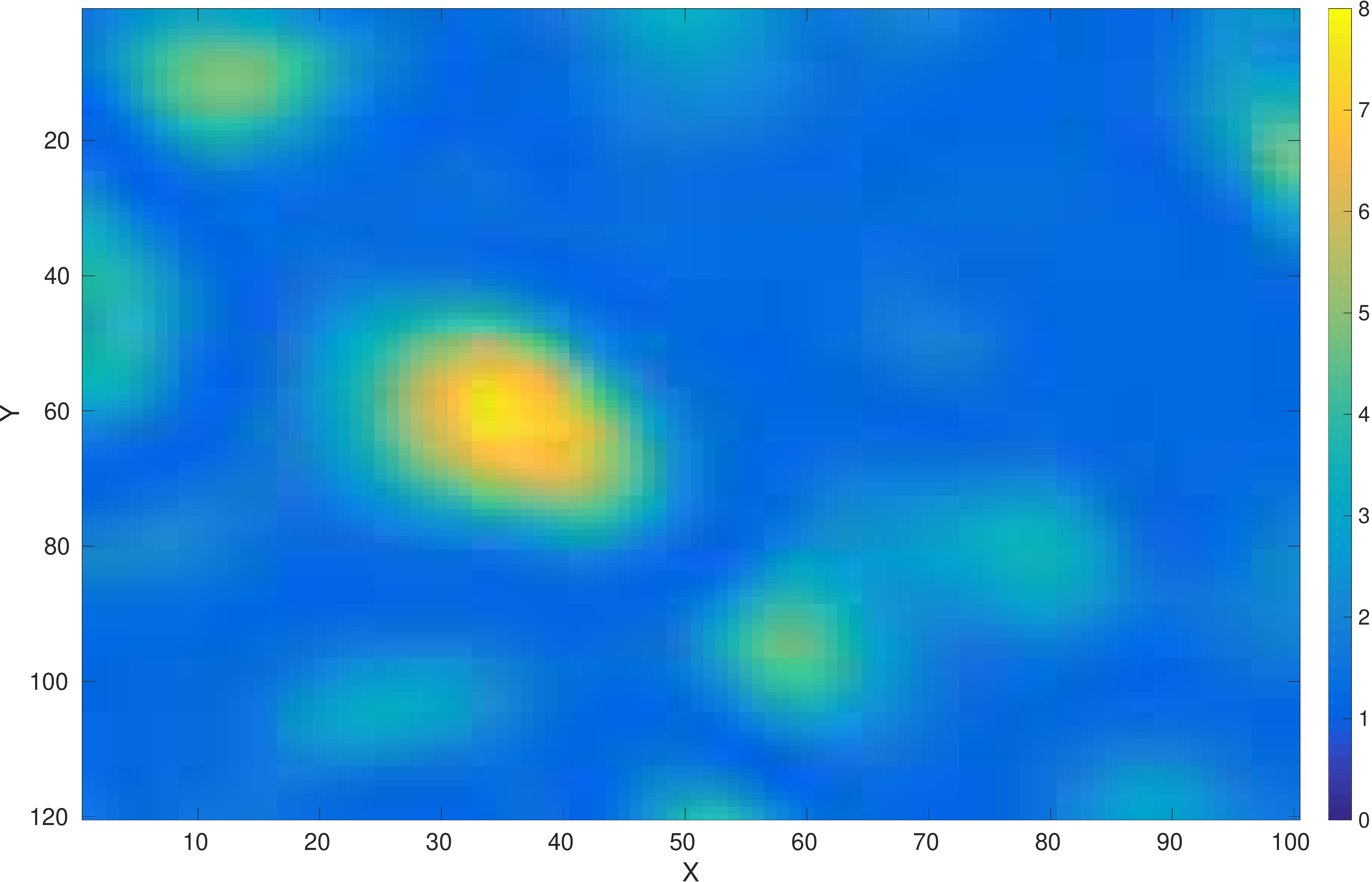}} \\
	\end{tabular}
	\caption{\label{fig:DA_obs_perfect} Real and simulated observations in the perfect scenario. Top row: Real observations (Panel (a)) generated using the reference model in Figure \ref{fig:DA_models_perfect}(a), and the mean of simulated observations obtained by applying the forward simulator to the initial ensemble of model variables (Panel (b)). Bottom row: As in Panel (b), but for the mean of simulated observations with respect to the final ensemble obtained without (Panel (c)) and with (Panel (d)) MEC in data assimilation, respectively.}
\end{figure} 

\renewcommand{\nScale}{0.2}
\begin{figure} 
	\centering
	\begin{tabular}{rr}
		\subfloat[Box plots of data mismatch \textit{without} MEC]{\includegraphics[width=0.45\textwidth]{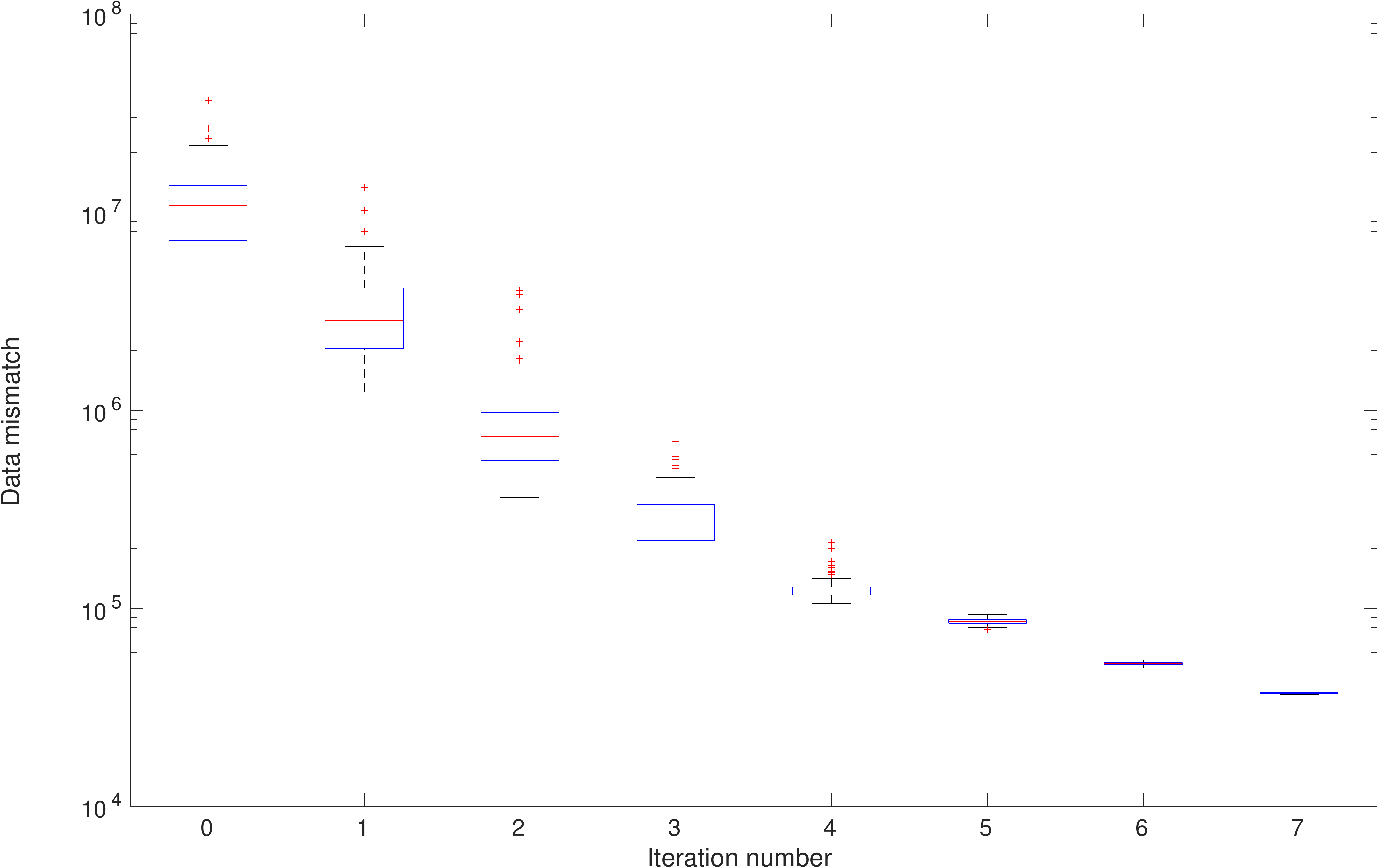}} &
		\subfloat[Box plots of data mismatch \textit{with} MEC]{\includegraphics[width=0.45\textwidth]{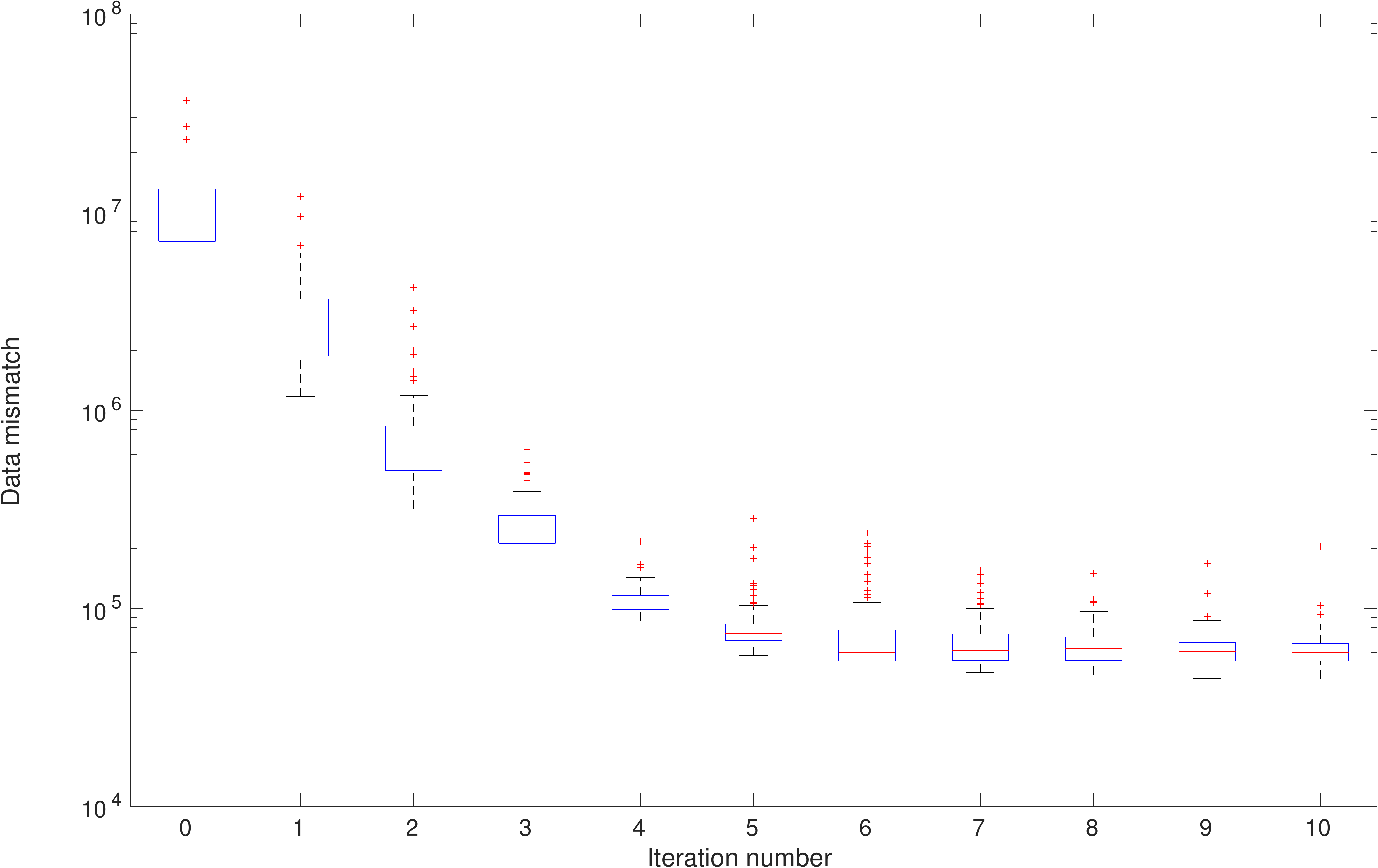}} \\[+0.5cm]
		\subfloat[Box plots of RMSE \textit{without} MEC]{\includegraphics[width=0.45\textwidth]{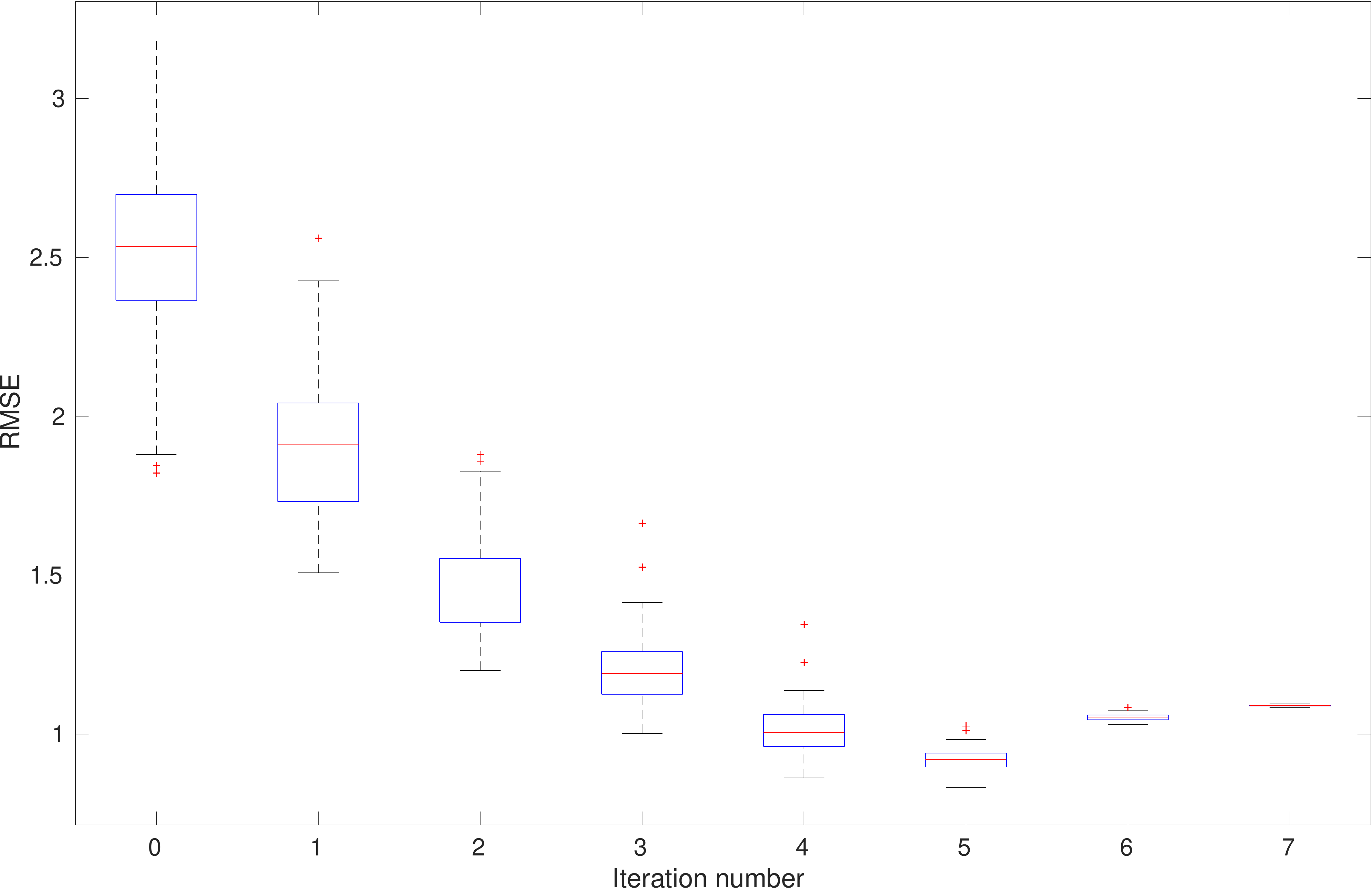}} &
		\subfloat[Box plots of RMSE \textit{with} MEC]{\includegraphics[width=0.45\textwidth]{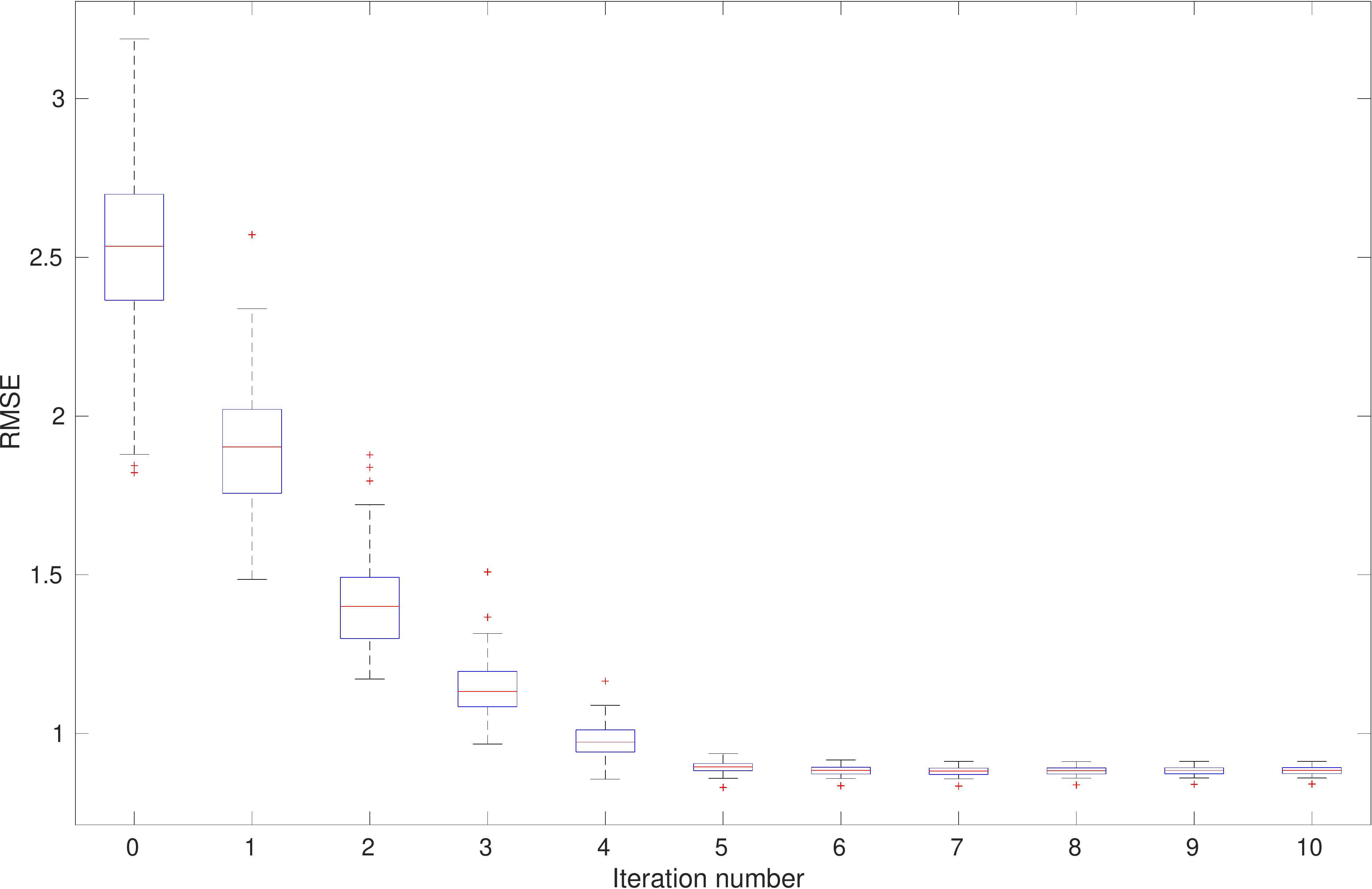}} \\
	\end{tabular}
	\caption{\label{fig:DA_boxplots} Box plots of data mismatch (top) and RMSE (bottom) with respect to the ensembles at different iteration steps in the perfect scenario. Results in Panels (a) and (c) correspond to the case without MEC adopted in data assimilation, whereas those in Panels (b) and (d) to the case with MEC. Unless otherwise stated, data mismatch in the experiment with MEC is always calculated using the modified forward simulator with a residual term, as in Eq. (\ref{eq:new_forward_simulator}). Note that in Panels (a) and (c), the iES terminates at the iteration step 7, due to the stopping criterion that the average data mismatch at this step is less than four times the number of observations (which is $4 \times 12,000$ in this case) for the first time.}
\end{figure}  

\renewcommand{\nScale}{0.3}
\begin{figure} 
	\centering
	\begin{tabular}{c}
		\subfloat{\includegraphics[width=0.8\textwidth]{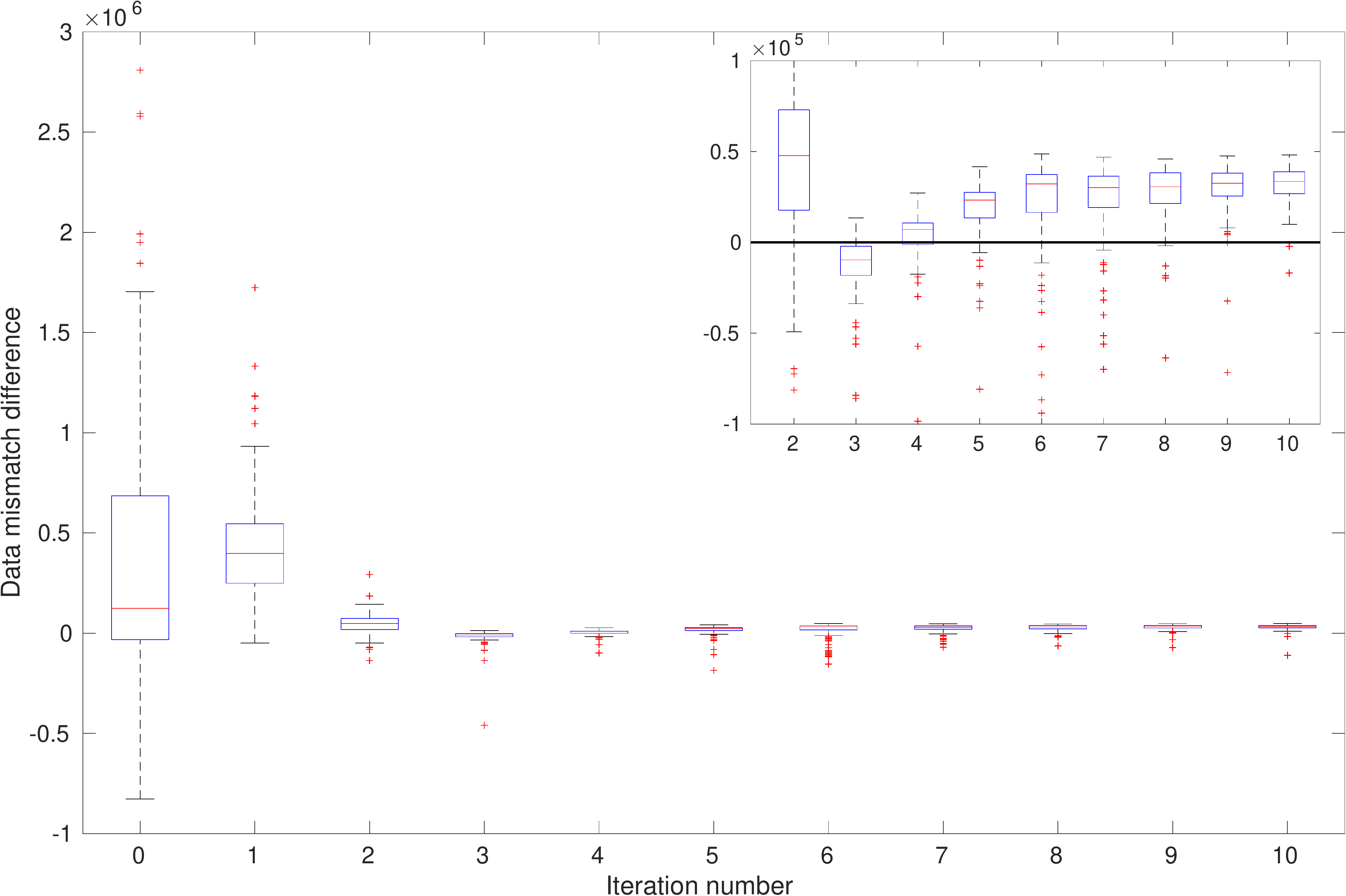}} \\
	\end{tabular}
	\caption{\label{fig:dmDiff_Boxplot_perfect} Box plots of data mismatch differences at different iteration steps, for the experiment with MEC in the perfect scenario. At a given iteration step, these differences are derived using data matching values that are calculated with the residual term excluded from Eq. (\ref{eq:new_forward_simulator}), minus data matching values that are computed with the residual term included in Eq. (\ref{eq:new_forward_simulator}), with respect to the corresponding ensemble of model variables at that iteration step. Therefore, positive data mismatch differences indicate that including the residual term helps match real observations better. For better visualization, we show the box plots from iteration steps 2 to 10 in a separate, zoomed-in subplot in the upper right corner.}
\end{figure}  

In the PS, we conduct a comparison study involving two experiments. In one of them, no MEC is adopted since the forward simulator is known to be perfect. In the other, kernel-based MEC is introduced to data assimilation, even though the forward simulator is perfect (in many places, we will simply say MEC when there is no confusion). Except for this difference, the other settings in these two experiments are identical. We note that, in the relevant experiment, MEC is conducted by combining Eqs. (\ref{eq:new_forward_simulator}) and (\ref{eq:predicted_residual}), whereas in Eq. (\ref{eq:predicted_residual}) the number $N_{cl}$ of cluster is set to $1$ in the current experiments. We will examine the impact of $N_{cl}$ on data assimilation in the IS. 

In comparison to the real observations in Figure \ref{fig:DA_obs_perfect}(a), the mean of simulated observations with respect to the initial ensemble (Figure \ref{fig:DA_obs_perfect}(b)) appears substantially different in many regions. As a result, the data mismatch values of the initial ensemble are relatively large, as reported in Table \ref{tab:rmse_dm_ps}. Through data assimilation using the iES, data mismatch values of the updated ensembles tend to decrease as the iteration process proceeds, whether MEC is introduced or not, as one can see in  Figures \ref{fig:DA_boxplots}(a) and \ref{fig:DA_boxplots}(b). Accordingly, after data assimilation, the means of simulated observations with respect to the final ensembles (with or without MEC), as shown in Figures \ref{fig:DA_obs_perfect}(c) and \ref{fig:DA_obs_perfect}(d), respectively, resemble the observations better than that with respect to the initial ensemble. 

In both experiments, the maximum iteration step of the iES is set to $10$. In the experiment without MEC introduced, however, the iES stops at the iteration step 7, due to an alternative stopping criterion that is triggered to terminate the iES, when the average data mismatch is lower than four times the number of observations (which is $4 \times 12,000$) for the first time. This early-stopping phenomenon indicates a higher risk of over-fitting observations, should the iteration process have continued after iteration step 7. On the other hand, in the experiment with MEC, since there are more parameters adopted in data assimilation, intuitively one might expect that the problem of over-fitting observations can be even more severe. Surprisingly, it turns out that over-fitting actually appears avoided, while the iteration stops at the  maximum step. As a result, the final mean data-mismatch value in the experiment with MEC is higher than that in the experiment without MEC, as reported in Table \ref{tab:rmse_dm_ps}. One possible explanation of the ability to avoid over-fitting may be related to the effect of localization \citep{luo2018correlation_spej}, which is rendered by the MEC mechanism, as will be discussed later. 

For quality check, in Figures \ref{fig:DA_boxplots}(c) and \ref{fig:DA_boxplots}(d) we also show the box plots of RMSEs of the ensembles of model variables at different iteration steps. When no MEC is introduced, the RMSEs tend to decrease at the first five iteration steps, and then bounce back to somewhat higher values at the last two iteration steps. This kind of ``U-turn'' behavior was also found in the earlier work of \cite{luo2016sparse2d_spej}, and can be mitigated or avoided by introducing a procedure of sparse data representation \citep{luo2016sparse2d_spej}, or localization \citep{luo2018correlation_spej,luo2018ECMOR_towards}, to the iES (an investigation on this issue, however, is beyond the scope of the current work). In contrast, with MEC introduced to the iES, the ``U-turn'' behavior seems vanished. Furthermore, the final mean RMSE value in the experiment with MEC is lower than that in the experiment without MEC, as one can see in Table \ref{tab:rmse_dm_ps}. In Figures \ref{fig:DA_models_perfect}(c) and \ref{fig:DA_models_perfect}(d), we show the mean of the final ensembles obtained in the experiments with or with MEC. Clearly, the mean models of the final ensembles appear more similar to the reference model than the mean model of the initial ensemble in Figure \ref{fig:DA_models_perfect}(b). The mean model of the final ensemble without MEC tends to do better than that with MEC in the regions on the left-hand side (e.g., for $\text{x} \leq 50$), but worse in the rest of the regions.

One can also observe an interesting phenomenon by comparing the spreads of box plots in Figure \ref{fig:DA_boxplots}, or the calculated ensemble STDs of data mismatch and RMSE in Table \ref{tab:rmse_dm_ps}. Recall that, in the experiments, no localization is introduced to the iES. As a result, it may not be surprising to see that, in the experiment without MEC, ensemble collapse seems to take place. In contrast, in the experiment with MEC, ensemble collapse does not appear to be a problem, or at least is mitigated. This seems to suggest that the kernel-based MEC mechanism can (partially) lead to the same effect on preventing ensemble collapse as localization does. A possible explanation to this phenomenon may be that, as aforementioned, since we use Gaussian RBF as the kernel function, the kernel parameters (scale and weight) associated with a certain center point would exhibit localized impacts, and mainly influence model variables that are sufficiently close to that  center point. 

As aforementioned, in SLP, typically one has many (matched) input-output pairs as the training data. In contrast, in data assimilation problems, we use a single realization (or one-shot) of the observations (at a given time instance and a given spatial location) to infer possible model variables. As a result, in SLP, one often has the luxury to split a dataset into two parts, one for training (and cross-validation) and one for test; whereas in data assimilation with MEC, this kind of luxury typically does not exist. This makes MEC a particularly challenging problem. Indeed, apart from the potential inconsistencies between the observations and the estimated model variables, there are only one-shot observations used for residual functional estimation, which makes it difficult for the updated forward simulator to generalize to other unseen training data (e.g., new input-output pairs), as our experiments indicate (results not shown). 

Bearing the above challenges in mind, when evaluating the performance of MEC, we do not particularly focus on inspecting the generalization ability of the updated forward simulator (after all, the goal of data assimilation is to estimate the ground truth corresponding to real observations). Instead, we adopt the following cross-validation procedure, namely, for a given ensemble of model variables in the experiment with MEC, we compare the corresponding data mismatch values, when the residual term $\hat{\mathbf{r}}\left(\mathbf{z}, \boldsymbol{\eta}\right)$ is used or not used in the forward simulator (cf. Eq. (\ref{eq:new_forward_simulator})). Such a comparison aims to examine whether the introduction of the residual term to the forward simulator helps match real observations better or not.   

Following this notion, Figure \ref{fig:dmDiff_Boxplot_perfect} shows the box plots of data mismatch differences at different iteration steps, with respect to the experiment with MEC. At a given iteration step (hence a given ensemble of model variables), these differences are obtained by subtracting data mismatch values which are calculated with the residual term in the modified forward simulator in Eq. (\ref{eq:new_forward_simulator}) (as in Figure \ref{fig:DA_boxplots}(b)), from the corresponding data mismatching values which are calculated without including the residual term in Eq. (\ref{eq:new_forward_simulator}). Positive difference values in the box plots thus imply that the presence of the residual term in Eq. (\ref{eq:new_forward_simulator}) is useful for helping match real observations better, and vice versa. From this perspective, Figure \ref{fig:dmDiff_Boxplot_perfect} suggests that, with the initial ensemble of kernel parameters, the effect of including the residual term in Eq. (\ref{eq:new_forward_simulator}) at iteration step $0$ is mixed, and there are substantial numbers of difference values residing on both sides of zero (although overall the number of positive values does seem to dominate). After one iteration (at iteration step $1$), the model qualities are improved in terms of RMSE (cf. Figure \ref{fig:DA_boxplots}(d)), meanwhile the number of positive difference values also increases. However, as models are further improved, the number of positive difference values does not necessarily always dominate, as one can spot in the box plot at iteration step $3$. Nevertheless, as the iteration process continues, this kind of ``over-correction'' diminishes. Eventually, the number of positive difference values dominates at the final stage, while the RMSEs of estimated models tend to gradually reduce.                    

Based on the experiment results in the PS, we conclude that, in this particular case study, even though the forward simulator is perfect, it appears still beneficial to integrate kernel-based MEC into data assimilation for performance improvements.

\subsection*{Results in the imperfect scenario (IS)}  
%
%

\begin{table*} 
	\centering
	\caption{\label{tab:rmse_dm_IS} Means and STDs of data mismatch and RMSE with respect to the initial ensemble, and the final ensembles with or without model-error correction (MEC), in the \textit{imperfect} scenario.}
	\begin{tabular}{||c||c||c||c||}
		\hline  
		& Initial ensemble & Final ensemble (no MEC) & Final ensemble (with MEC)  \\ 
		\hline 			
		Data mismatch (mean $\pm$ STD)   & $6.5372 \pm 5.2423 (\times 10^7)$ & $4.5211 \pm 0.0590 (\times 10^5)$ &  $1.3248 \pm 1.2528 (\times 10^5)$  \\
		\hline		
		RMSE (mean $\pm$ STD)		& $2.5240 \pm 0.3070$  & $1.2053 \pm 0.0091$ & $1.0696 \pm 0.0174$ \\
		\hline 
	\end{tabular}
\end{table*}   

\renewcommand{\nScale}{0.2}
\begin{figure} 
	\centering
	\begin{tabular}{cc}
		\subfloat[Reference model]{\includegraphics[width=0.45\textwidth]{DA_truth.pdf}} &
		\subfloat[Mean of the initial ensemble]{\includegraphics[width=0.45\textwidth]{DA_initMean.pdf}} \\
		\subfloat[Mean of the final ensemble \textit{without} MEC]{\includegraphics[width=0.45\textwidth]{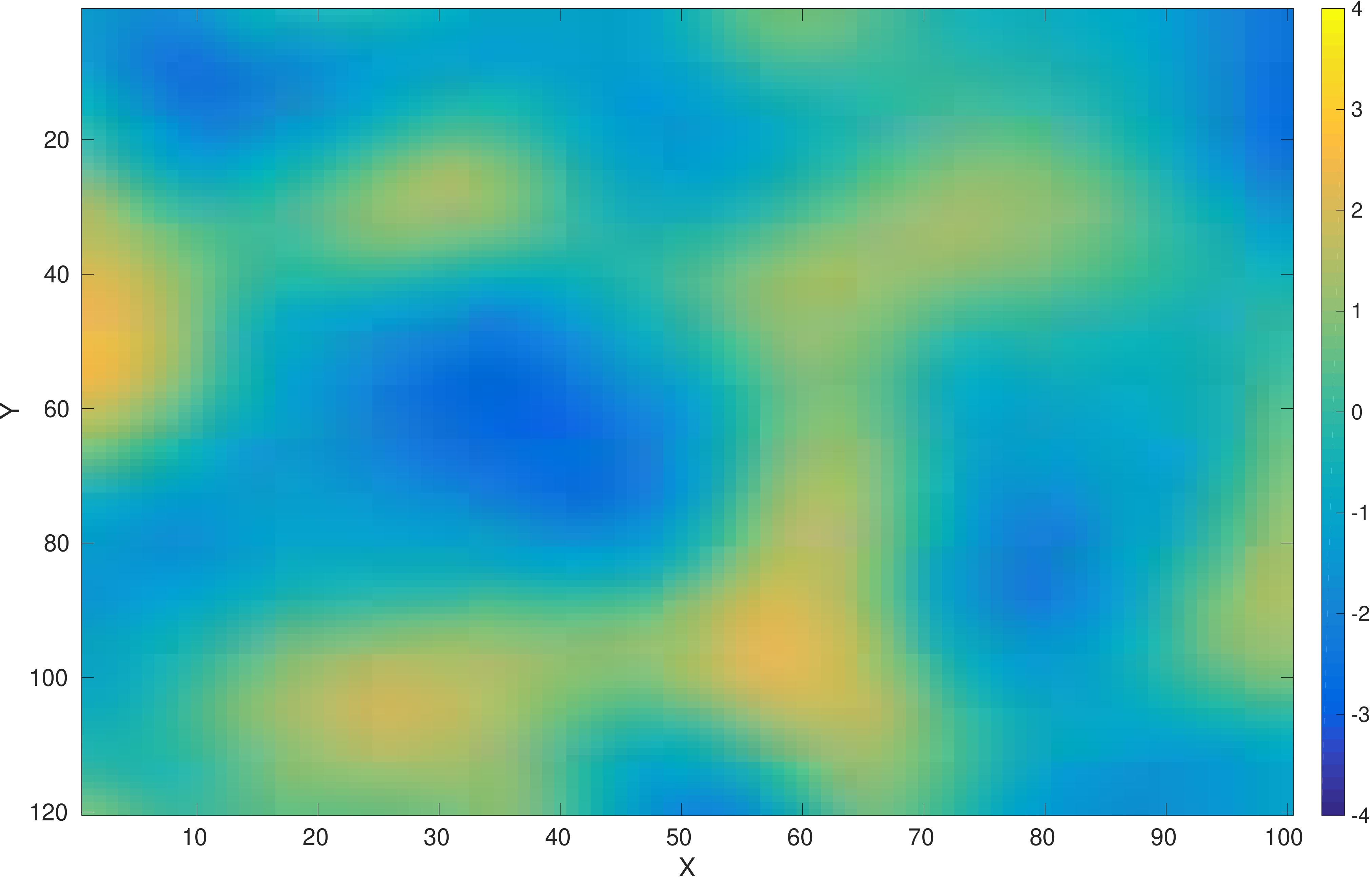}} &
		\subfloat[Mean of the final ensemble \textit{with} MEC]{\includegraphics[width=0.45\textwidth]{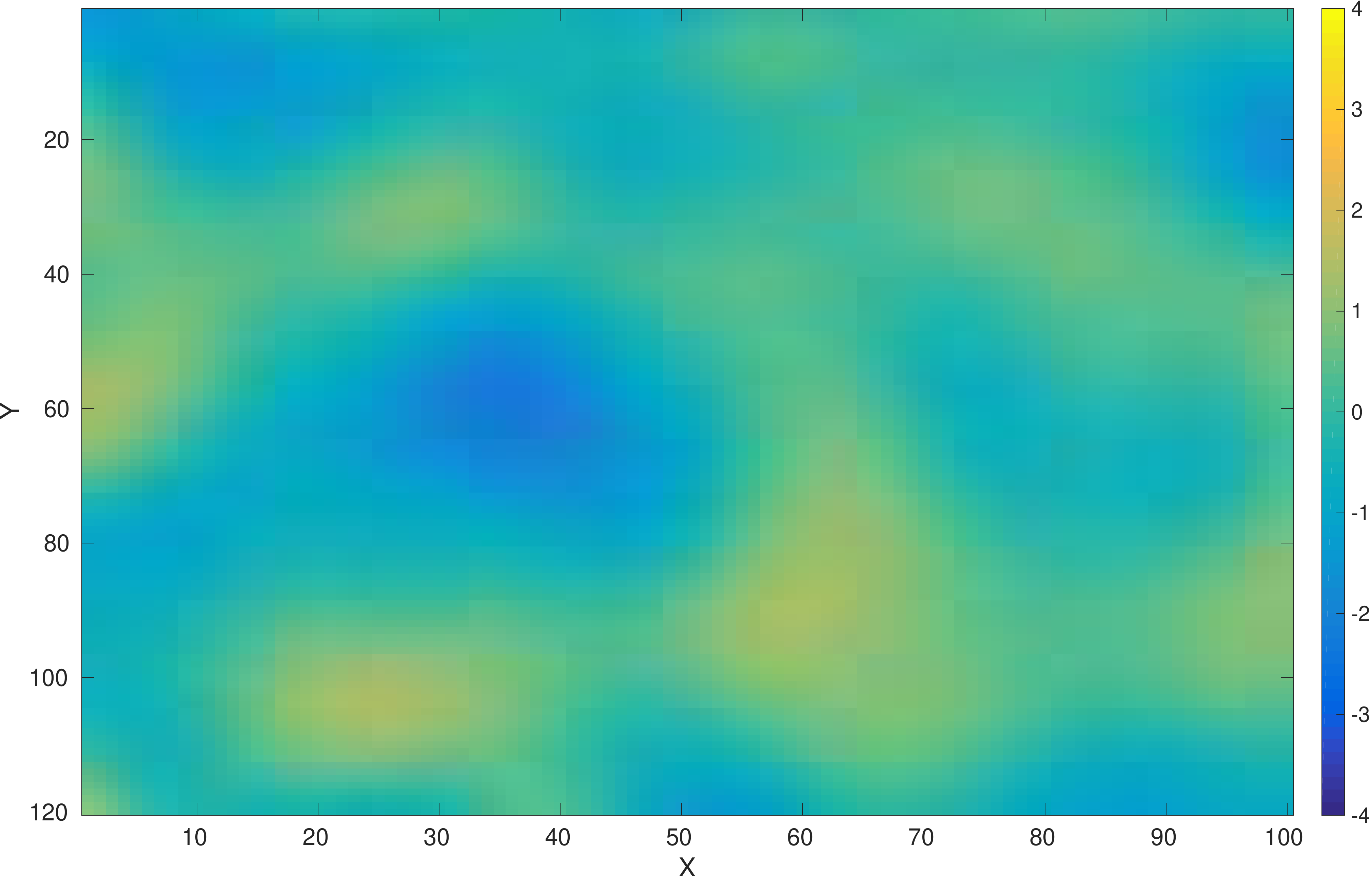}} \\
	\end{tabular}
	\caption{\label{fig:DA_models_imperfect} As in Figure \ref{fig:DA_models_perfect}, but for the experiment results in the \textit{imperfect} scenario ($N_{cl} = 1$).}
\end{figure} 

\renewcommand{\nScale}{0.2}
\begin{figure} 
	\centering
	\begin{tabular}{cc}
		\subfloat[Real observations]{\includegraphics[width=0.45\textwidth]{DA_trueObs.pdf}} &
		\subfloat[Mean of initial simulated observations]{\includegraphics[width=0.45\textwidth]{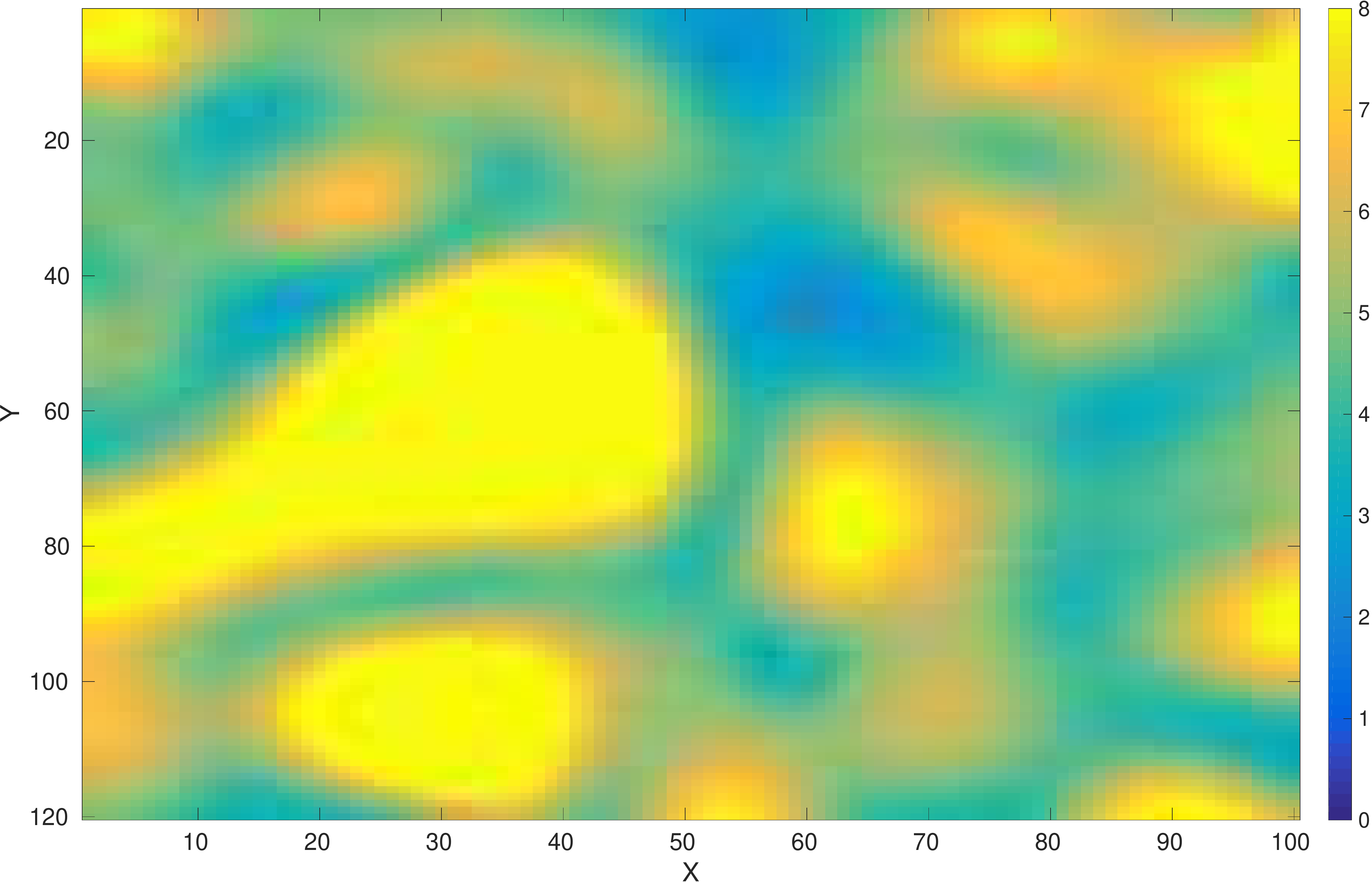}} \\
		\subfloat[Mean of final simulated observations \textit{without} MEC]{\includegraphics[width=0.45\textwidth]{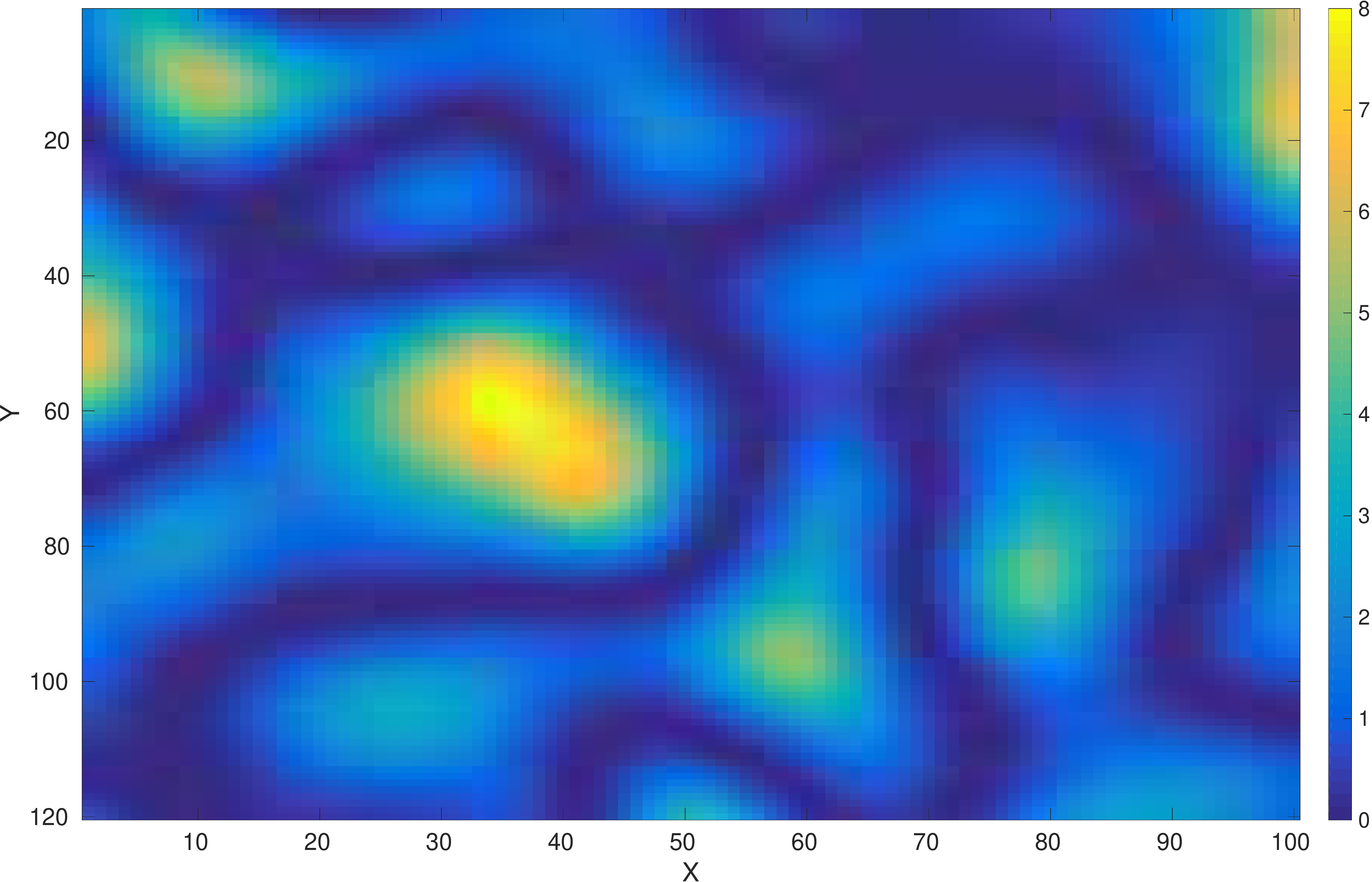}} &
		\subfloat[Mean of final simulated observations \textit{with} MEC]{\includegraphics[width=0.45\textwidth]{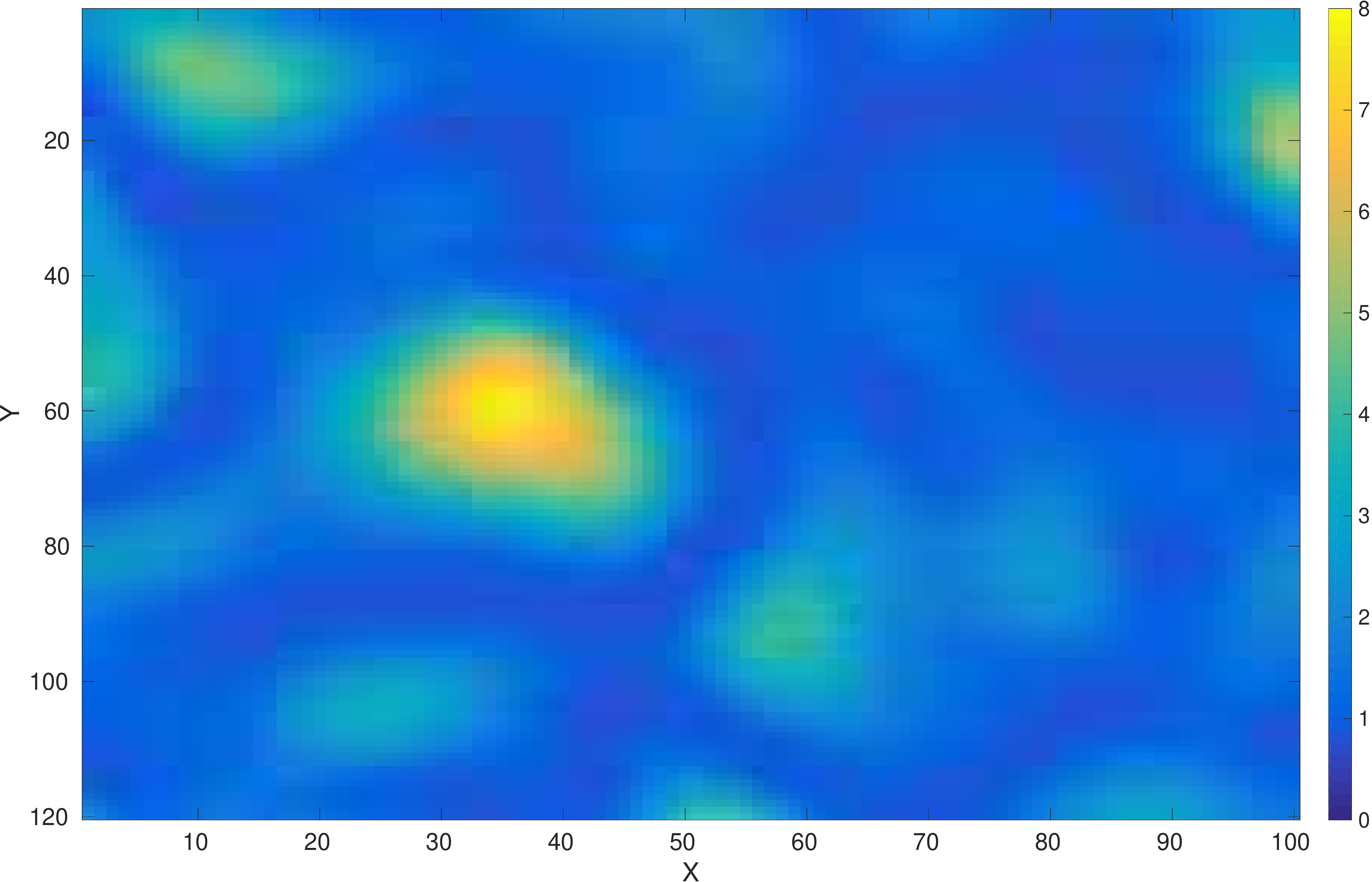}} \\
	\end{tabular}
	\caption{\label{fig:DA_obs_imperfect} As in Figure \ref{fig:DA_obs_perfect}, but for the experiment results in the \textit{imperfect} scenario ($N_{cl} = 1$).}
\end{figure} 

\renewcommand{\nScale}{0.2}
\begin{figure} 
	\centering
	\begin{tabular}{rr}
		\subfloat[Box plots of data mismatch \textit{without} MEC]{\includegraphics[width=0.45\textwidth]{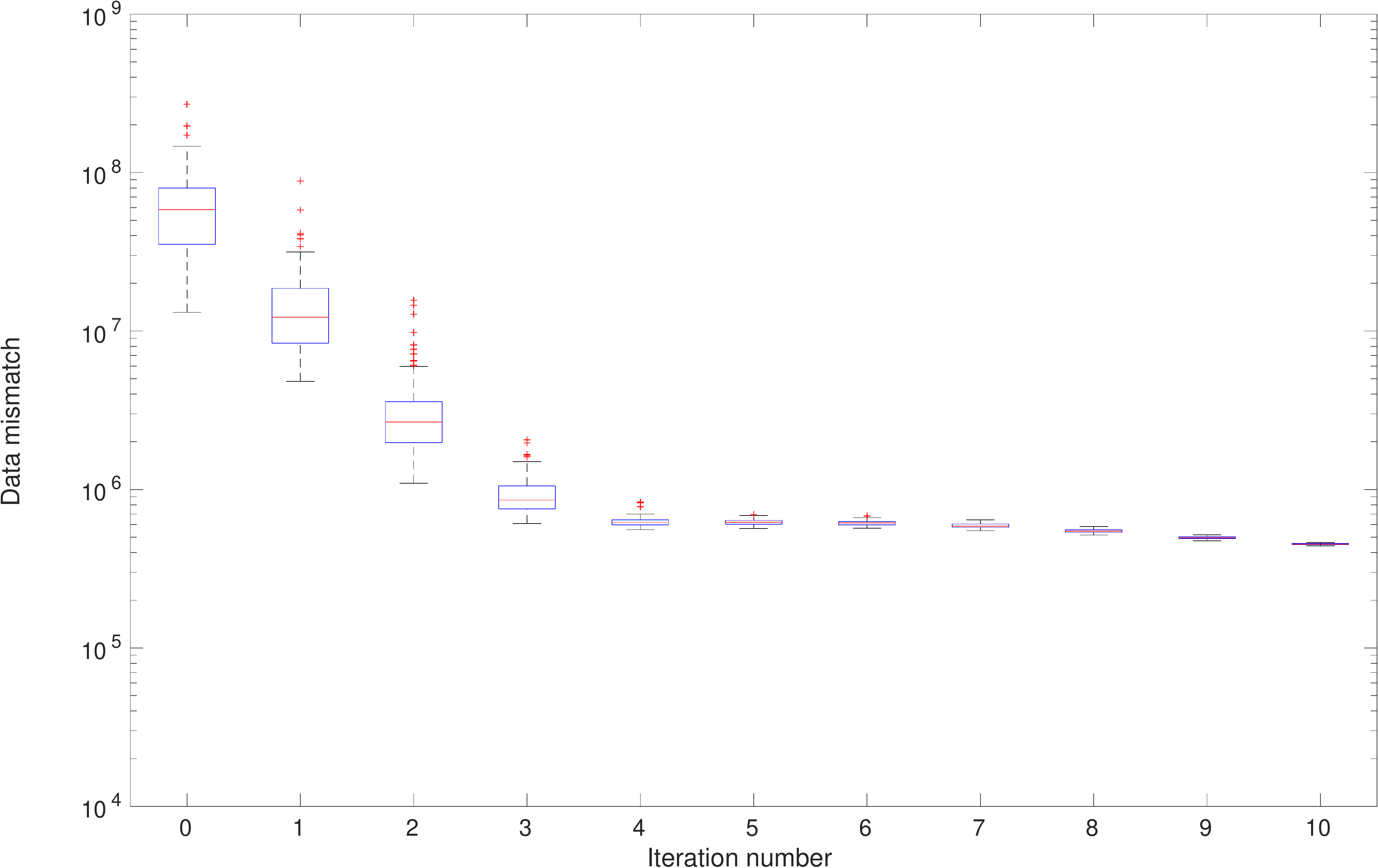}} &
		\subfloat[Box plots of data mismatch \textit{with} MEC]{\includegraphics[width=0.45\textwidth]{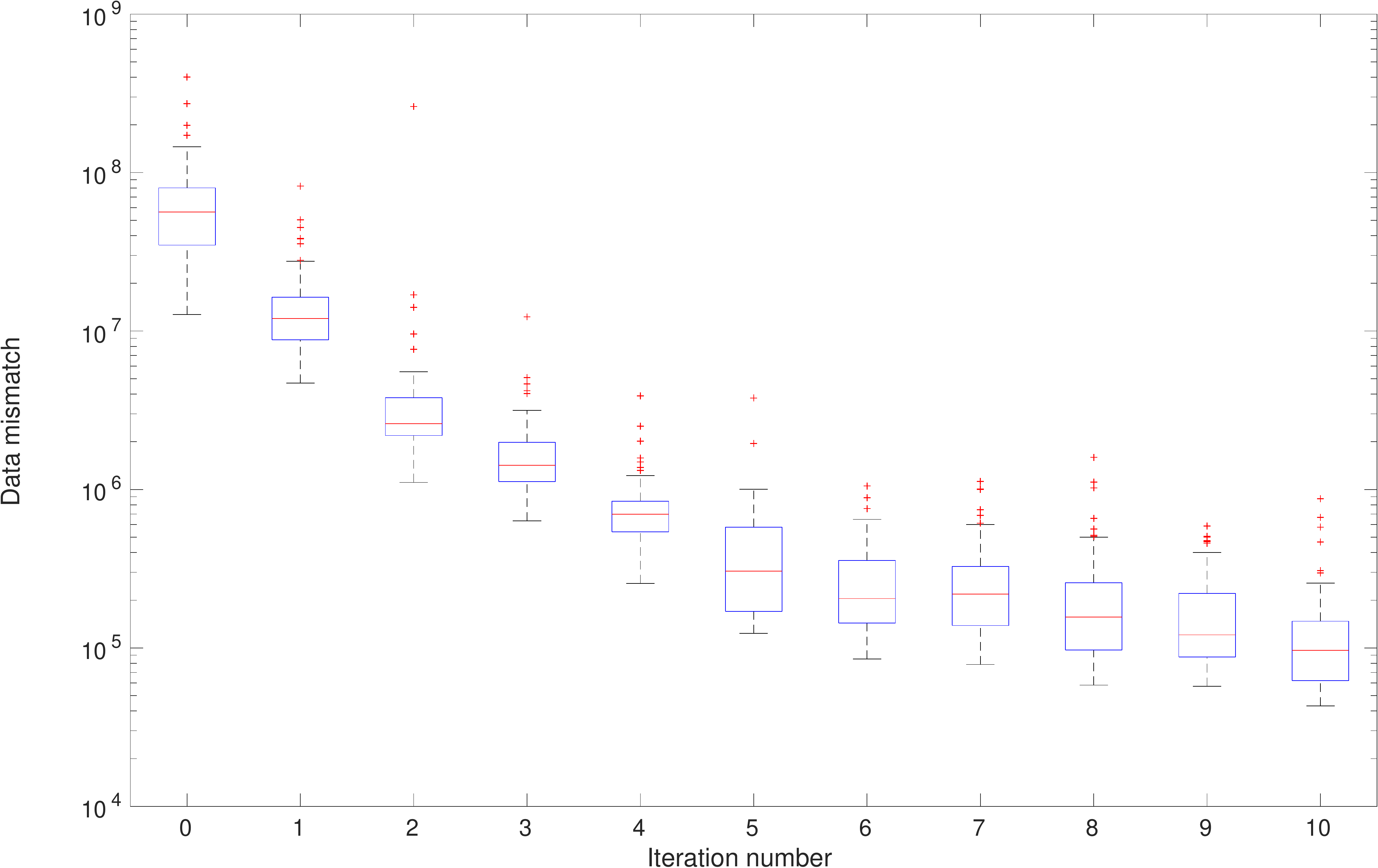}} \\[+0.5cm]
		\subfloat[Box plots of RMSE \textit{without} MEC]{\includegraphics[width=0.45\textwidth]{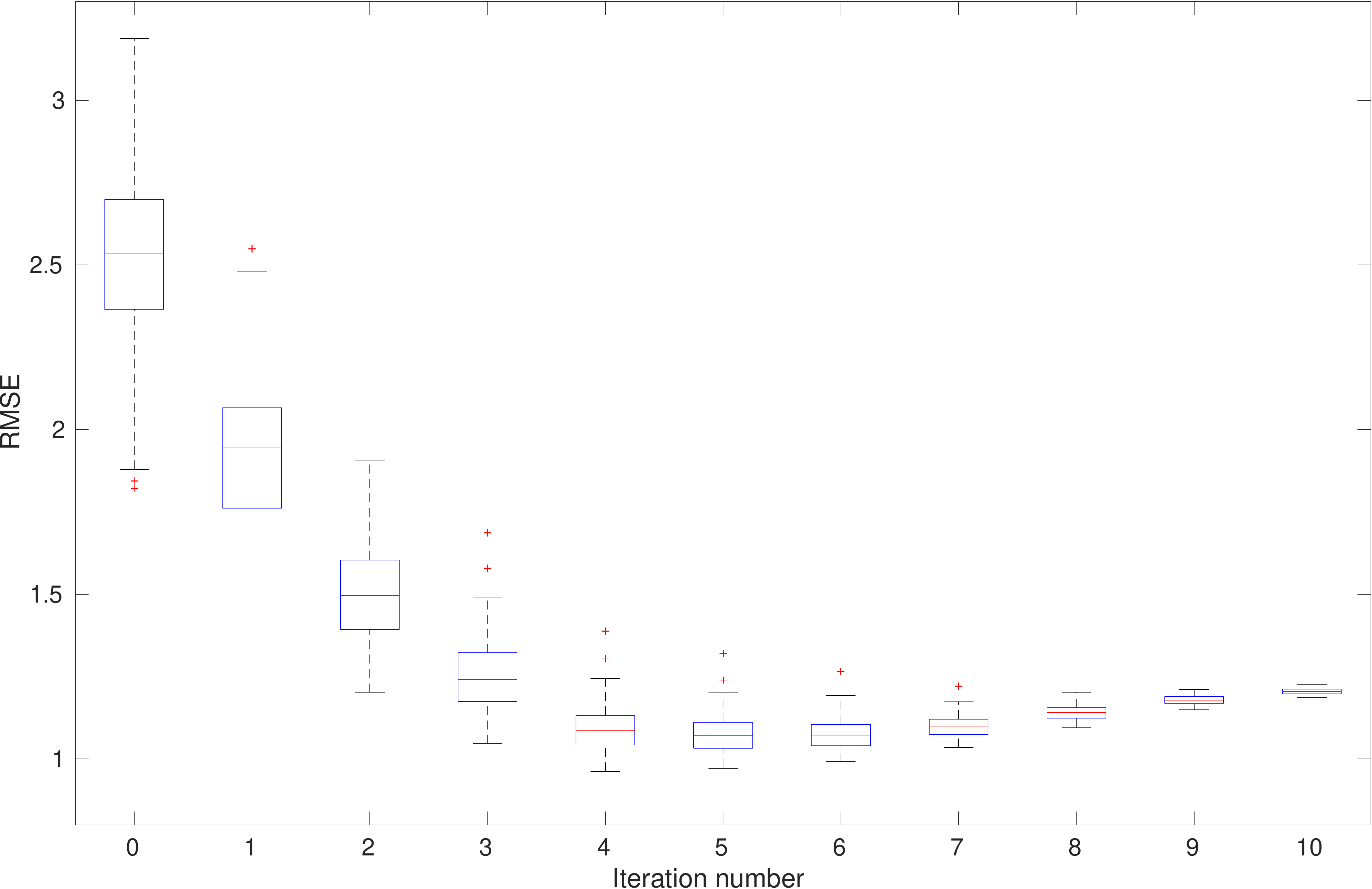}} &
		\subfloat[Box plots of RMSE \textit{with} MEC]{\includegraphics[width=0.45\textwidth]{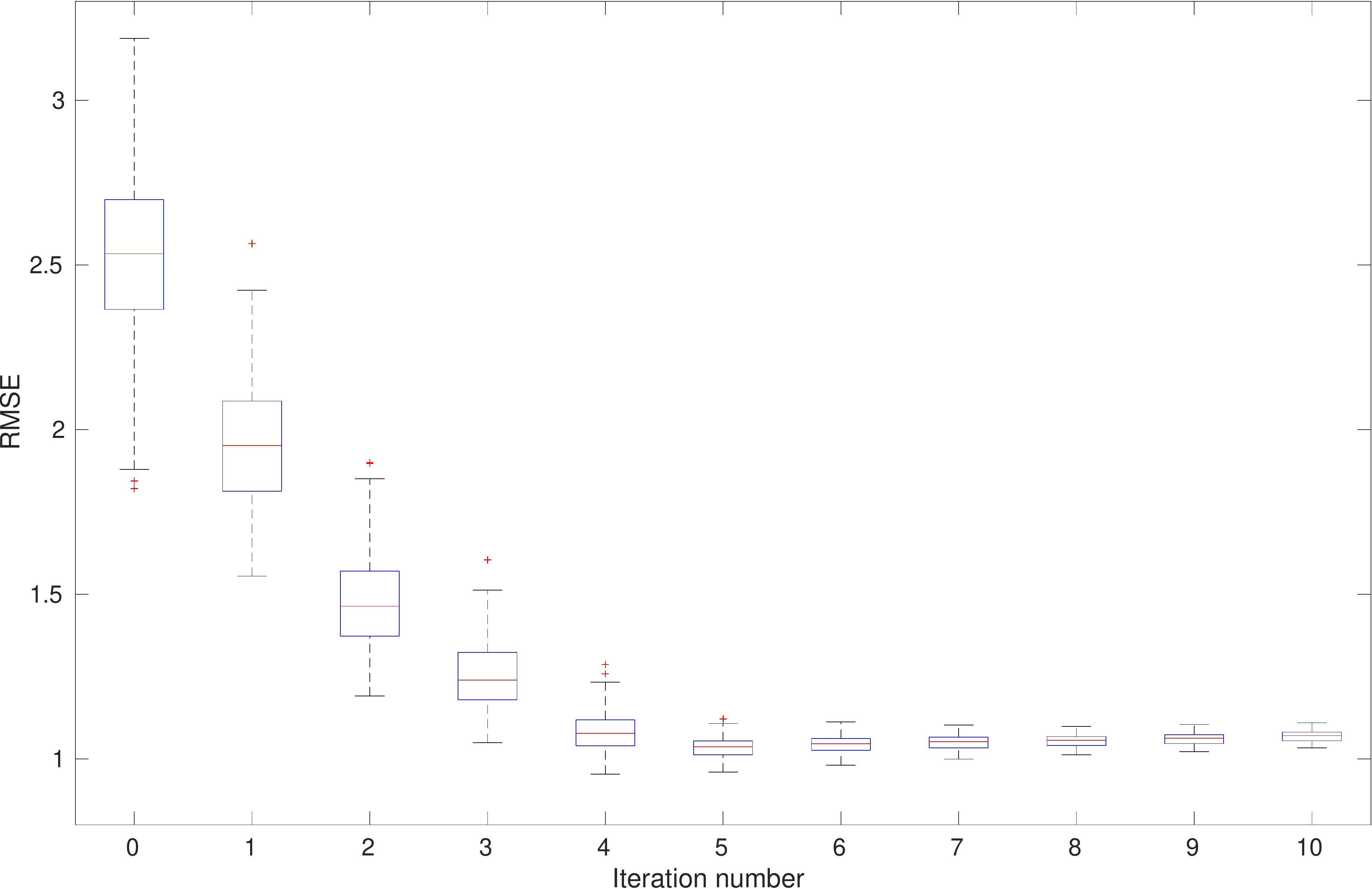}} \\
	\end{tabular}
	\caption{\label{fig:DA_boxplots_imperfect} As in Figure \ref{fig:DA_boxplots}, but for the experiment results in the \textit{imperfect} scenario ($N_{cl} = 1$).}
\end{figure}  

\renewcommand{\nScale}{0.4}
\begin{figure} 
	\centering
	\begin{tabular}{c}
		\subfloat{\includegraphics[width=0.8\textwidth]{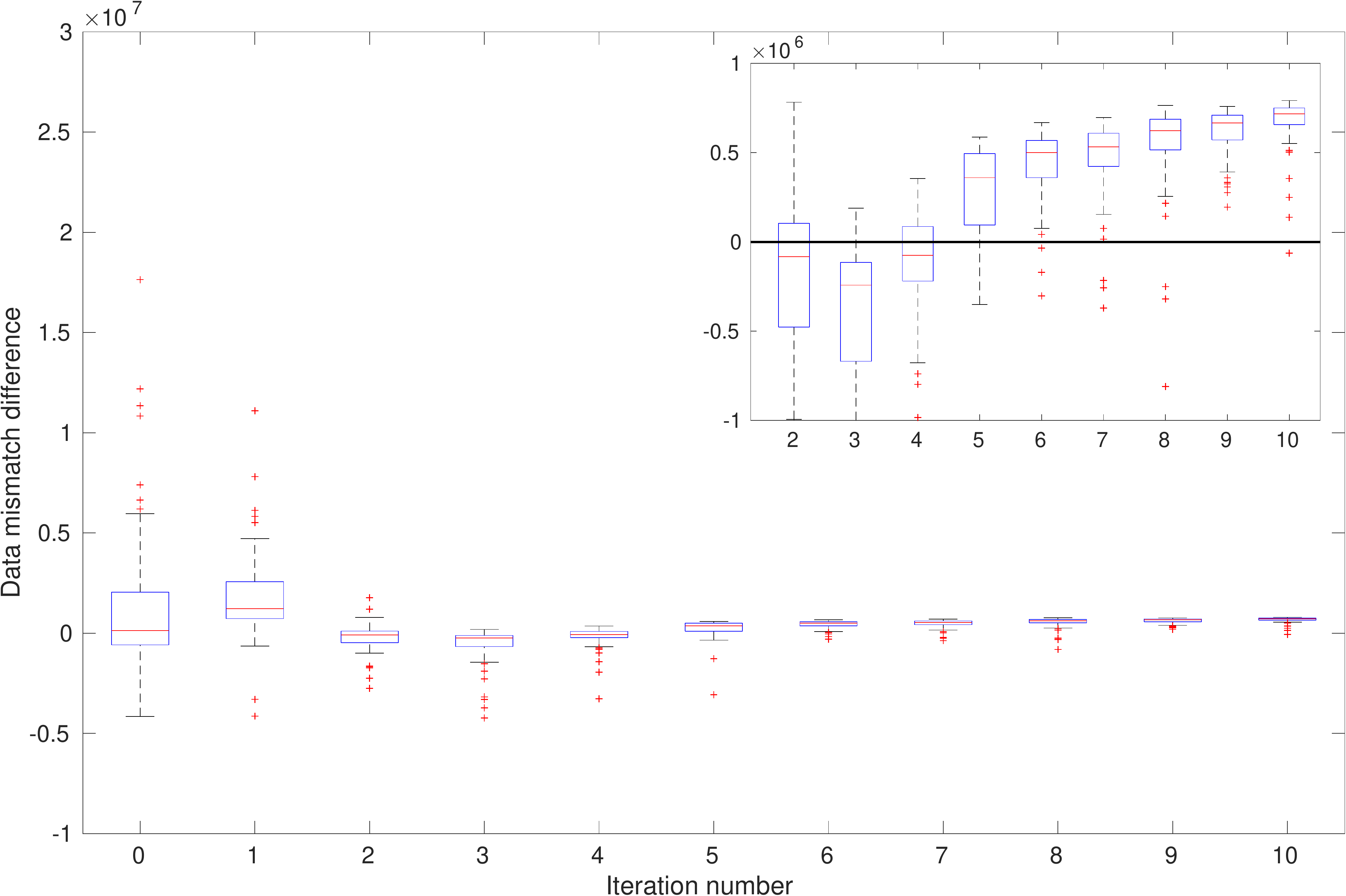}} \\
	\end{tabular}
	\caption{\label{fig:dmDiff_Boxplot_imperfect} As in Figure \ref{fig:dmDiff_Boxplot_perfect}, but for the experiment results in the \textit{imperfect} scenario ($N_{cl} = 1$).}
\end{figure}  

\subsubsection*{Results with $N_{cl} = 1$}       
In parallel to the results in the PS, we first report the results with $N_{cl} = 1$ in the IS. In this case, we also compare the assimilation performance with respect to one experiment where there is no MEC introduced, and another experiment where kernel-based MEC is adopted, with the number of cluster $N_{cl} = 1$. The initial ensemble of kernel parameters is generated in the way as in the case study of SLP.

Table \ref{tab:rmse_dm_IS} reports both data mismatch and RMSE (in terms of mean $\pm$ STD) for the initial ensemble, and the final ensembles obtained when MEC is or is not adopted. For the purpose of comparison, we adopt the same initial ensemble as in the PS. From Table \ref{tab:rmse_dm_IS}, one can again see that the use of MEC helps reduce mean values of both data mismatch and RMSE, while retaining higher ensemble spreads in the final ensemble, in comparison to the choice in which MEC is not used. In addition, by comparing Tables \ref{tab:rmse_dm_ps} and \ref{tab:rmse_dm_IS}, one also spots the impact of imperfection on data assimilation: In the presence of imperfection, the performance of data assimilation is worsened, with mean values of both data mismatch and RMSE in the IS becoming larger than those in the PS.

The subsequent results in Figures \ref{fig:DA_models_imperfect} -- \ref{fig:dmDiff_Boxplot_imperfect} are shown in analogy to their counterparts, Figures \ref{fig:DA_models_perfect} -- \ref{fig:dmDiff_Boxplot_perfect}, in the PS, respectively. A comparison between these figures are largely consistent with our observations stated in the preceding paragraph. In particular, the box plots of data mismatch differences at different iteration steps, as shown in Figure \ref{fig:dmDiff_Boxplot_imperfect}, also indicate that kernel-based MEC is useful for improving data match to real observations. 

Overall, the experiment results presented here confirm again that, in this particular case study, kernel-based MEC helps improve the performance of data assimilation in the presence of imperfection in the forward simulator.   

\subsubsection*{Comparison to an alternative MEC mechanism} 
\renewcommand{\nScale}{0.2}
\begin{figure} 
	\centering
	\begin{tabular}{rr}
		\subfloat[Mean of final simulated observations with bias-based MEC]{\includegraphics[width=0.45\textwidth]{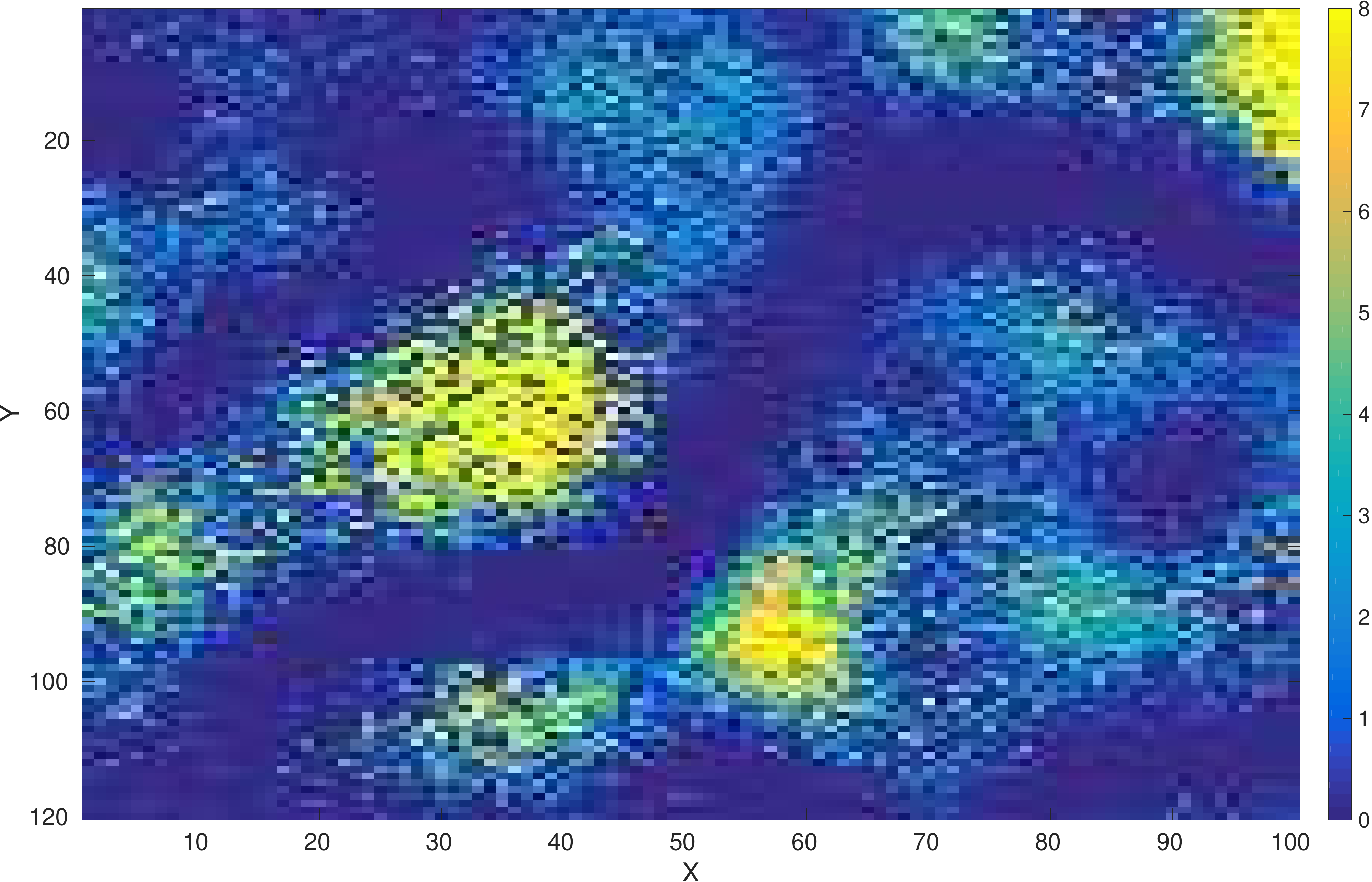}} &
		\subfloat[Mean of the final ensemble with bias-based MEC]{\includegraphics[width=0.45\textwidth]{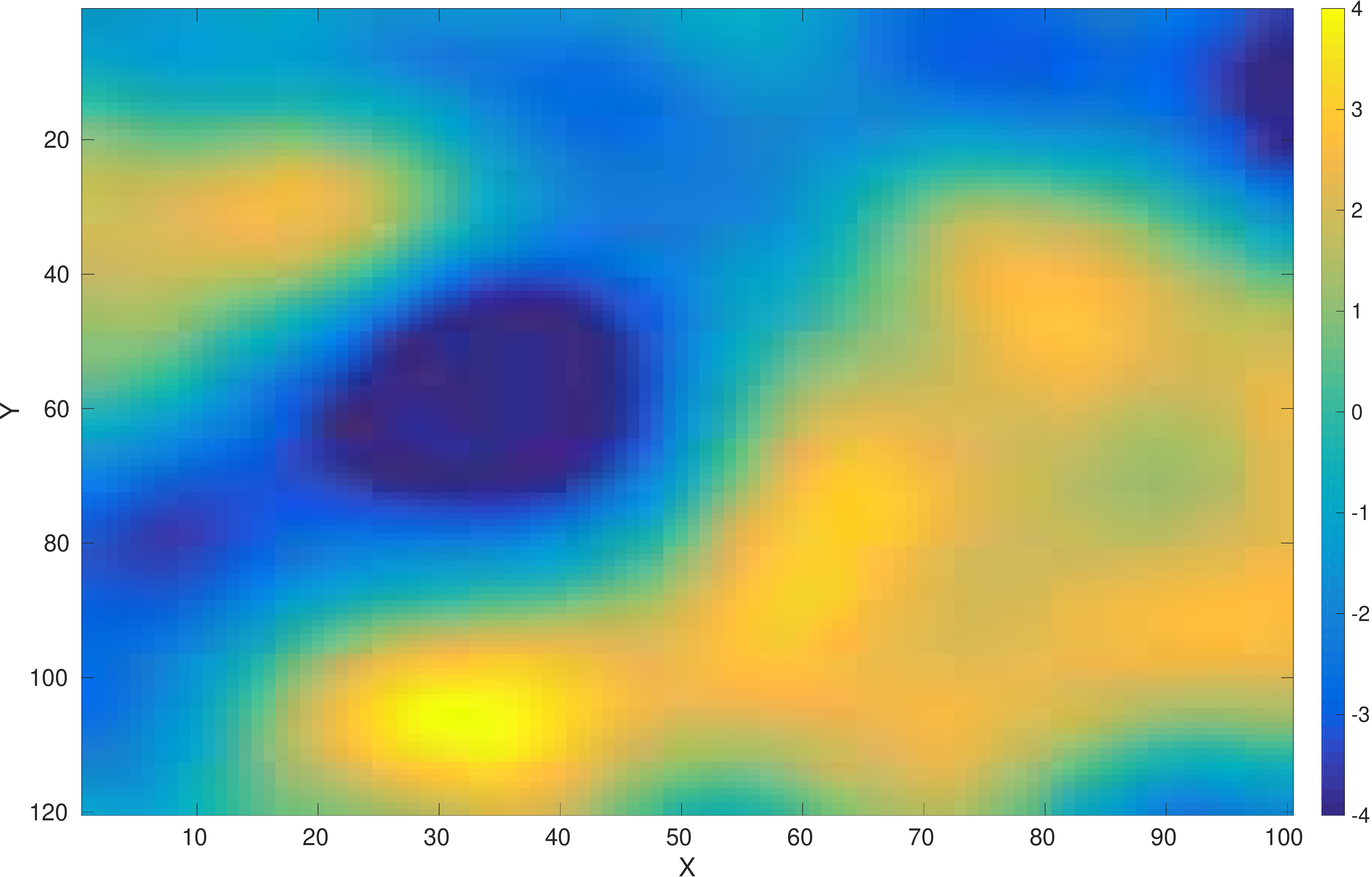}} \\[+0.5cm]
		\subfloat[Box plots of data mismatch with bias-based MEC]{\includegraphics[width=0.45\textwidth]{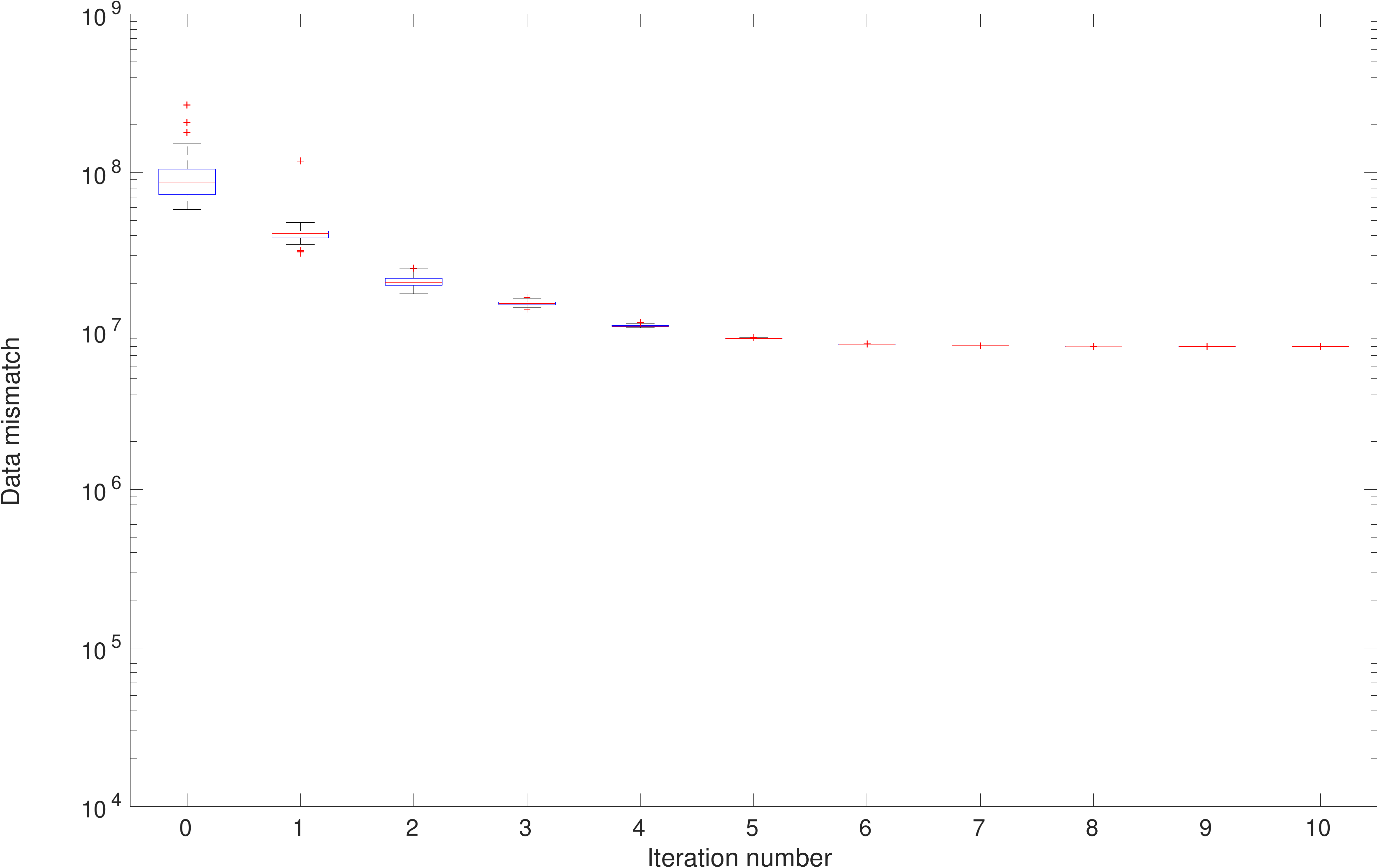}} &
		\subfloat[Box plots of RMSE with bias-based MEC]{\includegraphics[width=0.45\textwidth]{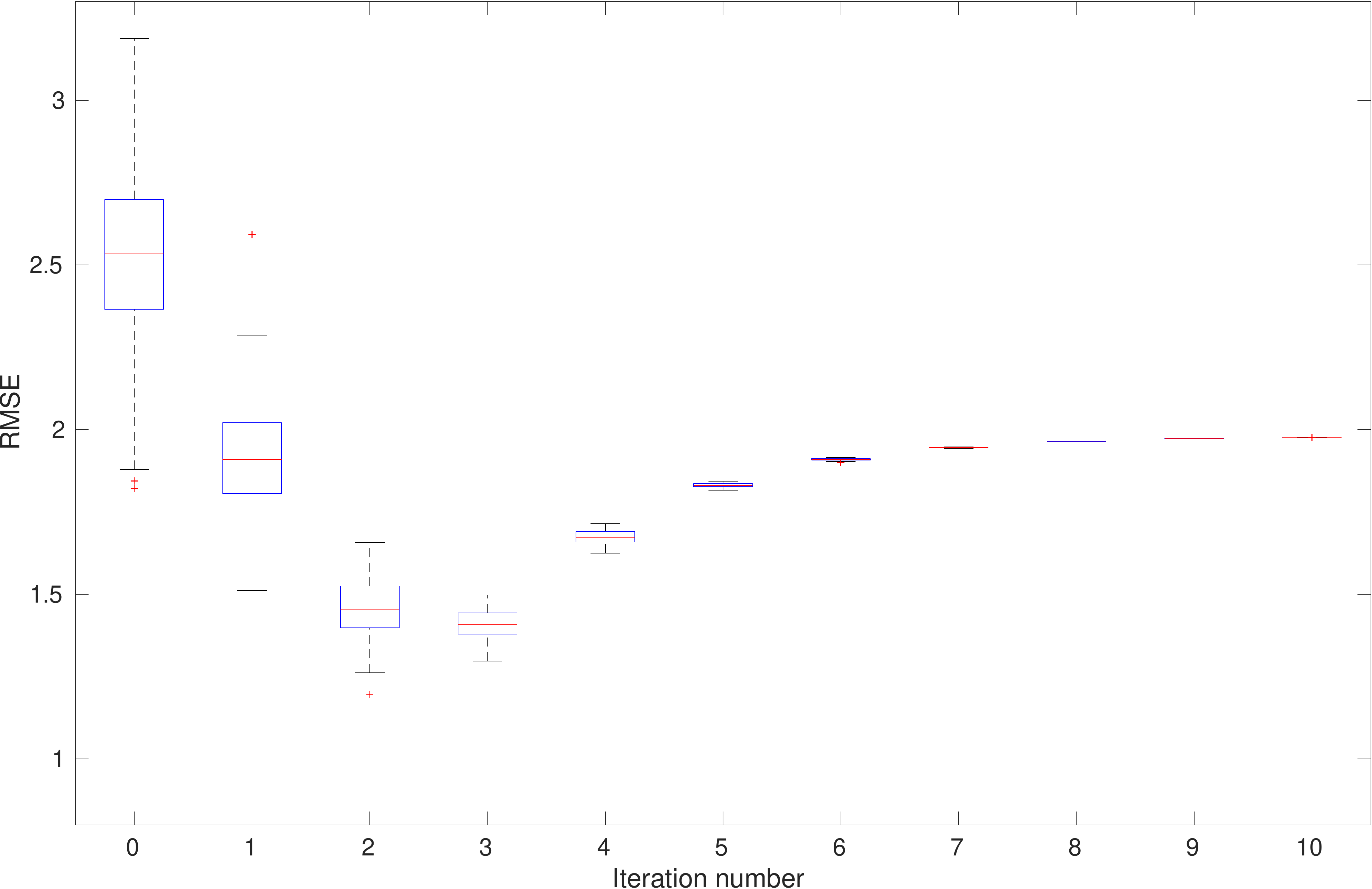}} \\
	\end{tabular}
	\caption{\label{fig:bias-based MEC} Experiment results with bias-based MEC in the \textit{imperfect} scenario.}
\end{figure}  

\renewcommand{\nScale}{0.4}
\begin{figure} 
	\centering
	\begin{tabular}{c}
		\subfloat{\includegraphics[width=0.8\textwidth]{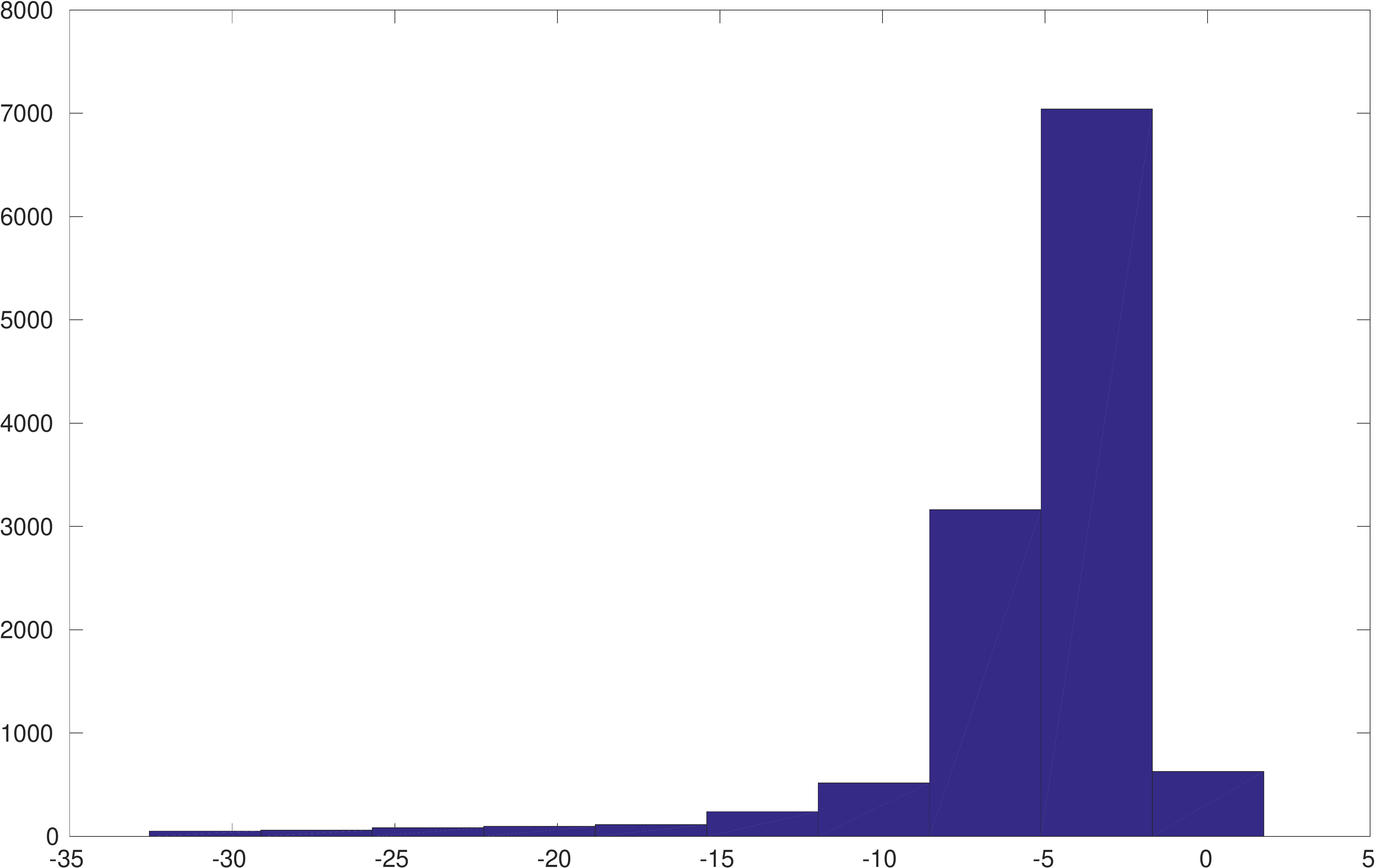}}
	\end{tabular}
	\caption{\label{fig:residual_hist} Histogram of the mean of the residuals with respect to the initial ensemble.}
\end{figure}

An alternative idea for MEC in data assimilation would be that, in Eq. (\ref{eq:residual_rbf_approximation}), instead of adopting kernel-based functional approximation, one may simply approximate the residual term by an unknown bias term, similar to the strategy adopted in, e.g., \cite{dee1995line}. It would then be of interest to see how this alternative MEC method performs, in comparison to kernel-based MEC. For reference later, we call this alternative method bias-based MEC.

In the experiment, we also choose to integrate this bias-based MEC into ensemble-based data assimilation. To initialize an ensemble of biases, we first compute an ensemble of residuals between real observations and simulated observations with respect to the initial ensemble. We then calculate the mean and covariance of the residual, and use these statistics to draw an (initial) ensemble of biases, in a way similar to that we adopted to generate the initial ensemble of model variables. After that, similar to the setting in Eq (\ref{eq:HM_with_SLP}), we augment both model variables and biases, and use the iES to update them in the course of data assimilation.

Figure \ref{fig:bias-based MEC} summarizes the experiment results with respect to bias-based MEC. In comparison to the results with kernel-based MEC in Figures \ref{fig:DA_models_imperfect} -- \ref{fig:dmDiff_Boxplot_imperfect}, it is clear that bias-based MEC tends to result in higher data mismatch and RMSE. In terms of mean $\pm$ STD, the data mismatch values of the final ensemble for bias-based MEC are $(7.9847 \times 10^6) \pm 80.1288$, and the corresponding RMSEs are $1.9767 \pm (2.7584 \times 10^{-4})$. Relative to the mean values, the tiny STDs of the final data mismatch and RMSEs suggest that ensemble collapse is a severe issue in the experiment with bias-based MEC. 

The relative under-performance of bias-based MEC might be partially attributed to the simplifying assumptions, e.g., whiteness, stationarity, and normality \citep{dee1995line}, regarding simulator imperfection. To see this, Figure \ref{fig:residual_hist} shows the histogram of the mean of the residuals with respect to the initial ensemble. As one can see there, the distribution of the mean residuals does not seem to resemble a normal distribution well.

\subsubsection*{The impact of the number $N_{cl}$ of clusters}
\renewcommand{\nScale}{0.4}
\begin{figure} 
	\centering
	\begin{tabular}{c}
		\subfloat{\includegraphics[width=0.8\textwidth]{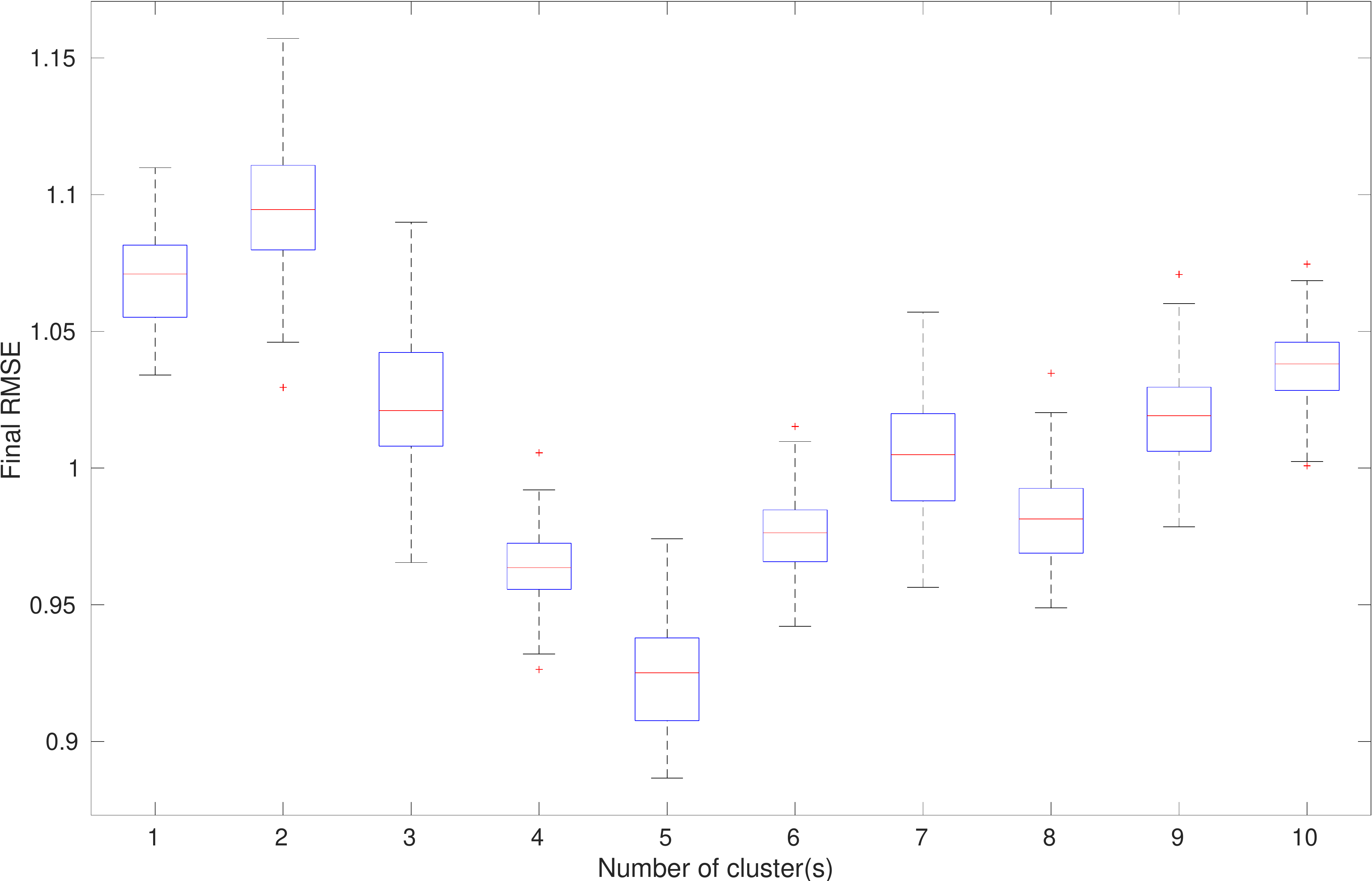}}
	\end{tabular}
	\caption{\label{fig:impact_nCluster} Box plots of RMSEs of the final ensembles, obtained with different numbers $N_{cl}$ of clusters.}
\end{figure} 

As in SLP, when using kernel-based functional approximation for MEC, one can also choose to first group model variables into different clusters, and then estimate an ensemble of kernel parameters for each cluster. The final residual functional is taken as the weighted average of the individual (kernel-based) approximation functional estimated from each cluster, similar to the idea described in Eq (\ref{eq:predicted_residual}). Note that, in this case study, we know that the initial ensemble of model variables is generated using fast Gaussian simulation \citep{lorentzen2017history,Luo2018CorrLoc_Norne}. Therefore, in principle, either the joint or the marginal  distribution of the model variables is unimodal, and intuitively there would be no need to consider an multi-modal-based approximation strategy. Nevertheless, as we will show below, the multi-modal strategy may help improve the performance of data assimilation.

Figure \ref{fig:impact_nCluster} reports the box plots of RMSEs with respect to the final ensembles that are obtained in data assimilation using different $N_{cl}$ values. In the experiment, $N_{cl}$ takes its value from the set $\{1,2,\dotsb,10\}$. As one can see in Figure \ref{fig:impact_nCluster}, except for the case with $N_{cl} = 2$, all other choices tend to result in lower RMSEs, in comparison to the choice of $N_{cl} = 1$. This thus suggests that, similar to the results in SLP (cf. Figure \ref{fig:prediction_ensemlbe_multicluster}), one may obtain better assimilation performance by using a relatively large value for $N_{cl}$ that exceeds the actual number of mode(s) in the distribution of model variables. On the other hand, though, the optimal choice of the value of $N_{cl}$ remains to be an open problem in the current work.    

\section*{Discussion and conclusion}\label{sec:conclusion}

This work focuses on addressing simulator imperfection in data assimilation from a perspective of functional approximation, which leads to an ensemble-based data assimilation framework that integrates functional approximation through a certain machine learning approach into an ensemble-based assimilation algorithm. For better comprehension of how such an integration can be established, we start from considering a class of supervised learning problems, and then discuss the similarity between supervised learning and variational data assimilation. This insight (of similarity) not only leads to an ensemble-based approach to solving supervised learning problems, but also sheds light on the development of an ensemble-based data assimilation framework that, in a natural way, merges machine learning and data assimilation methods to handle simulator imperfection. In the current work, we adopt a kernel-based learning approach to functional approximation. Nevertheless, as discussed in earlier texts, one may also employ other suitable machine learning methods for the purpose of functional approximation.

For performance demonstration, we first study a supervised learning problem. Through the investigations therein, we identify a challenge that may arise when using kernel-based ensemble learning in the presence of multi-modal training inputs. To overcome this problem, we consider a multi-modal learning strategy that helps achieve reasonably good results. Moreover, this multi-modal strategy can be transferred to the data assimilation problem later, also helping improve the performance of data assimilation. Apart from the multi-modal strategy, in the data assimilation problem, we also inspect the performance of the ensemble-based data assimilation framework with the integrated, kernel-based model-error correction (MEC) mechanism. The experiment results indicate that, in this particular case study, using kernel-based MEC tends to improve the data assimilation performance, no matter if simulator imperfection is present or not. In addition, the experiment results also show that kernel-based MEC tends to outperform an alternative, bias-based MEC mechanism.

As a proof-of-concept study, in the current work, we consider a relatively simple data assimilation problem, in which there is only one unknown parameter to estimate for each gridblock. Conceptually, based on Eqs. (\ref{eq:mt_rbf_approximation}) -- (\ref{eq:mt_kernel_para}), it will not be difficult to extend the integrated data assimilation framework to the case studies in which there are multiple unknown parameters on each gridblock (to some extent, this is partially investigated, with the use of Eq. (\ref{eq:residual_rbf_approximation}) in the experiments). Such an extension will be investigated in our future work, with the experiment settings being as close to real field case studies as possible.

\bibliographystyle{ametsoc2014}
\bibliography{./references}

\end{document}